\documentclass[final,11pt,onecolumn,a4]{IEEEtran}

\usepackage{bm,graphicx,epsfig,amsmath,amssymb,mathrsfs}

\usepackage{cite}

\newlength{\intwidth}
\DeclareRobustCommand{\fpint}[2]
  {\mathop{%
     \text{%
       \settowidth{\intwidth}{$\int$}%
       \makebox[0pt][l]{\makebox[\intwidth]{$-$}}%
       $\int_{#1}^{#2}$}}}

\newtheorem{theorem}{Theorem}

\newtheorem{corollary}{Corollary}

\newtheorem{definition}{Definition}[section]

\newtheorem{lemma}{Lemma}

\newtheorem{proposition}{Proposition}

\renewenvironment{proof}[1][Proof]{\flushleft \textit{#1: } }{\ \rule{0.0em}{0.0em}}

\newcommand{\q}{{\bf q}}
\newcommand{\f}{{\bf f}}
\newcommand{\te}{{\bf t}}
\newcommand{\g}{{\bf g}}
\newcommand{\x}{{\bf x}}
\newcommand{\y}{{\bf y}}
\newcommand{\be}{{\bf b}}

\begin{document}



\title{The Hyperanalytic Wavelet Transform}



\author{Imperial College Statistics Section\\
Technical Report TR-06-02\\
\today\\
Sofia~C.~Olhede \& Georgios Metikas\\
Department of Mathematics, Imperial College London, SW7 2AZ UK
\thanks{Manuscript received XXXXXXX XX, XXXX; revised XXXXXXX XX, XXXX.
G. Metikas was supported by an EPSRC grant.}
\thanks{G. Metikas and S. Olhede are with the Department of Mathematics,
Imperial College London, SW7 2AZ, London, UK (s.olhede@imperial.ac.uk). Tel:
+44 (0) 20 7594 8568, Fax: +44 (0) 20 7594 8517.}}

\markboth{STATISTICS SECTION TECHNICAL REPORT: TR-06-02}{Olhede
\& Metikas: Hyperanalytic Wavelet Transform}
\maketitle
\begin{abstract}
In this paper novel 
classes of 2-D vector-valued spatial domain wavelets are defined,
and their properties 
given. The wavelets are 2-D generalizations 
of 1-D analytic wavelets, developed from the Generalized
Cauchy-Riemann equations and represented as quaternionic
functions. 
Higher dimensionality complicates the issue of analyticity, more than one
`analytic' extension of a real function is possible, and an `analytic'
analysis wavelet will not necessarily construct `analytic' decomposition
coefficients. The decomposition of locally unidirectional and/or separable
variation is investigated in detail, and two distinct
families of hyperanalytic wavelet coefficients are introduced, the monogenic
and
the hypercomplex wavelet coefficients.
The recasting of the analysis in a different frame
of reference and its effect on the constructed coefficients is investigated,
important
issues for sampled transform coefficients. The magnitudes of the
coefficients are shown to exhibit stability with respect to shifts in phase.
Hyperanalytic 2-D wavelet coefficients
enable the retrieval of a phase-and-magnitude description of an image in
phase space, similarly to the description of a 1-D signal with the use of
1-D analytic wavelets, especially appropriate for oscillatory signals. 
Existing 2-D directional wavelet decompositions are related to the newly
developed framework, and new classes of mother wavelets
are introduced.
\end{abstract}

\vspace{0.1in}



\section{Introduction \label{Introduction}}
\PARstart{T}{his} paper constructs new classes of 2-D vector-valued mother
wavelet
functions.
The vector-valued wavelet functions are defined to produce 2-D extensions
of 1-D analytic
complex wavelet coefficients, denoted hyperanalytic wavelet coefficients.
The extension of 1-D analytic wavelet coefficients to 2-D is non-trivial.
Firstly more
than one definition of `analytic' exists in 2-D. The suitable choice of the
concept in 2-D
depends on the local directionality and dimensionality of the signal
under analysis, and local unidirectionality, as well as separability, will
explicitly be considered. Secondly, unlike the
1-D case, a distinction can be made between wavelet coefficients obtained from
using
wavelets that are `analytic', wavelet coefficients that are `analytic'
in their spatial index,
and wavelet coefficients
of `analytic' 2-D functions. In 2-D
it will furthermore become necessary to introduce
the notion of $\theta$-`analytic', or `analytic'
in a rotated frame of reference, to interpret the properties of the
wavelet
coefficients. We shall demonstrate that it is possible to choose a definition of `analytic' suitably so that the hyperanalytic wavelet coefficients have
tractable and useful properties in the rotated frame. The
mathematical framework developed for the construction of the hyperanalytic
wavelet coefficients also provides a context to some of the most
promising new complex/vector-valued 
wavelet filter decompositions in the literature \cite{Selesnick2005,Chan},
as well as allows us to construct new classes of vector-valued wavelets.

Over twenty year ago, Grossman and Morlet \cite{Grossman} developed
the Continuous Wavelet Transform (CWT) \cite{mallat}, using
continuous complex-valued mother wavelets. Initial analysis based on wavelet decompositions was implemented
using such mother wavelets. Both magnitude and phase descriptions of non-stationary
signals were determined, and
an early example of analysis include wavelet ridge methods proposed
by Delprat {\em et al.} \cite{Delprat}. 
However subsequently for many years interest focused on the Discrete Wavelet Transform (DWT)
and signal estimation. The DWT was developed to implement the wavelet transform
of time-compact mother wavelets and as compact discrete wavelet filters cannot be exactly analytic \cite{Auscher}, real wavelets were used.
A revival of interest in later years has occurred in both
signal processing
and statistics for the usage of complex wavelets \cite{Selesnick2005,Barber}, and in particular
complex analytic wavelets \cite{kingsbury1,kingsbury2,selesnick1,selesnick2}.
This revival of interest may be linked to the development of
complex-valued discrete wavelet filters \cite{Lina} and the clever dual filter bank \cite{kingsbury1}. The complex wavelet transform has been shown to provide a powerful
tool in signal and image analysis \cite{Selesnick2005}, where most of the
properties of the transform follow from the analyticity of the wavelet function.
This paper strives to build upon the strengths of recent developments, by
deriving large classes of wavelets generalizing
the concept of a 1-D local complex-valued
analytic decomposition to a 2-D vector-valued hyperanalytic decomposition. 

The wavelet decomposition is in 2-D only a single example of a local decomposition.
Local decompositions of an image in terms of spatial
structure
associated with a given scale and spatial position have with great success
formed the basis for many
different procedures of local image
analysis and estimation \cite{AMVA,mallat,starck2}. As the variational structure
in 2-D is much
richer than in 1-D, the term `local' may denote a variety of different forms
of spatially
limited structures, and so examples of decompositions include {\em i.e.}
wavelets
\cite{mallat}, curvelets \cite{starck2002} and
bandelets \cite{pennec}  {\em etc}. 
A local decomposition method
should be chosen
with easy interpretability of the decomposition coefficients and compression
in mind, as both of these characteristics
facilitate the analysis of observed images. 
For local image analysis, the local directionality and dimensionality
of variation 
are key features of the local visual appearance of the image. These two features
need to be
characterised,
and can be used 
for further processing,
such as image segmentation, feature detection and disparity estimation,
see for example work by von B{\"u}low and Sommer \cite{BulowSommer1}, as
well as that of von B{\"u}low
\cite{BulowPhd}, Hwang {\em et al.} \cite{Hwang2001}, Antoine
{\em et al.} \cite{Antoine1,Antoine3}
{\em etc.} To represent both local dimensionality and directionality at each
spatial point more
than
a scalar-valued
representation is required.  Based on a fixed local dimensionality the characterisation
of the image will of course vary.
For example, images that are locally 1-D, are well characterised by the direction
of the image variation, its magnitude and the scale of variation in the given
direction.
Texture features present in an image
may correspond to locally separable structures, characterised by the
direction
of the separable structure, the local magnitude of variation and the scale
of variation in two directions. Single texture components often correspond
to variations associated with the same scale in both perpendicular directions.

The observational axis of an image
cannot be assumed to coincide with the natural axis of the local structure.
It is therefore necessary that the local representations should
not substantially change form with rotations
and translations of the spatial index, but exhibit equivariance
to such operations. 
Of course the CWT exhibits equivariance
when calculated in continuous space at all values of the locality index
\cite{AMVA},
but
if the transformation is to be discretized, it needs to
satisfy some further constraints so that the representation is stable.
In particular we require that
the magnitude of the vector-valued decomposition should be stable to shifts
in phase of the signal, or
a sampled collection of coefficients of the process may exhibit high translation variance in comparison to a sampled collection of coefficients obtained from
the same image subjected to a small spatial shift in both axes simultaneously
\cite{Selesnick2005,gopinath2003}. Translation variance is clearly an undesirable
feature of a decomposition.
The simple
characterisation
of any shift in index is also important for analysis of multiple images.
Vector-valued representations, as will be demonstrated, can with a suitable
choice
of decomposition filter, both provide decompositions stable to small spatial
shifts (or phase shifts)
and allow for the easy characterisation of rotations and translations of
the spatial index.

As examples of successful existing
vector-valued local representations of image structure we note the 2-D
analytic
signal \cite{Hahnbook},
or the complex CWT \cite{Selesnick2005,gopinath2003}.
The 1-D analytic signal at a given value of $x_1$ provides a meaningful magnitude
and phase
description of a signal, if only a single
oscillatory
component is present.
For each
fixed scale the complex wavelet coefficients correspond to spatially local representations
of the image in terms of two wavelet coefficients, but the magnitude
and phase of the complex
wavelet coefficients are not when an arbitrary complex-valued wavelet is
used either meaningful
nor  are they interpretable. For this reason complex wavelet transform coefficients
obeying suitable conditions are generally used \cite{Selesnick2005,Fernandes2}.
We aim to generalise the usage of spatially local two-vector representations
to local three- and four-vector representations, using hyperanalytic wavelet
decompositions.

Vector-valued representations of local structure in terms of quaternionic
representations, or equivalently three- and four-vectors,
have been developed by 
von B{\"u}low and Sommer \cite{BulowSommer1}, Felsberg and Sommer \cite{FelsbergSommer}
as well as Hahn and Snopek
\cite{HahnSnopek2}. In particular Hahn \cite{Hahn}, von B{\"u}low and Sommer
\cite{BulowSommer1},
as well as Felsberg and Sommer \cite{FelsbergSommer}, have developed 2-D
extensions of the 1-D analytic signal.
These representations give
descriptions of local variation in terms of a local magnitude
and phase/s.
Phase should represent local variations
and structure, whilst
magnitude is usually interpreted as local signal presence
\cite{Wakin}. A local well-constructed magnitude and phase description of
a signal is often
very
informative. However if more than a single narrowband component is present
at a given
spatial location in an image, the 2-D `analytic' signal, no matter which
of the above generalizations are chosen, will completely lack interpretation,
just like in the case of the 1-D analytic signal of a multi-component signal
\cite{Boashash2,Cohen}.

The work by Hahn and Snopek \cite{HahnSnopek2} consists of
forming a suitable
2-D quaternionic (four-vector valued)
extension of the Wigner-Ville distribution 
\cite{Cohen}. In 1-D the Wigner-Ville distribution is known to be uninformative when describing
signals contaminated by noise \cite{Flandrin2}, and is unsuitable as a description of signals consisting of more than a
single component
at any spatial position
due to interference between components \cite{Cohen}.
To deal both with the addition of noise and the presence of multiple components in 1-D, 
local analytic representations
have been constructed, such as analytic wavelet coefficients 
\cite{Holschneider,Olhede,Delprat},
or other local representations of the signal calculated in conjunction with a localisation
procedure \cite{Olhede3,icohen}. In this paper we aim to construct 2-D `analytic wavelet coefficients'
or, locally `analytic' representations of images, by using vector-valued
wavelet functions. 
The wavelets
will be constructed to merge the 2-D
wavelet
transform with existing 2-D `analytic' representations.
The synthesis of local image representations and magnitude
and phase descriptions of local structure will be shown to enjoy many special
properties, but requires several non-trivial generalizations of existing
1-D concepts.

Before proceeding to the construction of 2-D wavelets to serve the same purpose
as 1-D analytic wavelets, let us first briefly revisit
the properties of the 1-D analytic wavelet decomposition. The decomposition
is constructed from a single mother wavelet function 
$\psi(x_1)\in L^2\left({\mathbb{R}}\right),$
that satisfies the admissibility condition \cite{AMVA}. The mother wavelet
is translated
in time by $b,$ and scaled by factor $a$ to form a family of wavelets denoted
by 
$\left\{\psi_{a,b}(x_1)
\right\}.$ For $\psi(\x)\in L^2\left({\mathbb{R}}^d\right),$
the dilation and scaling operators for any index dimension
$d$ are 
${\cal D}_a\psi\left(\x\right)=\left|a\right|^{-d/2} \psi\left(\x/a\right),$
and $
{\cal T}_{\be}\psi\left(\x\right)=\psi\left(\x-\be\right),$ respectively.
A member of the 1-D wavelet family is $\psi_{a,b}(x_1)={\cal
D}_a{\cal T}_{b}\psi\left(x_1\right).$ 
To represent
a 1-D signal $g(x)$ in terms of local contributions the CWT, and its inversion,
are defined by
\begin{equation}
w_{\psi}(a,b;g)=\int_{-\infty}^{\infty} g(x_1) \psi_{a,b}^{*}(x_1)\;dx_1,
\quad g(x_1)=\frac{1}{c_{\psi}}\int_{-\infty}^{\infty} \int_{-\infty}^{\infty}
w_{\psi}(a,b;g)\psi_{a,b}(x_1)\;\frac{da}{a^2}db.
\end{equation}
$c_{\psi}\in{\mathbb{R}}$ is a constant specific to the function $\psi(x_1),$
its existence guaranteed by the admissibility condition.
If $\psi(\cdot)$ can be interpreted as local in time to $x_1=0$ and frequency
to $f_1=f_0$, the signal is reconstructed in terms of the local contributions
$w_{\psi}(a,b;g):$ furthermore $w_{\psi}(a,b;g)$ is associated with a function local to position $x_1=b$ and
frequency $f_1=f_0/a.$
The wavelet coefficients are used to characterise the
local properties of the signal and the decomposition is often described as
a `local Fourier'
transform, even if this is a somewhat inaccurate description.

If $\psi(x_1)$ is real-valued, then the magnitude of $w_{\psi}(a,b;g)$
can be used to determine local signal content, or be used as the basis of
more
complicated analysis methods. If $w_{\psi}(a,b;g)$ is complex-valued then it can be represented by its
magnitude and phase.  If $\psi(x_1)$ is analytic, denoting an analytic
mother wavelet \footnote{An `analytic' mother wavelet in this context is
a complex-valued mother
wavelet that is also an analytic signal \cite{Cohen}.}
by $\psi^{(+)}(x_1),$
then the wavelet is a complex-valued
function. The CWT of a real-valued signal using a complex-valued analytic
wavelet yields
complex-valued wavelet coefficients, where the coefficients are an analytic signal in $b.$ 
In 1-D, if we define the complex variable $t_1=x_1+\bm{j}y_1$ and restrict
$y_1>0,$
then
an analytic signal $g^{(+)}(x_1)=g(x_1)+\bm{j}g^{(1)}(x_1),$ evaluated in the complex argument, that is $g^{(+)}(t_1),$
is an analytic function in the upper half-plane of $y_1>0$ \cite{Hahnbook}.
$g^{(+)}(t_1)$ thus satisfies the Cauchy-Riemann equations in the upper half-plane and this determines
the relationship between $g(t_1)$ and $g^{(1)}(t_1).$ 
The fact that a real and imaginary part of $g^{(+)}(x_1+\bm{j}y_1)$
satisfy the Cauchy-Riemann equations imply that the phase of $g^{(+)}(x_1)$ will be well defined, once some additional restrictions are placed
on $g(x_1).$ Given the analytic wavelet coefficients are an analytic signal,
and local in frequency,
the instantaneous
frequency
of an oscillatory component present at scale $a,$ may for example be determined
using its
analytic wavelet coefficients, and multi-component
oscillatory
signals may be characterised via wavelet ridge analysis \cite{Delprat}. 
Note that if $w_{\psi}(a,b;g)\in{\mathbb{C}}$ but
$\psi(x_1) \neq \psi^{(+)}(x_1),$ then $w_{\psi}(a,b;g)$ can still be represented
in terms of a magnitude and phase, but the phase will not in general be interpretable
in terms of local structure.
Furthermore wavelet
coefficients constructed from analytic mother wavelets have a magnitude invariant
to phase-shifts
and this alleviates observed problems  with coefficients
exhibiting strong shift-dependence \cite{Selesnick2005,Olhede5}. The analytic
wavelet
has
one major drawback, longer essential time-support, but provides the framework
for analysis shorter time-support filters may approximately provide,
see for example work by Selesnick {\em et al.} \cite{Selesnick2005} as well
as Fernandes {\em et al.} \cite{Fernandes2}. Selesnick {\em et al.} give a comprehensive
discussion of the many important properties of analytic wavelet coefficients
when calculated in both 1-D and higher dimensions, and these properties justify
our interest in developing a single mathematical framework for the interpretation
of 2-D locally `analytic' representations.

In 2-D the CWT of $g(\x)\in L^2(\mathbb{R}^2)$
is implemented by choosing an admissible mother \cite{AMVA} wavelet $\psi(\x)\in
L^2(\mathbb{R}^2),$ and from this function, a family of functions
$\left\{\psi_{\bm{\xi}}(\x)\right\}$ for $\bm{\xi}=\left[a,\theta,\be\right]^T$
is defined by:
\begin{equation}
\psi_{\bm{\xi}}(\x)={\cal{D}}_a {R}_{\theta}{\cal{T}}_{\be}\left\{\psi\right\}
(\x)=U_{\bm{\xi}}\psi(\x)
\label{famofwav}.
\end{equation}
$R_{\theta}g(\x)=g(\bm{r}_{-\theta}x),$ is the rotation operator given in
matrix notation as:
\[\bm{r}_{\theta}=\left[ \left\{ \cos(\theta), -\sin(\theta) \right\},\left\{\sin(\theta),
\cos(\theta)
\right\} \right].\] 
The purpose of including $\theta$ in the decomposition is to identify local
behaviour like unidirectional, or separable structures, not aligned with
the observational axes by rotating the
analysis function over the full set of possible orientation, when necessary
\cite{AMVA}.
The wavelet coefficients of function $g(\x)$ are then defined by \cite{AMVA}:
\begin{eqnarray}
\nonumber
w_{\psi}(\bm{\xi},g)&=&\langle \psi_{\bm{\xi}}(\x),g(\x)\rangle=
\int_{-\infty}^{\infty}\int_{-\infty}^{\infty}g(\x)\psi_{\bm{\xi}}^*(\x)\;d^2\x,
\label{wavtrans}
\end{eqnarray}
where $*$ denotes the act of conjugation. 
Just like in 1-D the $w_{\psi}(\bm{\xi},g)$
coefficients can be used
to reconstruct the function $g(\x)$ from a weighted average of the 
$\left\{\psi_{\bm{\xi}}(\x)\right\}$ by:
\begin{eqnarray}
g(\x)=\frac{1}{c_{\psi}}\int_{-\infty}^{\infty} \int_{-\infty}^{\infty} 
\int_{0}^{2\pi} \int_{0}^{\infty} w_{\psi}(\bm{\xi},g)\psi_{\bm{\xi}}(\x)\;\frac{da
d\theta d^2\be}{a^3}.
\label{reconny5}
\end{eqnarray}
To form $w_{\psi}(\bm{\xi},g),$ $g(\x)$ has been localised
in both space and spatial frequency depending on the form of $\psi(\x).$
$\psi(\x),$ the mother wavelet function, thus needs to be chosen with care,
so that $w_{\psi}(\bm{\xi},g)$ has suitable properties. We wish to derive
$w_{\psi}(\bm{\xi},g)$ that are locally `analytic'.

Firstly we need to specify what we mean by `analytic' in 2-D. A reasonable starting
point to generalise $g^{(+)}(x_1)$ to 2-D is by considering 2-D versions
of analytic
functions and then forming the appropriate limit, as $\lim_{y_1\rightarrow
0^+} g^{(+)}(t_1)=g^{(+)}(x_1).$
The definition of a 2-D analytic
function has previously been examined in detail in pure mathematics, within
the field of Harmonic Analysis, by Stein and Weiss, see for example \cite{Stein,SteinWeiss}.
Similar concepts have appeared in applied mathematics and geophysics
\cite{dixon,nabighian}, and in signal processing 
\cite{Hahn,FelsbergSommer,BulowSommer1}. An important fact to note is that
the generalisation of analyticity
to higher dimensions is not unique \cite{SteinWeiss}, and that more than
a
single generalisation of the Cauchy-Riemann system, and thus to the definition
of a hyperanalytic function, is possible. 
Any function that satisfies a generalization of the Cauchy-Riemann system
will in fact be referred to as a {\em hyperanalytic} function in this paper.
A hyperanalytic
signal is, in contrast to the hyperanalytic function, defined as the suitable limit of such a
function.

Two possible hyperanalytic signals are relevant in this context: the {\em hypercomplex
signal} of von B{\"u}low and Sommer \cite{BulowSommer1}, and the
{\em monogenic signal} of Felsberg and Sommer \cite{FelsbergSommer}.
For clarity we will
provide some discussion of the two hyperanalytic signals
that are both possible 
and useful 2-D extensions of the 1-D analytic signal,
corresponding in turn to
four- and three-vector
valued functions respectively. 
The hyperanalytic signal is usually represented
as a quaternionic function, see for example \cite{FelsbergSommer,BulowSommer1}, {\em i.e.}
as an object
taking values in the skew field of quaternions, denoted ${\mathbb{H}}.$ The
advantage
of using
the quaternion representation is that convenient polar forms (magnitude and phase/s) are naturally introduced \cite{FelsbergSommer,BulowSommer1}, representing
magnitude and variation of either separable or unidirectional variation.
Once the choice of `analytic' has been specified, or convenient local polar
representation, the `analytic' mother wavelet is formed.
We approach the choice
of mother wavelet by striving to separate disparate structures
present at the same spatial position and to represent
the isolated local structure, in terms of local structure (phase) and magnitude.
The interpretation of the phase functions of these coefficients must be established,
and this will depend on the form 
of $\psi(\x),$ and assumptions regarding the local structure of $g(\x).$
We
shall construct quaternionic wavelets starting from a real 2-D wavelet that
will
be augmented into a quaternionic object, in analogue to 1-D analytic wavelets.
The construction will
rely on the choice of the starting point of a real mother wavelet, and the
chosen form of hyperanalytic extension of this real function.

One starting point to the construction of the quaternionic mother wavelet
is the anisotropic, or elongated,
real-valued mother wavelet. 
Such a wavelet can construct coefficients representing plane wave or separable
structures well. The wavelet would be extended into a quaternionic
function so that a hyperanalytic wavelet coefficient 
representation can be formed of unidirectional or separable local structure in terms of the hyperanalytic
phase/s and magnitude.
Using quaternionic
coefficients enables the easier parameterisation of local variation, and
will rely on the transform coefficients being hyperanalytic, rather
than the wavelet functions themselves. 
The strategy is related to work by Selesnick {\em et al.} \cite{Selesnick2005} (who use directional
complex wavelets), as well as Chan {\em et al.} \cite{Chan} (who use separable quaternionic
wavelets).
Both the construction by Selesnick {\em et al.} and Chan {\em et al.}
are implicitly based on the hypercomplex signal, and we explicitly construct
a
formal framework for interpreting such representations. We also form the
monogenic
extension of real directional wavelets.

Another possibility is to separate the signal into components
present at given scales using an isotropic real wavelet, and then {\em assuming} locally at that spatial
position
and scale either plane waves, or
potentially separable structures only are present, represent these by forming
the hyperanalytic signal of the wavelet coefficients. This corresponds to the 2-D version of the wavelet projection
representations proposed by Olhede and Walden \cite{Olhede3}.

One of the main questions to tackle in the construction is the appropriate
interpretation that can be given
to the
quaternionic
wavelet
coefficients depending on the choice of mother wavelet function.
The answer is complicated in 2-D by the 
fact that rotations do {\em not} commute with the operators that construct
hyperanalytic
images. 
Given the complication introduced by the rotations, we start by considering
analysis when no such operation is required:
we assume that the image is aligned with the observational
axes, and set the angle of rotation in the decomposition equal to zero. We
then construct wavelet functions such that the CWT
coefficients are hyperanalytic signals in their spatial indexing, and refer
to such mother
wavelets as {\it hyperanalyticizing}. The mother wavelet
function that constructs monogenic coefficients, is a monogenic signal, but we refer to any quaternionic analysis filter that
produces hypercomplex wavelet coefficients as a {\em hypercomplexing
wavelet}, thus explicitly noting that the wavelet produces hypercomplex
coefficients but is not necessarily itself a hypercomplex signals. This then,
when no rotation is implemented, produces hypercomplex wavelet coefficients
from a real-valued image. 

To deal with
the act of rotation
we introduce the definition of a $\theta$-hyperanalytic signal,
that is a function if observed in some rotated frame
of reference, is a hyperanalytic signal. 
Once the monogenic signal has been extended to $\theta$-monogenic, it follows
that the quaternionic wavelet producing the $\theta$-monogenic wavelet coefficients
is itself $\theta$-monogenic. In contrast, the hypercomplexing wavelet
is not a $\theta$-hyperanalytic signal as the property does not even hold
for $\theta=0.$ We discuss in detail the interpretation of
the wavelet coefficients constructed from the family of wavelet functions
constructed by scaling, rotating and translating the hyperanalyticizing mother
wavelet. In analogue to the 1-D analytic and anti-analytic
decomposition of a real-valued signal, see \cite{Papoulis}, we introduce the decomposition of a
real-valued image into $\theta$-hyperanalytic
and anti-hyperanalytic components. The hyperanalyticizing wavelet annihilates/converts
the anti-components
to construct an object equal to the real CWT of a $\theta$-hyperanalytic function.
This equivalence is important not only for the interpretation of the constructed
coefficient but also allows for the derivation of additional useful properties of
the hyperanalytic CWT.

Specific features of local structure may be extracted
directly and easily from the defined wavelet coefficients, for instance potential local unidirectionality and separability in the rotated frame of reference
may be characterised. For aggregations of oscillatory signals, the wavelets we develop
in this paper, can be used to separate the multiple components, and
simultaneously give simple descriptions of the components. 
Previous work on local analysis of oscillatory signals is clearly distinct,
a number of authors, see for example
\cite{Hwang2001,Hwang1998,GonnetTorresani,Antoine3}, have discussed wavelet
ridge analysis using complex wavelets. The
advantage of using quaternionic rather than complex wavelets, as suggested
in this paper, is that separable
oscillations admit a simple and compact representation. Furthermore for local
plane
waves using a suitable choice of quaternionic wavelets, it is not necessary
to calculate the transform at all directions to determine the local plane
wave description. Thus the usage of quaternionic wavelet functions either introduces
the simplified description of additional texture structure, or simplifies
the analysis procedure.

Coefficient equi- and invariance to local changes of phase and
rotation alignment between the signal and the wavelet are also discussed.
For the CWT this is not of immediate importance, as the transform coefficients exhibit equivariance
with respect to changes in axes
\cite{AMVA}, but as discrete space decompositions
may be related to continuous space decompositions, the stability of the decomposition
to small spatial misalignment (phase-shifts)
is very important \cite{gopinath2003}. 
For the hypercomplex and monogenic wavelet coefficients, we show that
with a suitable choice of wavelet, magnitude invariance to phase shifts
is achieved. 
For the monogenic wavelet coefficients, mother wavelet functions may be selected,
such that the coefficients also exhibit invariance of magnitude with respect
to local rotations. 

We discuss generic structures of mother wavelet constructions, i.e. combining separable/ isotropic/ directional scale localising procedures, with both the hypercomplex and monogenic
signal representation, but also give specific families of mother wavelets
falling into these classes. 
Existing quaternionic
filters can be related to the hypercomplex and monogenic constructions 
\cite{Metikas,Chan,Chan2,Chan3,FelsbergSommer2}.
In analogue to 1-D wavelets approximating a local Fourier transform, most
of the existing
quaternionic
wavelets have been constructed to serve the role of a local Quaternion
Fourier Transform (QFT).
In addition to \cite{Chan}, work in a similar vein is due to 
\cite{Bayro,Buelow1998,Traversoni,Traversoni2} and corresponds to constructing
local versions of the QFT. He \& Yu \cite{He,He2} constructed quaternionic
multi-variate
decompositions,
and also work by Hsieh \cite{Hsieh,Hsieh2} on motion
smoothing should be noted, as well as work by Felsberg and Sommer on the monogenic scale space \cite{FelsbergSommer2}. 
We provide a single complete framework that relates all
of these constructions, and link the existing decompositions to Harmonic
Analysis. We also construct completely
new families of quaternionic wavelets that belong to the same family, and
determine their properties.

Central to the developments in this paper is the concept
of hyperanalyticity and local versions thereof. The hyperanalytic signal is developed
as a generalization of the 1-D analytic signal. The analytic signal is best
introduced via the 1-D
Fourier
transform, 1-D local oscillations and analytic functions,
and so we start by a discussion of these aspects of the analytic signal.

\section{Oscillations, dimensionality \& analyticity}
\subsection{Analyticity}
An oscillatory 1-D signal $c(x_1)$ is modelled as
an amplitude and frequency modulated signal \cite{Boashash2}:  
\[
c(x_1)=a_c(x_1)\cos(2\pi \varphi_c(x_1))
=\Re\left\{a_c(x_1)\exp(2\pi \bm{j}\varphi_c(x_1))\right\}.\] 
To determine the frequencies present in signal $c(x_1)$
the Fourier Transform (FT) of the signal could be calculated.
The FT of a generic $d$-D signal $g(\x),$ denoted $G(\f),$ is defined
by:
\begin{equation}
G(\f)=\int_{\mathbb{R}^d} 
g(\x) \exp(2\bm{j}\pi \f^{T} \x)\;d^d\x,\quad 
G(\f)=\left|G(\f)\right|\exp(-2\bm{j}\pi\varphi_G(\f)),
\end{equation}
and $g(\x)$ can be reconstructed via
\begin{equation}
g(\x)=\int_{\mathbb{R}^d}  \left|G(\f)\right|\exp(2\bm{j}\pi (\f^{T} \x-\varphi_G(\f)))\;d^d\f.
\end{equation}
In 1-D at $x=x_1$ if $g(x_1)=c(x_1)$ is oscillatory, and has a clearly defined content at frequency $f_1=f_0(x_1)$ approximately
we then find
$
c(x_1)=\left|C(f_0(x_1))\right|\cos(2\pi (f_0(x_1) x_1-\varphi_C(f_0(x_1)))).
$
To calculate the amplitude and phase of a real-valued signal,
or $\left|C(f_0(x_1))\right|$ and $2\pi (f_0(x_1) x_1-\varphi_C(f_0)),$
from $c(x_1),$
the analytic extension of a signal could be used \cite{Hahnbook}. This is
defined for an arbitrary real signal $g(x_1)$ by:
\begin{equation}
\label{analyticsignal}
g^{(+)}(x_1)=2 \int_{0}^{\infty} G(f) \exp(2\pi\bm{j}f x_1)\;df=
g(x_1)+\bm{j}g^{(1)}(x_1).
\end{equation}
$g^{(1)}(x_1)$ is the {\em Hilbert Transform} (HT) of $g(x_1)$. The HT in the time domain and the analytic signal in the frequency domain respectively, are
\begin{equation}
\label{HT}
g^{(1)}(x_1)={\cal H}\left\{g\right\}(x_1)=\frac{1}{\pi}\fpint{-\infty}{\infty}
\frac{g(y)}{x_1-y}\;dy,\quad
G^{(+)}(f_1)=G\left(f_1\right)\left(1+{\mathrm{sgn}}(f_1)\right).
\end{equation}
The HT is as usual defined as a Cauchy principal value integral \cite{Hahnbook}.
Using a stationary phase approximation \cite{Boashash2} the analytic signal of $c(x_1)$ is then approximately
$c^{(+)}(x_1)=c(x_1)+\bm{j}c^{(1)}(x_1)
=\left|C(f_0(x_1))\right|\exp(2\bm{j}\pi (f_0(x_1) x_1-\varphi_C(f_0(x_1))))+o(1)$. Clearly the
amplitude and phase of $c(x_1)$ are retrieved by $a_c(x_1)=\left|c^{(+)}(x_1)\right|,$
and $\varphi_c(x_1)=\frac{1}{2\pi}\tan^{-1}\left(
c^{(1)}(x_1)/c(x_1)\right).$ For a generic signal $g(x_1)$ the
local amplitude and phase may be defined by:
\begin{equation}
a_g(x_1)=\left|g^{(+)}(x_1)\right|,\quad \varphi_g(x_1)=\frac{1}{2\pi}\tan^{-1}\left(\frac{
g^{(1)}(x_1)}{g(x_1)}
\right).
\label{analyticpolar}
\end{equation}
For $a_g(x_1)$ and $\varphi_g(x_1)$
at point $x_1$ to be meaningful representations of $g(x_1),$ the latter signal
is assumed to be at $x_1$ mainly limited in frequency to oscillations
with a single period
The key ingredient in defining a local amplitude and phase via equation (\ref{analyticpolar})
is the calculation of the analytic signal and the definition of
$g^{(1)}(x_1).$

The 1-D analytic signal is the limit of an analytic function, where a {\em function} is denoted analytic when it satisfies the Cauchy-Riemann
equations in the upper half of the complex plane. To construct an analytic
function in a complex argument from a real function in a real argument,
the Poisson
kernels and Poisson convolutions are introduced
\cite{SteinWeiss}.
The Poisson kernels, and the convolutions of $g(\cdot)$ with these kernels, are using an auxiliary variable $y$ defined by:
\begin{equation}
\label{1Dpoisson}
p_H(t_1)=\frac{1}{\pi}\frac{y}{x_{1}^2+y^2},\quad q_H(t_1)=
\frac{1}{\pi}\frac{x_1}{x_1^2+y^2},\quad u_g(t_1)=g(\cdot)*p_H(\cdot,y),\quad
v_g^{(1)}(t_1)=g(\cdot)*q_{H}(\cdot,y).
\end{equation} 
Define complex-valued variable $t_1=x_1+\bm{j}
y,$ and note
as $y\rightarrow 0^{+},$ $u_g(t_1)\rightarrow g(x_1)$ and $v_g^{(1)}(t_1)
\rightarrow 
{\mathcal{H}}\left\{g\right\}(x_1)$, see \cite{SteinWeiss} for details.
$k_g^{(\pm)}(t_1)=u_g(t_1)\pm\bm{j}v_g^{(1)}(t_1)$ satisfies the
Cauchy-Riemann
equations in $\pm y_1>0$ \cite{SteinWeiss}, and are denoted analytic (+) and anti-analytic
(-) functions, respectively. As $y \rightarrow 0^{\pm}$
$k_g^{(\pm)}(t_1)\rightarrow g^{(\pm)}(x_1)=g(x_1) \pm \bm{j} {\cal H}g(x_1),$
and the limit of an analytic function is an analytic signal.
Note that $g^{(+)}(t_1),$ {\em i.e.} the analytic signal
evaluated at the complex argument, is still an analytic function for $y_1>0,$
see also \cite{Hahnbook}[p.~5].
More importantly the analytic and anti-analytic signals are represented in polar form by
\begin{equation}
g^{(\pm)}(x_1)=\left|g^{(\pm)}(x_1)\right|e^{\pm 2\pi \bm{j}\phi_g(x_1)}.
\label{polaranalytic}
\end{equation}
Any real signal $g(x_1)\in L^2({\mathbb{R}})$ can be decomposed into an analytic and an anti-analytic signal \cite{Papoulis}:
\begin{eqnarray}
\nonumber
g(x_1)&=&\frac{1}{2}
\left[g^{(+)}(x_1)+g^{(-)}(x_1) \right] =
 \frac{1}{2}\left[
g(x_1)+\bm{j}g^{(1)}(x_1)+
g(x_1)-\bm{j} g^{(1)}(x_1)\right]\\
&=&\left|a_g(x_1)\right|\cos(2\pi \phi_g(x_1)),
\label{realphase}
\end{eqnarray}
and so a local magnitude $a_g(x_1)$ and phase $\phi_g(x_1)$ can for any real-valued
function $g(x_1)$ be defined
from
the modulus and phase of $g^{(\pm)}(x_1),$ whose forms are given by equation
(\ref{polaranalytic}). We interpret equation (\ref{realphase}) as the analytic and anti-analytic decompositions of $g(x_1),$ and this decomposition is important
to derive the properties of the 1-D analytic wavelet transform, as demonstrated
in \cite{Olhede1}[p.~426]. These well-known properties of the analytic signal
need to be suitably extended to the multi-dimensional choice of `analytic'.

\subsection{Hyperanalyticity}
In higher dimensions to determine equivalents of $a_g(x_1)$ and $\phi_g(x_1)$
from an observed
real-valued function $g(\x)$ the equivalent/s of function $g^{(1)}(x_1)$
must be defined and calculated. Based on these definitions the equivalents
of $a_g(x_1)$ and $\phi_g(x_1)$ will then be given some suitable interpretation. For this purpose d-D extensions of the Cauchy-Riemann
equations are used, and any given extension corresponding to a set of equations
will be denoted a set of {\em Generalized Cauchy-Riemann equations}. In 1-D the spatial variable is 1-D, and there is
a single associated auxiliary variable $y_1,$ where an analytic function satisfies
the Cauchy-Riemann equations in $y_1>0.$ We consider, in {2-D} and higher dimensions,
augmenting the spatial variable $\x$ by a set
of auxiliary variables collected in vector-valued variable $\y,$ of dimension $p$ (some restrictions apply to
the choice of $\y,$ and its dimensionality), where $p$ will be related to
the local structure of the image to be represented. 
We define $\Gamma^{(p)}=\left\{\y,\;y_i>0,\;i=1,\dots,p\right\},$
and consider hyperanalyticity for  $\y\in \Gamma^{(p)},$ and $\x\in{\mathbb{R}}^d.$
Examples of such spaces include $p=1$ with $y_1>0$ and $\x\in{\mathbb{R}}^d,$
as well as $\x$ and $\y$ of the same dimension, and $\te=\x+\bm{j}\y$
restricted to
the cross-product of $d$ Euclidean upper
half-planes, usually referred to as a {\em tube} see \cite{SteinWeiss}[p.~90],
denoted $T_{\Gamma}^{(d)}.$  

We shall be quite careful to distinguish
between hyperanalytic functions, and their limits as $\y\rightarrow {\bf
0}^+.$
For future reference we make the following set of definition:
\begin{definition}{The hyperanalytic function. \label{hyperanfunction}}\\
Any vector-valued function in spatial variable $\x,$ with associated $p$
dimensional auxiliary
variable $\y,$ that satisfies a given d-D generalisation of the Cauchy-Riemann equations for $\y\in \Gamma^{(p)}$ is denoted a hyperanalytic function.
\end{definition}
In analogue to the 1-D analytic signal, we define the hyperanalytic signal.
\begin{definition}{The hyperanalytic signal. \label{hyperansignal}}\\
Any vector-valued function denoted $g^{(+\dots+)}(\x),$ that can be written as the limit
as $\y\rightarrow {\bf 0}^+$ of a hyperanalytic function $k_g^{(+\dots+)}(\x,\y),$
with $p$-dimensional auxiliary variable $\y$ is denoted a hyperanalytic signal.
\end{definition}
$g^{(+\dots+)}(\x)$ will have $p$ number of $+$'s to denote the dimensionality
of the auxiliary variable $\y.$ 
In 1-D we found that evaluating $g^{(+)}(x_1)$ at $x_1$ taking a complex value $t_1,$
the analytic {\em signal} still corresponded to an analytic {\em function} in either the upper half plane. 
This will not
in general be the case in higher dimensions and it will not be natural to
represent the spatial variable as a complex-valued quantity: furthermore
in
some instances the dimension of $\x$ and $\y$ will not be the same.
For future reference, to treat the rotation of the CWT, we also define a hyperanalytic function and signal in a rotated
frame of reference.
\begin{definition}{The $\theta$-hyperanalytic function. \label{thetahyperanfunction}}\\
Any vector-valued function in spatial variable $\x^{\prime}=\bm{r}_{-\theta}\x,$
with associated $p$
dimensional auxiliary
variables $\y,$ that satisfies a given d-D generalisation of the Cauchy-Riemann
equations in variable $\x^{\prime}$ when $\y$ is in $\Gamma^{(p)}(\theta)$ is denoted
a $\theta$-hyperanalytic function. 
\end{definition}
The $\theta$-hyperanalytic function, once the limit $\y\rightarrow {\bf 0}^+$
is taken, yields a $\theta$-hyperanalytic signal. Note that the choice of
$\y$ may depend on $\theta.$
The suitable choice of hyperanalytic to represent local image phenomena depends
on the local dimensionality of the
analysed image. To be able to make a decision on a suitable local representation,
the local dimensionality of a narrowband oscillatory signal
deserves some further discussion. 

\subsection{2-D Oscillations \& narrowband signals}
In 2-D, assuming oscillations associated with a
given period are observed, the 
local dimensionality of the image may be of many different forms. 
Define $\f_0(\phi_0)=\left[f_0\cos(\phi_0)\quad
f_0\sin(\phi_0)\right]^T,$ and note that, if components are only present with
period
$1/f_0$, this would imply the signal takes the form:
\begin{equation}
\label{continuum}
g(\x)=\int_0^{2\pi} a_g(\phi_0) \cos(2\pi  \x^T\f_0(\phi_0))d\phi_0.
\end{equation}
Thus the general representation of the image corresponds to a
continuum of
directions, and $a_g(\phi_0)$ is non-zero for many $\phi_0\in\left[0,2\pi\right).$
Often locally the variation
of the image is far more structured than the
general form of equation (\ref{continuum}) and admits a much simplified representation.
For example, typical patterns such as plane waves would correspond to a single direction
in the signal, or with a slight generalization of slow spatial modulation
{\em i.e.}:
\begin{equation}
\label{planeoscillation}
g(\x)=a_g(\phi_0(\x),\x) \cos(2\pi  \x^T\f_0(\phi_0(\x),\x)),\end{equation}
whilst textured features \cite{Bovik1991} often take the form of an aggregation of plane waves:
\begin{equation}
\label{separableoscillation}
g(\x)=\sum_{l=1}^{L} a_g(\phi_0^{(l)}(\x),\x) \cos(2\pi  
\x^T\f_0(\phi_0^{(l)}(\x),\x)).\end{equation}
Not infrequently we find that $L=2,$ and if this is the case then with $a_g(\phi_0^{(2)}(\x),\x)
=a_g(\phi_0^{(1)}(\x),\x)$ and $\phi_0^{(2)}(\x)=\phi_0^{(2)}(\x)+\pi/2$
and writing $\x=x\left[\cos(\chi)\quad
\sin(\chi)\right]^T,$ we obtain
\begin{equation}
\label{separablexxx}
g(\x)=2a_g(\phi_0^{(1)}(\x),\x)\cos(f_0 x \cos(\chi+\phi_0^{(1)}(\x)))\cos(f_0(\x) x \sin(\chi+\phi_0^{(1)})),\end{equation}
and the signal is separable.
Key local characterising features of the image
will correspond to $L$, the local dimensionality of the oscillation at spatial
point $\x$, $f_0(\x)$, the local period, and
the form of the $L$ pairs $\left\{a_g(\phi_0^{(l)}(\x),\x),\phi^{(l)}(\x)\right\}$,
that is
the magnitude
and direction of the $l$th component. We consider three special forms of
local structure: the local plane wave, that is $L=1$,
local separable behaviour, that is
$L=2$ and $\phi_0^{(2)}(\x)=\phi_0^{(1)}(\x)+\pi/2$, and $L>1$ 
with no special relationship between  $\phi_0^{(2)}(\x)$ and $\phi_0^{(1)}(\x).$
In the latter case there is
no local separable structure present. 

The analysis problem then corresponds to, 
if confronted with a single real image that naturally
fits into the model of equation (\ref{planeoscillation}), or for that matter
equation (\ref{separableoscillation}), 
determining
$a_g^{(l)}(\phi_0(\x),\x)$, $f_0(\x)$, and $\phi_0^{(l)}(\x)$ from the single observed image, or equivalently defining the 2-D set of functions playing
the role of $g^{(+)}(x_1),$ from
$g(\x).$ Once this has been achieved, the characterising parameters may be
determined from $g^{(+)}(\x).$ In analogue with 1-D, a natural starting point
for such determination,
is a sinusoidal decomposition.

\subsection{Quaternionic Fourier Transforms}
The FT decomposes structure into plane waves of $\exp(2\pi \bm{j} \f^T\x),$
\cite{BulowPhd} and does not represent separable structures as a single coefficient.
A natural tool for 
analysis of separable oscillations is instead the Quaternionic Fourier Transform
(QFT) \cite{BulowPhd}. Instead of using complex numbers, the QFT is defined in terms of quaternionic
units \cite{Deavours,sudbery}.
An arbitrary quaternionic object
takes the form $q=q_1+q_2
\bm{i}+q_3\bm{j}+q_4 \bm{k}\in {\mathbb{H}}$ where $q_i\in{\mathbb{R}},\;i=1,\dots
4$ and
$\bm{i}^2=\bm{j}^2=\bm{k}^2=\bm{i}\bm{j}\bm{k}=-1,$
whilst $\bm{ij}=-\bm{ji}=\bm{k},$ $\bm{ik}=-\bm{ki}=-\bm{j},$ 
and finally $\bm{jk}=-\bm{kj}=\bm{i}.$ $q$ has conjugate
$q^{*}=q_1-q_2
\bm{i}-q_3\bm{j}-q_4 \bm{k},$ and $q$ can be written as
$q=\Re\left\{q\right\}+{\mathscr{P}}u\left\{q\right\},$ where
$\Re\left\{q\right\}=\frac{1}{2}\left[q+q^*\right],$ is the real part and
${\mathscr{P}}u\left\{q\right\}=\frac{1}{2}\left[q-q^*\right],$ is the pure quaternion part. The QFT is then defined in 2-D by \cite{PeiDingChang2001}:
\begin{equation}
\label{qftdefin}
G_Q(\q)={\cal{F}}_{\cal Q}\left\{g
\right\}(\q)=
\int_{\mathbb{R}^2}\exp(-2\pi \bm{i} x_1 q_1)
g(\x) \exp(-2\pi \bm{j} x_2 q_2)\;d^2\x,
\end{equation}
and inverted by
\begin{equation}
g(\x)={\cal{F}}^{-1}_{\cal Q}\left\{G_Q
\right\}(\x)=\int_{\mathbb{R}^2}\exp(2\pi \bm{i} x_1 q_1)
G_Q(\q) \exp(2\pi \bm{j} x_2 q_2)\;d^2\x.
\end{equation}
The QFT can recover separate parity structure in $x_1$ and $x_2$ as it separately
records four real values at each quaternion frequency $\q,$ thus forming
a natural analysis tool for structure separable over $x_1$ and $x_2$ in terms
of $\cos(2\pi x_1 q_1)\cos(2\pi x_2 q_2),$ $\cos(2\pi x_1 q_1)\sin(2\pi x_2 q_2),$ $\sin(2\pi x_1 q_1)\cos(2\pi x_2 q_2),$ and $\sin(2\pi x_1 q_1)\sin(2\pi x_2 q_2),$ each corresponding to a separable function in $x_1$ and $x_2,$
that is either even or odd.
A local representation of this form would require defining more than a single additional
component $g^{(1)}(\x)$ at $\x$ as both local rates corresponding to $q_1$
and $q_2$ would
need determination.

The Fourier
transform of a signal $g(\x)$ can also be defined in an arbitrary pure unit
quaternion $\bm{e}=
\bm{i}e_1+\bm{j}e_2+\bm{k}e_3,$ where $\sum e_j^2=1,$ can also be calculated,
and is given
by:
\begin{equation}
\label{dirFourier}
G_{\bm{e}}(\f)=\int_{-\infty}^{\infty} \int_{-\infty}^{\infty} g(\x) e^{-2\pi
\bm{e} \f^T \x}\;d^2\x,\quad
g(\x)=\int_{-\infty}^{\infty} \int_{-\infty}^{\infty} G_{\bm{e}}(\f) e^{2\pi
\bm{e} \f^T \x}\;d^2\f.
\end{equation}
We will refer to the above transformation (\ref{dirFourier}) as the Unit Quaternion Fourier Transform (UQFT), but it is also known as the Type 3 QFT,
see \cite{PeiDingChang2001}[p.~2784]. The UQFT is still interpretable as
a sinusoidal decomposition of structure, as any pure unit $\bm{e}$ satisfies
the De
Moivre relationships \cite{Hamilton}. The representation is now in terms
of plane waves, {\em i.e.} $\cos(2\pi\f^T \x)$ and $\sin(2\pi\f^T \x),$ and would
require more than a single component $g^{(1)}(\x)$ as both the direction and period of $\f$
would require determination from the signal.
For signals with simple local dimensionality, the QFT and the UQFT form
natural representations of the signal. Using the FT allowed us to obtain
the analytic signal/function representation of $c(x_1)$ via the stationary
phase approximation and we shall
discuss the hyperanalytic signals/functions naturally linked with the QFT
and UQFT.

\subsection{The hypercomplex signal}
To extend the analytic signal to 2-D von B{\"u}low and Sommer \cite{BulowSommer1}
used partial
HTs \cite{Hahn} in both of the spatial
directions, thus defining the hypercomplex signal.
The partial HT of $g(\x)$ in either of $x_l,\;l=1,\;2,$ when
the second
argument is fixed, is denoted ${\mathcal{H}}_l\left\{g\right\}(\x),\;l=1,\;2$
and the total
Hilbert transform is formed by
${\mathcal{H}}_T\left\{g\right\}(\x)=
{\mathcal{H}}_2\left\{{\mathcal{H}}_1\left\{g\right\}\right\}(\x).$ ${\mathcal{H}}_T\left\{g\right\}(\x)$
corresponds to the
consecutive operation of two
partial HTs \cite{Hahn}.
\begin{definition}{The hypercomplex signal.}\\
The hypercomplex signal extension of an image $g(\x)$ is defined
by
\begin{equation}
\label{hypercdef}
g^{(++)}(\x)={\cal HC}g(\x)=g(\x)+\bm{i}g^{(1)}_H(\x)+\bm{j}g^{(2)}_H(\x)+
\bm{k}g^{(3)}_H(\x),\end{equation}
where
\[
g^{(1)}_H(\x)={\mathcal{H}}_1\left\{g\right\}(\x),\quad
g^{(2)}_H(\x)={\mathcal{H}}_2\left\{g\right\}(\x)\quad
g^{(3)}_H(\x)={\mathcal{H}}_T\left\{g\right\}(\x).
\]
\end{definition}
To confirm that the hypercomplex signal corresponds to a valid hyperanalytic
signal
we firstly need to define a hyperanalytic
function and then
introduce
two auxiliary variables, $y_1$ and $y_2.$ We
define the set of functions 
\begin{eqnarray*}
u_{g;H}(\x,\y)& =p(x_1,y_1)*p(x_2,y_2)*g(\x),\quad
v^{(1)}_{g;H}(\x,\y)&=q_1(x_1,y_1)*p(x_2,y_2)*g(\x)\\
v^{(2)}_{g;H}(\x,\y)& =p(x_1,y_1)*q_1(x_2,y_2)*g(\x),\quad
v^{(3)}_{g;H}(\x,\y) &=q_1(x_1,y_1)*q_1(x_2,y_2)*g(\x),
\end{eqnarray*}
where $p(x,y)$ and $q_1(x,y)$ are given by equation (\ref{1Dpoisson}).
Collect the spatial and auxiliary variables in the vector-valued variable $\te=\x+\bm{j}\begin{pmatrix}
\bm{k}y_1 & y_2 \end{pmatrix}^T$ \footnote{Formally $t_1=x_1+\bm{i}y_1,$ but any number system of the
form of a real number plus a scaled pure unit quaternion is isomorphic to ${\mathbb{C}},$
and so
we will relax this formality and not distinguish ${\mathbb{C}}=\left\{x+\bm{j}y \right\}$
from $\left\{x+\bm{i}y \right\}.$ \label{footie}}.
Define the
quaternionic (or vector-valued) function:
\begin{eqnarray}
k_{g;H}^{(++)}(\x,\y)&=&u_{g;H}(\x,\y)+\bm{i}v^{(1)}_{g;H}(\x,\y)
+\bm{j}v^{(2)}_{g;H}(\x,\y)+\bm{k}v^{(3)}_{g;H}(\x,\y).
\end{eqnarray}
It is convenient to express $k_{g;H}^{(++)}(\x,\y)$ in terms of its Cayley-Dickson forms \cite{Moxey}[p.~1952]. 
We define $u_{g;H}^{(i)}(\x,\y)=u_{g;H}(\x,\y)+\bm{j}v^{(2)}_{g;H}(\x,\y),$
$v^{(i)}_{g;H}(\x,\y)=v_{g;H}^{(1)}(\x,\y)+\bm{j}v^{(3)}_{g;H}(\x,\y),$
$u_{g;H}^{(j)}(\x,\y)=u_{g;H}(\x,\y)+\bm{j}v^{(1)}_{g;H}(\x,\y)$ and
$v^{(j)}_{g;H}(\x,\y)=v_{g;H}^{(2)}(\x,\y)+\bm{j}v^{(3)}_{g;H}(\x,\y).$
The functions $u_{g;H}^{(\cdot)}(\x,\y)$ and
$v_{g;H}^{(\cdot)}(\x,\y),$ with $\cdot=i,\;j,$ are the simplex and perplex
respectively,
of the Cayley-Dickson forms.
Thus it in fact follows that
\begin{equation}
k_{g;H}^{(++)}(\x,\y)=u_{g;H}^{(i)}(\x,\y)+\bm{i}v^{(i)}_{g;H}(\x,\y)=
u_{g;H}^{(j)}(\x,\y)+\bm{j}v^{(j)}_{g;H}(\x,\y).
\end{equation}
$u_{g;H}^{(i)}(\te),\;v^{(i)}_{g;H}(\te)\in {\mathbb{C}}$ and the
pair of functions
satisfies the Cauchy-Riemann equations in $t_1$ for fixed values of $t_2,$
in the upper half plane of $y_1>0,$ whilst $u_{g;H}^{(j)}(\te)\in {\mathbb{C}},\;v^{(j)}_{g;H}(\te)\in {\mathbb{C}},$
and the pair satisfies the Cauchy-Riemann equations in $t_2$ for fixed values of $t_1,$
in the upper half plane of $y_2>0.$ The hypercomplex system of equations
is defined as:
\begin{eqnarray}
\label{hypercomplexsystem}
&&\frac{\partial u_{g;H}^{(i)}(\te)}{\partial y_1}+\frac{\partial v_{g;H}^{(i)}(\te)}{\partial x_1}=0,\quad\frac{\partial u_{g;H}^{(i)}(\te)}{\partial x_1}-\frac{\partial v_{g;H}^{(i)}(\te)}{\partial y_1}=0\\
&&\frac{\partial u_{g;H}^{(j)}(\te)}{\partial y_2}+
\frac{\partial v_{g;H}^{(j)}(\te)}{\partial x_2}=0,
\quad\frac{\partial u_{g;H}^{(j)}(\te)}{\partial x_2}-\frac{\partial v_{g;H}^{(j)}(\te)}{\partial y_2}=0.
\nonumber
\end{eqnarray}
This set of equations provide a valid 2-D generalisation of the Cauchy-Riemann equations, and thus, as the set of functions constituting 
$k_{g;H}^{(++)}(\x,\y)$ satisfy equations (\ref{hypercomplexsystem}) for $(\x,\y)$ taking values in the tube $T_{\Gamma},$ we may denote the
latter function
hyperanalytic, {\em cf } definition (\ref{hyperanfunction}).
Clearly as $\y\rightarrow \mathbf{0}^+$,
$u_{g;H}(\te)\rightarrow g(\x),$ $v^{(1)}_{g;H}(\te)\rightarrow g^{(1)}_H(\x),$
$v^{(2)}_{g;H}(\te)\rightarrow g^{(2)}_H(\x),$ and
$v^{(3)}_{g;H}(\te)\rightarrow g^{(3)}_H(\x),$ and thus 
$\lim_{\y\rightarrow \mathbf{0}^+} k_{g;H}^{(++)}(\x,\y) \rightarrow g^{(++)}(\x),$
and so
$g^{(++)}(\x)$ corresponds to a hyperanalytic signal 
({\em cf } definition (\ref{hyperansignal})).\\
Furthermore, the hypercomplex signal $g^{(++)}(\x),$ when viewed in the quaternionic frequency
domain \cite{BulowSommer1}[p.~2849], satisfies $
G^{(++)}_{Q}(\q)=0\quad \forall \q\notin {\mathbb{R}}^{+2},$
this providing a 2-D version of equation (\ref{HT}).

To be able to consider separable structures in an arbitrary direction we extend the definition of a hypercomplex signal 
to a $\theta$-hypercomplex signal.
\begin{definition}{The $\theta$-hypercomplex signal.}\\
The $\theta$-hypercomplex signal of a real-valued image $g(\x)$ is defined
for any $\theta\in\left[0,2\pi\right),$ as
\begin{eqnarray}
\nonumber
g_{\theta}^{(++)}(\x)&=&{\cal HC}_{\theta}g(\x)
=g(\bm{r}_{-\theta}\x)+\bm{i}g^{(1)}_H(\bm{r}_{-\theta}\x)
+\bm{j}g^{(2)}_H(\bm{r}_{-\theta}\x)+
\bm{k}g^{(3)}_H(\bm{r}_{-\theta}\x)\\
\nonumber
&=&R_{\theta}g(\x)+\bm{i}R_{\theta}{\cal H}_1 g(\x)
+\bm{j}R_{\theta}{\cal H}_2 g(\x)
+\bm{k}R_{\theta}{\cal H}_1{\cal H}_2 g(\x)\\
&=&g_{\theta}(\x)+\bm{i}g^{(1)}_{\theta;H}(\x)+\bm{j}g^{(2)}_{\theta;H}(\x)
+\bm{k}g^{(3)}_{\theta;H}(\x),
\label{thetahyper}
\end{eqnarray}
\end{definition}
thus defining the functions $g^{(s)}_{\theta;H}(\x),$ for $s=1,\;2,\;3.$
The $\theta$-hypercomplex signal is suitable for analysing structure that
is separable,  or aligned with variations,
in a rotated frame of reference. To be able to compactly note the QFT of
$g^{(++)}_{\theta}(\x)$
define the reflection operator ${\cal J}_{\theta}.$ The reflection of $g(\x)$ in the line that makes an angle $\theta$ with the $x_1$ axis is given by
${\cal J}_{\theta}
g(\x)=g(\bm{J}_{\theta}\x),$ using matrix 
$\bm{J}_{\theta}= \left[\left\{ \cos(2\theta), \sin(2\theta) \right\}, \left\{\sin(2\theta), -\cos(2\theta) \right\} 
\right],$ where  $\theta \in \left[0,2\pi\right).$ 
Then:
\begin{eqnarray}
\label{thetahyperft}
{\cal{F}}\left\{g^{(++)}_{\theta}\right\}(\f)&=&
\left[1+{\mathrm{sgn}}\left(\left[\bm{r}_{-\theta}\f\right]_1\right)
\right]\left[1-\bm{k}{\mathrm{sgn}}\left(\left[\bm{r}_{-\theta}\f\right]_2\right)
\right]G(\bm{r}_{-\theta}\f)
\\
\nonumber
{\cal{F}}_{\cal Q}\left\{g^{(++)}_{\theta}\right\}(\q)
&=&\frac{1-\bm{k}}{2}\left[1+{\mathrm{sgn}}\left(
\left[\bm{r}_{-\theta}\q\right]_1 \right)\right]\left[1+{\mathrm{sgn}}\left(
\left[\bm{r}_{-\theta}\q\right]_2 \right)\right]G\left(\bm{r}_{-\theta}\q\right)\\
&&+\frac{1+\bm{k}}{2}\left[1+{\mathrm{sgn}}\left(
\left[\bm{r}_{\theta}\q\right]_1 \right)\right]\left[1+{\mathrm{sgn}}\left(
\left[\bm{r}_{\theta}\q\right]_2 \right)\right]G\left(\bm{r}_{-\theta}\bm{J}_{\pi/2}\q\right)
\label{thetahypercquat}
\end{eqnarray}
The result follows by Appendix \ref{qftofg++}, combined with the fact that
the QFT of a real-valued signal $g(\x)$ can be calculated from its Fourier transform (see Hahn \& Snopek \cite{HahnSnopek2} as well as Pei, Ding and
Chang \cite{PeiDingChang2001}). Thus the FT of $g^{(++)}_{\theta}(\x)$ is
supported on $\left[\bm{r}_{-\theta}\f\right]_2>0,$ whilst as noted in the
appendix, the QFT if implemented in a rotated frame is supported wholly on
one quadrant.
Most importantly, as the QFT can be inverted into
the spatial domain signal, if a signal admits a QFT representation
of equation (\ref{qffft}) then it is the QFT of a $\theta$-hypercomplex signal
identified by the form of $G(\bm{r}_{-\theta}\q).$ This
result can be used to interpret the properties of the quaternionic wavelet coefficients.

We shall define analysis methods suitable for signals observed
in an alternate frame of reference, as the local orientation of a signal
is a very important local characteristic.
Define the signal $\widetilde{g}_{-\theta}(\x)$ via
\begin{equation}
\label{rotatedsignal}
\widetilde{g}_{-\theta}(\x)=g(\bm{r}_{\theta}\x),\quad
g(\x)=\widetilde{g}_{-\theta}(\bm{r}_{-\theta}\x),\;\theta\in[0,2\pi),
\end{equation}
so that
$
\widetilde{g}_{-\theta;\theta}^{(++)}(\x)=\widetilde{g}_{-\theta}(\bm{r}_{-\theta}\x)+
\bm{i}\widetilde{g}^{(1)}_{-\theta;H}(\bm{r}_{-\theta}\x)
+\bm{j}\widetilde{g}^{(2)}_{-\theta;H}(\bm{r}_{-\theta}\x)+
\bm{k}g^{(3)}_{-\theta;H}(\bm{r}_{-\theta}\x),$ or
$\widetilde{g}_{-\theta;\theta}^{(++)}(\x)=g(\x)+
\bm{i}\widetilde{g}^{(1)}_{-\theta;H}(\bm{r}_{-\theta}\x)\\
+\bm{j}\widetilde{g}^{(2)}_{-\theta;H}(\bm{r}_{-\theta}\x)+
\bm{k}g^{(3)}_{-\theta;H}(\bm{r}_{-\theta}\x).
$
Of course,
$\widetilde{g}_{-\theta}(\bm{r}_{-\theta}\x)=R_{\theta}R_{-\theta}g(\x)=g(\x),$
but
$
\widetilde{g}_{-\theta,\theta;H}^{(s)}=\widetilde{g}^{(s)}_{-\theta;H}(\bm{r}_{-\theta}\x)=R_{\theta}{\cal H}_s R_{-\theta}g(\x)
\neq  g^{(s)}_{H}(\x),\;s=1,\;2,$ {\em etc.}
Hence $\widetilde{g}_{-\theta}(\bm{r}_{-\theta}\x)$
corresponds to the observed signal, whilst $
\widetilde{g}_{-\theta,\theta;H}^{(s)}$ corresponds to the HTs
in a rotated frame of reference observed in our frame of reference. 

\begin{definition}{$\theta$-hypercomplex: hypercomplex vs anti-hypercomplex decomposition.}\\
The decomposition of real-valued $g(\x)$ into four quaternionic components
$\widetilde{g}_{-\theta;\theta}^{(++)}(\x),$ $\widetilde{g}_{-\theta;\theta}^{(-+)}(\x),$
$\widetilde{g}_{-\theta;\theta}^{(+-)}(\x),$ and $\widetilde{g}_{-\theta;\theta}^{(--)}(\x):$
\begin{eqnarray}
g(\x)
&=&\frac{1}{4}\left(\widetilde{g}_{-\theta;\theta}^{(++)}(\x)+
\widetilde{g}_{-\theta;\theta}^{(-+)}(\x)+
\widetilde{g}_{-\theta;\theta}^{(+-)}(\x)+
\widetilde{g}_{-\theta;\theta}^{(--)}(\x)
\right),
\label{hypersplit}
\end{eqnarray}
where we define $\widetilde{g}_{-\theta;\theta}^{(++)}(\x)=g(\x)+
\bm{i}\widetilde{g}^{(1)}_{-\theta;H}(\bm{r}_{-\theta}\x)
+\bm{j}\widetilde{g}^{(2)}_{-\theta;H}(\bm{r}_{-\theta}\x)+
\bm{k}g^{(3)}_{-\theta;H}(\bm{r}_{-\theta}\x),$ 
$\widetilde{g}_{-\theta;\theta}^{(-+)}(\x)=g(\x)-
\bm{i}\widetilde{g}^{(1)}_{-\theta;H}(\bm{r}_{-\theta}\x)
+\bm{j}\widetilde{g}^{(2)}_{-\theta;H}(\bm{r}_{-\theta}\x)-
\bm{k}g^{(3)}_{-\theta;H}(\bm{r}_{-\theta}\x),$
$\widetilde{g}_{-\theta;\theta}^{(+-)}(\x)=g(\x)+
\bm{i}\widetilde{g}^{(1)}_{-\theta;H}(\bm{r}_{-\theta}\x)
-\bm{j}\widetilde{g}^{(2)}_{-\theta;H}(\bm{r}_{-\theta}\x)-
\bm{k}g^{(3)}_{-\theta;H}(\bm{r}_{-\theta}\x)
$ and $\widetilde{g}_{-\theta;\theta}^{(--)}(\x)=g(\x)-
\bm{i}\widetilde{g}^{(1)}_{-\theta;H}(\bm{r}_{-\theta}\x)
-\bm{j}\widetilde{g}^{(2)}_{-\theta;H}(\bm{r}_{-\theta}\x)+
\bm{k}g^{(3)}_{-\theta;H}(\bm{r}_{-\theta}\x).$ 
 will be denoted as the $\theta$-Hypercomplex analytic/anti-analytic decomposition of $g(\x).$ The
four components $\widetilde{g}_{-\theta;\theta}^{(++)}(\x),$ $\widetilde{g}_{-\theta;\theta}^{(-+)}(\x),$
$\widetilde{g}_{-\theta;\theta}^{(+-)}(\x),$ and $\widetilde{g}_{-\theta;\theta}^{(--)}(\x),$
are denoted as the $\theta$-hypercomplex
signal, the $\theta$ first anti-hypercomplex
signal, the $\theta$ second anti-hypercomplex
signal and the $\theta$ third anti-hypercomplex
signal, respectively. 
\end{definition}
In
the rotated frame of reference, equation (\ref{hypersplit}) provides a 2-D
analogue of the analytic/anti-analytic decomposition of a 1-D real-valued
signal \cite{Papoulis}. We 
may use the $\theta$-hypercomplex signal to decompose any real image $g(\x)$ into
four components, where each component in a rotated frame of reference is
analytic/anti-analytic in each spatial variable {\em separately}. 
We characterize local structure in terms of magnitude and phase of 
$\widetilde{g}_{-\theta;\theta}^{(++)}(\x),$ interpretable in 
the rotated frame of reference. Any of
$\widetilde{g}_{-\theta,\theta}^{(\pm \pm)}(\x)$ if calculated at all 
$\theta=\left[0,\pi/2\right],$ will enable us to parameterise
separable oscillations.
The $\theta$-hypercomplex signal is introduced in this article for the purpose
of recognizing separable or univariate local structure that is not aligned
with the axis of observation; the act of rotations is integral to defining
the CWT in subsequent sections, and so its introduction is
a necessity.

\subsection{Properties of the hypercomplex signal \label{prophyp}}
A hypercomplex signal can be decomposed in terms of a modulus and a set of
phase/s,
using quaternionic polar coordinates,
see discussion by \cite{BulowSommer}.
The modulus is still a real number and can be interpreted as local energy
content whereas the phase/s are given by a unit energy quaternion and
contain additional information about the structure of the 
image. This, extended to the $\theta$-hypercomplex signal, corresponds to a polar representation of variation in orientation $\theta$.

\begin{definition}{The polar representation of the $\theta$-hypercomplex
signal.}\\
The $\theta$-hypercomplex image will be represented in terms of its
amplitude, and phases. Define the unit quaternion by
$q^{(U)}_{\theta}(\x)=q_{\theta}^{(1)}(\x)+\bm{i}q_{\theta}^{(2)}(\x)+
\bm{j}q_{\theta}^{(3)}(\x)+\bm{k}q_{\theta}^{(4)}(\x)=
e^{2\pi\bm{i}\alpha_{\theta}(\x)}
e^{2\pi\bm{k}\gamma_{\theta}(\x)}e^{2\pi\bm{j}\beta_{\theta}(\x)},$ then
the polar representation
corresponds to:
\begin{eqnarray}
g_{\theta}^{(++)}(\x)&=&
\left|g_{\theta}^{(++)}(\x)\right|e^{2\pi\bm{i}\alpha_{\theta}(\x)}
e^{2\pi\bm{k}\gamma_{\theta}(\x)}e^{2\pi\bm{j}\beta_{\theta}(\x)}
=\left|g_{\theta}^{(++)}(\x)\right|q^{(U)}_{\theta}(\x).
\label{polarhyper}
\end{eqnarray}
\end{definition}
Von B{\"u}low and Sommer considered
only the hypercomplex signal ({\em i.e.} $\theta=0$)
but
the exact algebraic forms of the angles for both $\alpha_{\theta}(\x)$ and $\beta_{\theta}(\x)$ in terms of $g^{(l)}_{\theta,H}(\x),\;l=1,\;2,\;3,$
can be found by substituting in the form for the quaternion components
into the formulae given by \cite{BulowSommer1}[p.~2849]:
\begin{eqnarray}
\nonumber
\alpha_{\theta}(\x)&=&\frac{1}{4\pi}\tan^{-1}\left(2\frac{q_{\theta}^{(3)}(\x)q_{\theta}^{(4)}(\x)
+q_{\theta}^{(1)}(\x)q_{\theta}^{(2)}(\x)}{q_{\theta}^{(1)2}(\x)+
q_{\theta}^{(3)2}(\x)-
\left(q_{\theta}^{(2)2}(\x)+
q_{\theta}^{(4)2}(\x)\right)}
\right), \nonumber \\
\nonumber
\beta_{\theta}(\x)&=&\frac{1}{4\pi}\tan^{-1}\left(2\frac{q_{\theta}^{(2)}(\x)q_{\theta}^{(4)}(\x)
+q_{\theta}^{(1)}(\x)q_{\theta}^{(3)}(\x)}{q_{\theta}^{(1)2}(\x)+
q_{\theta}^{(2)2}(\x)-
\left(q_{\theta}^{(3)2}(\x)+
q_{\theta}^{(4)2}(\x)\right)}
\right), \nonumber \\
\label{angles}
\gamma_{\theta}(\x)&=&-\frac{1}{4\pi}\sin^{-1}\left(2\left[
q_{\theta}^{(2)}(\x)q_{\theta}^{(3)}(\x)-q_{\theta}^{(1)}(\x)q^{(4)}(\x)
\right] \right).
\end{eqnarray}
More than one phase function are used, as
variations in either axes are given by $\alpha_{\theta}(\x)$ and 
$\beta_{\theta}(\x)$ respectively,
whilst $\gamma_{\theta}(\x)$ is considered by von B{\"u}low and Sommer as
a measure of
`degree of separability'. If the signal analysed takes the form
$g(\x)=g_1(x_1)g_2(x_2),$ {\em i.e.} the
signal is separable, then they note that \cite{BulowSommer1}[p.~2849]:
$
\alpha_{0}(\x)=\tan^{-1}\left(\frac{{\cal{H}}_1\left\{g\right\}(\x)}{g(\x)}\right),$
$\beta_{0}(\x)=\tan^{-1}\left(\frac{{\cal{H}}_2\left\{g\right\}(\x)}{g(\x)}\right),$
and
$\gamma_{0}(\x)=0.$ Even if $g(\x)$ is not separable
the polar representation of equation (\ref{polarhyper}) may be used to represent
an arbitrary real-valued image $g(\x).$
\begin{definition}{The $\theta$-polar representations of g(\x).}\\
The $\theta$-polar representation of a real signal $g(\x)$ is given by:
\begin{eqnarray}
g(\x)&=&\Re\left\{\left|\widetilde{g}_{-\theta,\theta}^{(++)}(\x)\right|
e^{2\pi\bm{i}\widetilde{\alpha}_{-\theta,\theta}(\x)}
e^{2\pi\bm{k}\widetilde{\gamma}_{-\theta,\theta}(\x)}e^{2\pi\bm{j}\widetilde{\beta}_{-\theta,\theta}(\x)}
\right\} \nonumber \\
&=&\left|\widetilde{g}_{-\theta,\theta}^{(++)}(\x)\right|
\left(\cos(2\pi\widetilde{\alpha}_{-\theta,\theta}(\x))\cos(2\pi\widetilde{\gamma}_{-\theta,\theta}(\x))
\cos(2\pi\widetilde{\beta}_{-\theta,\theta}(\x))\right. \nonumber \\
&&\left.+\sin(2\pi\widetilde{\alpha}_{-\theta,\theta}(\x))
\sin(2\pi\widetilde{\gamma}_{-\theta,\theta}(\x))
\sin(2\pi\widetilde{\beta}_{-\theta,\theta}(\x))
\right). 
\label{thetapolargen}
\end{eqnarray}
A separable structure in orientation $\theta$ will have $\gamma_{\theta}(\x)=0$
and the signal is represented by:
\begin{eqnarray}
\label{seppyhyper}
g(\x)&=&\left|\widetilde{g}_{-\theta,\theta}^{(++)}(\x)\right|\cos(2\pi\widetilde{\alpha}_{-\theta,\theta}(\x))
\cos(2\pi\widetilde{\beta}_{-\theta,\theta}(\x)).
\end{eqnarray}
\end{definition}
Note that $\widetilde{\alpha}_{-\theta,\theta}(\x)$ {\em etc} are determined
by equation (\ref{angles}) but where now the $\theta$-hypercomplex extension
of $\widetilde{g}_{-\theta}(\x)$ is used as the basis of the representation.
To see the utility of such a representation consider analysis of a signal
separable in the rotated frame of reference, {\em i.e.} 
\begin{equation}
\label{sepintheta}
g(\x)=\widetilde{g}_{-\theta}(\bm{r}_{-\theta}\x)=\widetilde{g}_{1,-\theta}
\left(\left[\bm{r}_{-\theta}\x\right]_1\right)
\widetilde{g}_{2,-\theta}\left(\left[\bm{r}_{-\theta}\x\right]_2\right).
\end{equation}
Under the assumption that $g(\x)$ satisfies equation (\ref{sepintheta}),
we determine that the $\theta$-hypercomplex extension of the signal is
given (in terms of $\x^{\prime}=\bm{r}_{-\theta}\x$ ) by:
\begin{eqnarray}
\widetilde{g}_{-\theta,\theta}^{(++)}(\x)&=&\widetilde{g}_{-\theta}(\x^{\prime})+
\bm{i}{\cal H}_1\left\{\widetilde{g}_{-\theta}\right\}
(\x^{\prime}) +\bm{j}{\cal H}_2\left\{\widetilde{g}_{-\theta}\right\}
(\x^{\prime})
+\bm{k}{\cal H}_2{\cal H}_1\left\{\widetilde{g}_{-\theta}\right\}
(\x^{\prime}).
\end{eqnarray}
${\cal H}_1\left\{\widetilde{g}_{-\theta}\right\}(\x)={\cal H}\left\{\widetilde{g}_{1,-\theta}\right\}
(x_1)
\widetilde{g}_{2,-\theta}(x_2),$\\ ${\cal H}_2\left\{\widetilde{g}_{-\theta}\right\}(\x)={\cal H}\left\{\widetilde{g}_{2,-\theta}\right\}
(x_2)
\widetilde{g}_{1,-\theta}(x_1)$ and ${\cal H}_2{\cal H}_1\left\{\widetilde{g}_{-\theta}\right\}(\x)={\cal H}\left\{\widetilde{g}_{2,-\theta}\right\}
(x_2)
{\cal H}\left\{\widetilde{g}_{1,-\theta}\right\}(x_1).$ 
Hence, in our extended definition, if $g(\x)\in{\mathbb{R}}$ is separable in some rotated frame of reference $\x^{\prime}=\bm{r}_{-\theta}\x,$ then
$\widetilde{\gamma}_{-\theta,\theta}(\x)=0$
for some $\theta\in[0,\pi/2).$ 

If a single component is present in $g(\x)$
at $\x$ and the aim of the analysis was to determine the orientation of the
separable structure in the image at this location, denoted $\theta^{\prime},$
then
we calculate $\widetilde{g}_{-\theta,\theta}^{(++)}(\x)$ at all angles 
$\theta\in[0,\pi/2),$ and determine the value of $\theta^{\prime}$ by 
$\widetilde{\gamma}_{-\theta^{\prime},\theta^{\prime}}(\x)=0.$ Once this has
been done the signal is well described by its local magnitude, {\em i.e.}
\[\left|\widetilde{g}_{-\theta,\theta}^{(++)}(\x)\right|=\sqrt{\left[\left(
\widetilde{g}_{1,-\theta}^2(x^{\prime}_1)+
{\cal H}^2\left\{\widetilde{g}_{1,-\theta}\right\}(x^{\prime}_1)
\right) \left(
\widetilde{g}_{2,-\theta}^2(x^{\prime}_2)+
{\cal H}^2\left\{\widetilde{g}_{2,-\theta}\right\}(x^{\prime}_2)
\right)\right]},\]
whilst local variation is best described by the two separable
phase functions 
\[\widetilde{\alpha}_{-\theta,\theta}(\x)=
\tan^{-1}\left( {\cal H}\left\{\widetilde{g}_{1,-\theta}\right\}(x^{\prime}_1)/
\widetilde{g}_{1,-\theta}(x^{\prime}_1)\right),\quad\widetilde{\beta}_{-\theta,\theta}(\x)=
\tan^{-1}\left( {\cal H}\left\{\widetilde{g}_{2,-\theta}\right\}(x^{\prime}_2)/
\widetilde{g}_{2,-\theta}(x^{\prime}_2)\right).\] 
The local frequency is found
in each of the two axes as 
$\frac{\partial}{\partial x_1^{\prime}}\widetilde{\alpha}_{-\theta,\theta}(\x)$
and 
$\frac{\partial}{\partial x_2^{\prime}}\widetilde{\beta}_{-\theta,\theta}(\x).$

\begin{proposition}{Orthogonality of the signal and the HTs.}\\
The original signal $g(\x)$ and $g_H^{(s)}(\x),\;s=1,2,3,$ are all mutually orthogonal.
\end{proposition}
\begin{proof} This is a trivial consequence of the orthogonality of a 1-D
HT.
\end{proof}
Furthermore we note that as $\widetilde{g}_{-\theta}(\x)$ 
is orthogonal to $\widetilde{g}_{-\theta;H}^{(s)}(\x),
\;s=1,2,3,$ we find that $g(\x)$ is orthogonal 
to $\widetilde{g}^{(s)}_{-\theta;H}(\bm{r}_{-\theta}\x),\;s=1,2,3.$
\begin{proposition}{Norm of the hypercomplex components.}\\
The  $g_H^{(s)},\;s=1,2,3,$ have the same norm as $g(\x).$
\end{proposition}
\begin{proof} 
This is a direct result following from the norm of the 1-D HT
of a real-valued function \cite{Hahnbook}[p.~4]. 
\end{proof}\\
Hence in summary: ${\mathscr{P}}u\left\{g^{(++)}(\x)\right\}$ defines a quaternionic
object where each real-valued function multiplying the quaternions $\bm{i},$
$\bm{j}$ and $\bm{k}$ can be considered as objects with the same spatial/spatial
frequency support as $g(\x),$
as the HT constructs an object with the same time-frequency
support as the original signal \cite{Boashash1}. 
Note that all of $\left\{g^{(s)}_{-\theta,\theta}
(\x)\right\}$
for $s=1,\;2,\;T$ are mutually orthogonal. ${\mathscr{P}}u\left\{g^{(++)}_{-\theta,\theta}
(\x)\right\}$
constructs a quaternionic object with the same spatial and spatial frequency support
as $g(\x),$ where the relative magnitudes of the components specifies the
structure of $g(\x)$ in a rotated frame.

\subsection{The monogenic signal}
The hypercomplex signal, once extended to the $\theta$-hypercomplex signals
recognizes locally separable structure, and can if the signal is calculated
for all $\theta\in\left[0,2\pi\right),$ determine locally unidirectional
structure at a given orientation. However, this does require the calculation of the hypercomplex extension at all $\theta.$ A hyperanalytic extension
for locally unidirectional
structure, or plane waves, at any given spatial point can also be made.  Felsberg
\& Sommer, \cite{FelsbergSommer}, introduced a hyperanalytic extension in image processing suitable for this purpose. See also an early work by Dixon \cite{dixon}, where monogenic functions were denoted `pure Hamiltonian' functions
and the discussion by Duffin in \cite{duffin}.
This approach starts from the $d$-D Poisson kernels \cite{SteinWeiss}, a
$d$-D
version of equation (\ref{1Dpoisson}), defined for any dimension $d\ge 1:$
\begin{equation}
p({\bf x},y)  = c_{d}\  \frac{y}{ \left[ |{\bf x}|^{2} + y^{2} \right]^{\frac{d+1}{2}} }, \quad
q_{l}({\bf x},y)   =  c_{d}\ \frac{x_{l}}{ \left[ |{\bf x}|^{2} + y^{2}
\right]^{\frac{d+1}{2}}
} \quad c_{d}=\frac{\Gamma(\frac{n+1}{2})}{\pi^{ \frac{n+1}{2}}},\;l=1,\dots,d.
\label{ddversionhilbert}
\end{equation}
An extra set of $d+1$ functions
\begin{equation}
\label{convolution}
u_{R;g}(\x,y)=\left(p(\cdot,y) \ast \ast 
g(\cdot)\right)(\x),\quad v_{R;g}^{(l)}(\x,y)=\left(q_l \ast \ast g\right)(\x,y),
\end{equation}
are defined from $g(\x).$ For $d=2,$
$\left\{u_{R;g}(\x,y),v_{R;g}^{(1)}(\x,y),v_{R;g}^{(2)}(\x,y)\right\}$
satisfy an alternative set of generalized Cauchy-Riemann equations,
with only a single auxiliary variable $y.$ These generalized
Cauchy-Riemann equations are called 
the Riesz system, see \cite{SteinWeiss}[p.~234].
\begin{definition}{The monogenic function.}\\
Any solution of the Riesz system in the upper half-space ($y>0$) is denoted a monogenic function \cite{FelsbergPhd}[p.~35]. Note that if
$k^{(+)}(\x,y)$
is a monogenic function then
$k^{(+)*}(\x,y)$ is a solution of the Riesz system in the lower half-space ($y<0$),
and will be denoted an anti-monogenic function, with notation 
$k^{(-)}(\x,y).$\end{definition}
For $d=2$, if $k_{R;g}^{(+)}(\x,y)=u_{R;g}(\x,y) + \bm{i} v^{(1)}_{R;g}(\x,y) + 
\bm{j} v^{(2)}_{R;g}(\x,y)$ then it is a monogenic function \cite{SteinWeiss}[p.~235]. 
In the limit of $y \rightarrow 0^{+}$ \cite{SteinWeiss}:
\begin{eqnarray}
\lim_{y \rightarrow 0^{+}} u_{R;g}({\bf x},y) = g({\bf x}),\quad
\lim_{y \rightarrow 0^{+}} v^{(l)}_{R;g}({\bf x},y) = g({\bf x}) * ... * r_{l}({\bf x}) = {\cal R}_{l}g({\bf x}) = g^{(l)}_R({\bf x}),
\label{limits}
\end{eqnarray}
where the quantities $r_l(\x)$ and $R_l(\f)$ are given by:
\begin{equation}
r_{l} =  c_{n}\ \frac{x_{l}}{ |{\bf x}|^{n+1} },\quad R_{l}(\f) = -\bm{j} \ \frac{f_{l}}{f}, \hspace{0.5cm} 
R_{1,Q}(\bm{q}) = -\bm{i} \ \frac{q_{1}}{q}, \hspace{0.5cm}
R_{2,Q}(\f) = -\bm{j} \ \frac{q_{2}}{q} \hspace{0.5cm}l=1,2.
\label{FTRiesz}
\end{equation}
$\left\{r_l(\x)\right\}$ are the Riesz kernels for $d \geq 2,$ and the FT and QFTs of the Riesz kernels are given by \cite{SteinWeiss,Hahn}.
\begin{definition}{The Riesz transform.}\\
The Riesz transform of a signal $g(\x)\in L^2({\mathbb{R}}^2)$ is defined
by operator ${\cal R}$
\begin{eqnarray}
\label{ris1}
{\cal R}\left\{g\right\}(\x)&=&\bm{i}{\cal R}\left\{g\right\}(\x)+
\bm{j}{\cal R}_2\left\{g\right\}(\x),\quad
{\cal R}_l\left\{g\right\}(\x)=g^{(l)}_R(\x)=\left(r_l ** g\right)(\x),\;l=1,\;2.
\end{eqnarray}
\end{definition}
Felsberg
\& Sommer \cite{FelsbergSommer} recently introduced the monogenic signal
into image processing by combining $g(\x)$ with ${\cal R}\left\{g\right\}(\x):$
\begin{definition}{The monogenic signal.}\\
The monogenic, and anti-monogenic signal of real signal
$g(\x),$ are defined by applying operator ${\cal M}^{\pm}$ to signal $g(\x):$
\begin{eqnarray}
g^{(\pm)}(\x) &=&
\label{analytic}
{\mathcal{M}}^{\pm}\left\{g\right\}(\x)
=g(\x)\pm {\mathcal{R}}\left\{g\right\}(\x)
=g(\x)\pm \left[\bm{i} g^{(1)}_R(\x)+
\bm{j} g^{(2)}_R(\x)\right].
\end{eqnarray}
Note that ${\mathcal{R}},$ and $g^{(l)}_R(\x)$ are defined by equation (\ref{ris1}).
\end{definition}
Letting $y\rightarrow0^+$ the monogenic signal is retrieved from the monogenic
function, and so we find that $g^{(+)}(\x)$ is a hyperanalytic signal,
{\em cf} definition (\ref{hyperansignal}). See also the careful discussion in \cite{FelsbergPhd}[p.~32--37] on monogenic signals and the generalized
Cauchy-Riemann equations. 

Only one auxiliary variable is introduced when constructing the Riesz
components, linked with the norm of $\x,$ and the monogenic and anti-monogenic signals
play the roles in 2-D of the analytic and anti-analytic signals
in 1-D. 
%
%
%
To accommodate that the observed image is not observed in the direction of
the local variation, we
introduce here the concept of a $\theta$-monogenic and anti-monogenic signal.

\begin{definition}{The $\theta$-monogenic signal.}\\
The $\theta$-monogenic, and $\theta$-anti-monogenic signal of real signal
$g(\x)$ are for any value of a rotation $\theta\in\left[0,2\pi\right)$ given by:
\begin{eqnarray}
\label{thetaanalytic}
g^{(\pm)}_{\theta}(\x) &=&R_{\theta}{\mathcal{M}}^{\pm}\left\{g\right\}(\x)=
{\mathcal{M}}_{\theta}^{\pm}\left\{g\right\}(\x)
=g(\bm{r}_{-\theta}\x)\pm \left[\bm{i} 
R_{\theta}{\cal R}_1 g(\x)+\bm{j}R_{\theta}{\cal R}_2 g(\x)
\right]\\
&=&g(\bm{r}_{-\theta}\x)\pm \left[\bm{i} g^{(1)}_R(\bm{r}_{-\theta}\x)+
\bm{j} g^{(2)}_R(\bm{r}_{-\theta}\x)\right]
=g_{\theta}(\x)\pm \left[\bm{i}g_{\theta;R}^{(1)}(\x)+\bm{j}g_{\theta;R}^{(2)}(\x)\right],
\nonumber
\end{eqnarray}
thus defining the components $g_{\theta;R}^{(s)}(\x)$ for $s=1,\;2.$
\end{definition}
Note that
\begin{eqnarray}
{\cal F} \left\{g^{(1)}_R(\bm{r}_{-\theta}\x)\right\}
=-\bm{j}\cos(\phi-\theta)G(\bm{r}_{-\theta}\f),\;
{\cal F} \left\{g^{(2)}_R(\bm{r}_{-\theta}\x)\right\}
=-\bm{j}\sin(\phi-\theta)G(\bm{r}_{-\theta}\f).
\label{riesssy5}
\end{eqnarray}
\begin{theorem}{$\theta$-monogenic functions and their monogenicity.\label{thetamonny2}}\\
The $\theta$-monogenic and the $\theta$-anti-monogenic signals $g^{(\pm)}_{\theta}(\x)$ are both the limits of 
sets of functions satisfying the Riesz system of equations in a rotated frame
of reference.
\end{theorem}
\begin{proof}
See section \ref{monoe}.
\end{proof}
Thus a $\theta$-monogenic signal can be considered to be the limit of a monogenic
function in the rotated frame of reference. 
For consistency with our notation
so far, 
if $\theta=0$, the signal is denoted
$g^{(\pm)}(\cdot)$ without subscript. Furthermore we may note that
\begin{eqnarray}
\nonumber
G_{\theta}^{(\pm)}&=&G(\bm{r}_{-\theta}\f)\pm\left(-\bm{k}\cos(\phi-\theta)G(\bm{r}_{-\theta}\f)
+\sin(\phi-\theta)G(\bm{r}_{-\theta}\f)\right)
\\
&=&\left[1\pm\sin(\phi-\theta)\mp\bm{k}\cos(\phi-\theta)
\right]G(\bm{r}_{-\theta}\f).
\label{weirdeqn}
\end{eqnarray}
If the Fourier transform of a function can be written in the form of
equation (\ref{weirdeqn})
it can be directly deduced that the function corresponds to the $\theta$-monogenic image of $g(\x).$ Thus just like we intend to use equation (\ref{thetahyperft}) to recognize a quaternionic object as corresponding to a $\theta$-hypercomplex
object, equation (\ref{weirdeqn}) can be used to determine if a quaternionic object is
$\theta$-monogenic, and what real-valued object plays the role of $g(\x).$

A real-valued image can be decomposed into a $\theta$-monogenic and
$\theta$-anti-monogenic component. Using $\widetilde{g}_{\theta}(\x)$
defined in equation (\ref{rotatedsignal}) we represent $g(\x)$ by:
\begin{eqnarray}
\nonumber
\widetilde{g}_{-\theta,\theta}^{(\pm)}(\x)&=&\widetilde{g}_{-\theta}(\bm{r}_{-\theta}\x)\pm\left(
\bm{i}\widetilde{g}^{(1)}_{-\theta;R}(\bm{r}_{-\theta}\x)
+\bm{j}\widetilde{g}^{(2)}_{-\theta;R}(\bm{r}_{-\theta}\x)\right)\\
&=&g(\x)\pm\left(
\bm{i}\widetilde{g}^{(1)}_{-\theta;R}(\bm{r}_{-\theta}\x)
+\bm{j}\widetilde{g}^{(2)}_{-\theta;R}(\bm{r}_{-\theta}\x)\right).
\label{thetamonogenic}
\end{eqnarray}
${\cal{R}}_1 R_{-\theta}g (\x)=\cos(\theta)R_{-\theta}{\cal{R}}_1 g (\x)+
\sin(\theta)R_{-\theta}{\cal{R}}_2 g (\x),$
and as noted by \cite{SteinWeiss}[p.~241]:
\begin{eqnarray}
\nonumber
\widetilde{g}^{(1)}_{-\theta;R}(\bm{r}_{-\theta}\x)&=&R_{\theta}{\cal{R}}_1 R_{-\theta}g(\x)
=\cos(\theta){\cal{R}}_1 g(\x)+\sin(\theta){\cal{R}}_2 g(\x)
\neq {\cal{R}}_1 g(\x)\\
\widetilde{g}^{(2)}_{-\theta;R}(\bm{r}_{-\theta}\x)&=&R_{\theta}{\cal{R}}_2 R_{-\theta}g(\x)=
-\sin(\theta){\cal{R}}_1 g(\x)+\cos(\theta){\cal{R}}_2 g(\x)
\neq {\cal{R}}_2 g(\x).
\label{notcommuating}
\end{eqnarray}
The rotation operator, $R_{\theta},$ and the Riesz transform operators, ${\cal
R}_l,\;l=1,\;2,$ do not commute, and so
$\widetilde{g}_{\theta,\theta}^{(\pm)}(\x)\neq g^{(\pm)}(\x).$ 
\begin{definition}{$\theta$-monogenic: monogenic/anti-monogenic decomposition.}\\
The decomposition of real-valued $g(\x)$ into two quaternionic components
$\widetilde{g}_{-\theta,\theta}^{(+)}(\x)$ and 
$\widetilde{g}_{-\theta,\theta}^{(-)}(\x):$ 
\begin{eqnarray}
g(\x)&=&\frac{1}{2}\left[
\widetilde{g}_{-\theta,\theta}^{(+)}(\x)+
\widetilde{g}_{-\theta,\theta}^{(-)}(\x)
\right],
\label{decomposition}
\end{eqnarray}
will be denoted the $\theta$-monogenic/anti-monogenic decomposition of $g(\x),$ where
the two components are given by: 
$\widetilde{g}_{-\theta,\theta}^{(\pm)}(\x)=
g(\x)\pm\left(
\bm{i}\widetilde{g}^{(1)}_{-\theta;R}(\bm{r}_{-\theta}\x)
+\bm{j}\widetilde{g}^{(2)}_{-\theta;R}(\bm{r}_{-\theta}\x)\right).$ This is a
natural decomposition of a locally 1-D signal more naturally observed in a rotated frame
of reference. Equation (\ref{decomposition}) forms a 2-D extension of the decomposition of a real-valued
object into a quaternionic form in analogue with the analytic/anti-analytic
decomposition of a 1-D signal, see \cite{Papoulis}.
\end{definition}

\subsection{Properties of the monogenic signal \label{propmon}}
Again, it is natural to represent a hyperanalytic function in terms
of amplitude and local structure, and so 
\cite{FelsbergPhd, FelsbergSommer1, FelsbergSommer2} represent the monogenic
signal in terms of an amplitude and a phase. A
single phase only is used to describe variation as the monogenic representation corresponds to 
a local plane wave structure. The phase represents
the local period in the direction of variation determined from the signal.
\begin{definition}{The polar representation of the $\theta$-monogenic signal.}\\
The $\theta$-monogenic and anti-monogenic images will be represented in terms of the
amplitude, phase and direction by
\begin{equation}
g^{(\pm)}_{\theta}(\x) = |g^{(+)}_{\theta}(\x)|\ e^{\pm 2\pi e_{\nu_{\theta}}(\x)
\phi_{\theta}(\x)},
\label{polarform}
\end{equation}
where the amplitude, phase and direction are given by
\begin{eqnarray}
|g^{(+)}_{\theta}(\x)| &=& 
\sqrt{[g_{\theta}(\x)]^{2} + 
[g^{(1)}_{\theta;R}(\x)]^{2} + [g^{(2)}_{\theta;R}(\x)]^{2}},
\label{modulus}, \quad
e_{\nu_{\theta}}(\x) = \left[  \bm{i} \ \cos{\nu_{\theta}(\x)} +  \bm{j} \ \sin{\nu_{\theta}(\x)} \right], \\
\nu_{{\theta}}(\x) &=&  \tan^{-1} \left(\frac{  g^{(2)}_{\theta;R}(\x) }{ g^{(1)}_{\theta;R}(\x) }\right), \quad
\phi_{\theta} (\x)= \frac{1}{2\pi} \tan^{-1} \left({\mathrm{sgn}}\left( g^{(2)}_{\theta;R}(\x)\right)\frac{ \sqrt{[g^{(1)}_{\theta;R}(\x)]^{2}
+ [g^{(2)}_{\theta;R}(\x)]^{2}}}{g_{\theta}(\x)}\right). \nonumber 
\nonumber
\end{eqnarray}
\end{definition}
The ${\mathrm{sgn}}(\cdot)$ added to the definition of $\phi_{\theta}(\x)$
to appropriately
determine
the range of $\phi_{\theta}(\x).$
Note that the $\theta$-monogenic and $\theta$-anti-monogenic images of a
real image $g(\x)$
have
the same modulus and phase. $g_{\theta}^{(+)*}(\x)=g_{\theta}^{(-)}(\x).$
\begin{definition}{The $\theta$-polar representation of g(\x).}\\
The $\theta$-polar representation of a real signal $g(\x)$ is given by:
\begin{eqnarray}
\nonumber
g(\x)&=&\Re\left\{\widetilde{g}_{-\theta,\theta}^{(\pm)}(\x)\right\}=\left|\widetilde{g}_{-\theta,\theta}^{(\pm)}(\x)\right|
\cos(2\pi \widetilde{\phi}_{-\theta,\theta}(\x)).
\label{polarmonrep}
\end{eqnarray}
\end{definition}
Thus locally a signal is described by its magnitude 
$\left|\widetilde{g}_{-\theta,\theta}^{(+)}(\x)\right|$ and the local structure
of the signal is determined by $\widetilde{\phi}_{-\theta,\theta}(\x),$ from
which an instantaneous frequency is determined by $\left|\nabla\widetilde{\phi}_{-\theta,\theta}(\x)
\right|.$
From equation (\ref{notcommuating}) it clearly follows that:
$
\left|\widetilde{g}_{-\theta,\theta}^{(\pm)}(\x)\right|=\left|g^{(\pm)}(\x)\right|,$
$\widetilde{\phi}_{-\theta,\theta}(\x)=\widetilde{\phi}(\x),$ but $
e_{\widetilde{\nu}_{-\theta,\theta}}(\x)\neq e_{\nu}(\x).$
The monogenic representation of $g(\x)$ is thus equivalent
to the $\theta$-monogenic representation, unlike the case of the hypercomplex and 
$\theta$-hypercomplex representation, that are clearly distinct. 
The intuitive understanding of these relationships is that the monogenic
signal is constructing a plane wave representation at spatial point $\x,$ where the direction of the plane wave is
determined from the
signal. No matter what axes we construct, the local plane wave representation in the amplitude and period (phase) remain the same, but the parameterization of the direction changes with the axes. 
The reason for introducing the $\theta$-monogenic
signal is that the calculations of the properties of the CWT
are simplified.
Hence for convenience the notation $\left|\widetilde{g}_{-\theta,\theta}^{(\pm)}(\x)\right|$
and $\widetilde{\phi}_{-\theta,\theta}(\x),$ are retained, even if both are in some sense redundant.
\begin{proposition}{Orthogonality of the signal and the Riesz transform.}\\ A real image $g(\x)$  is orthogonal to each component of its Riesz transform ${\cal R}$. However, the two components of the 
Riesz transform, $g^{(1)}_R(\x)$, $g^{(2)}_R(\x)$ are in general not
orthogonal.
\end{proposition} 
\begin{proof}
This is a direct generalisation of the $n=1$ case where the HT
is orthogonal to the original signal. Indeed, exploiting the Hermitian symmetry of the FT of a real image
 $G^{*}(\f ) =G(-\f )$, we have
\begin{equation}
\langle g^{(1)}_R, g \rangle =  \int d^{2} \f \ j \ \frac{f_{1}}{f} \ G^{*}(\f ) \ G(\f ) =
\int d^{2} \f \ j \ \frac{f_{1}}{f} \ G(- \f ) \ G(\f ) = 0.
\label{<gig>}
\end{equation}
The function $G(- \f ) \ G(\f ) / f$ is even-even (ee) with respect to $f_{1}$ and $f_{2}$ whereas $f_{1}$ itself is
obviously odd-even (oe). Since the integration is over the entire $E_{2}$
plane, it is zero. By the same argument
$\langle g^{(2)}_R, g \rangle = 0$ and therefore $\langle {\cal R}g, g \rangle=0$.
\end{proof}
Furthermore we note that $\widetilde{g}_{-\theta}(\x)$
is orthogonal to $\widetilde{g}_{-\theta}^{(s)}(\x),$ and thus $g(\x)$ is orthogonal
to $\widetilde{g}_{-\theta}^{(s)}(\bm{r}_{-\theta}\x),\;s=1,\;2.$
However, the two components of the Riesz transform itself are not in general orthogonal: 
\begin{equation}
\langle g^{(2)}_R, g^{(1)}_R \rangle = \int  \frac{f_{1} f_{2}}{f^{2}}  G(\f)
G(-\f )d^{2} \f \neq 0.
\label{<gigj>}
\end{equation}
The natural way to view the monogenic image is as consisting of two orthogonal
components, $g(\x)$ and ${\cal R}g(\x)$. 
\begin{proposition}{Orthogonality of the components of the Riesz transform
of an isotropic function.}\\
If the signal $g(\x)$ is radially symmetric and real, then the integrand
in (\ref{<gigj>}) is
\begin{equation}
\langle g^{(2)}_R, g^{(1)}_R \rangle =  0.
\label{<gigj>22}
\end{equation}
\end{proposition}
\begin{proof}
$G(\f)=G(f)$, and radially symmetric as well as real, whilst $f_1 f_2 $ is
odd with respect to either $f_1$ or $f_2.$ Thus $\langle g^{(2)}_R, g^{(1)}_R
\rangle =0 $. Note that for an isotropic signal $g(\x)=\widetilde{g}_{-\theta}(\x)$
and so $g(\x)$ is orthogonal to $\widetilde{g}_{-\theta,H}^{(s)}(\x).$
\end{proof}
\begin{proposition}{Norm of the Riesz components.}\\
The norm of the Riesz components is given by (as noted by Felsberg \& Sommer \cite{FelsbergSommer}[p.~3140]):
\begin{equation}
\langle g^{(s)}_R, g^{(s)}_R\rangle =\int \frac{f_{s}^2}{f^{2}}  G(- \f )
G(\f )d^{2} \f ,\;s=1,2,\quad\langle g^{(1)}_R, g^{(1)}_R\rangle+\langle g^{(2)}_R, g^{(2)}_R\rangle
=\langle g,g \rangle.
\end{equation}
If furthermore $g(\x)$ is radially symmetric then
\begin{equation}
\langle g^{(s)}_R, g^{(s)}_R\rangle =\frac{1}{2}\langle g, g\rangle,\;s=1,2.
\end{equation}
\end{proposition}
These relationships were proved and used in \cite{Metikas} in the special
case of the isotropic multiple Morse wavelets. In summary ${\cal R}g(\x)$
is orthogonal to $g(\x),$ has the same frequency support, as well as the
same total norm as $g(\x)$ and exhibits similar spatial decay as $g(\x)$
(or at least as similar as ${\cal H}g(x_1)$ does to $g(x_1).$) Hence as ${\cal H}g(x_1)$ is considered as having the same time-frequency structure presence
as $g(x_1)$ (see for example Boashash \cite{Boashash1}) and we consider 
${\cal R}g(\x)$ as having the same spatial and spatial frequency support
as $g(\x).$

\subsection{The monogenic signal and the UQFT}
The monogenic signal determines the directional preference of a given real
signal at a given point in space, $\x.$ Let us discuss this aspect in some
more detail. As the monogenic representation is suitable for a plane wave
field $g(\x),$
consider the representation given
to a purely unidirectional signal $g_U(\x).$ By assuming that $g_U(\x)$ purely
experiences variation in a given direction, implies that it can be represented
as a superposition of plane waves with direction $\nu.$ Then for some given constant $\bm{n}=\left[\cos(\nu) \quad \sin(\nu) \right]:$
\begin{equation}
g_{U}(\x)=\int_{0}^{\infty} G_{U}(f)\cos(2\pi f\bm{n}^T \x)\; d f.
\end{equation}
Define the unit quaternion 
$\bm{e}_{\nu}=\cos(\nu)\bm{i}+\sin(\nu)\bm{j}.$ We calculate the monogenic
signal of $g_{U}(\x),$ and determine that:
\begin{eqnarray}
\label{directionalFT}
g_{U}^{(+)}(\x)&=&\int_{0}^{\infty} G_{U}(f)\exp(2\pi \bm{e}_{\nu}f\bm{n}^T \x) df,\quad
\nu=\frac{g_{U;R}^{(2)}(\x)}{g_{U;R}^{(1)}(\x)} \nonumber \\
&=&\int_{0}^{\infty} G_{U}(f)\cos(2\pi f\bm{n}^T \x) df+
\left[\bm{i}\cos(\nu)+\bm{j}\sin(\nu)\right]\int_{0}^{\infty} G_{U}(f)\sin(2\pi f\bm{n}^T \x) df.
\end{eqnarray}
This result follows from the monogenic extension of a sinusoid given in
\cite{FelsbergSommer}, and
the global orientation of the globally directional signal may be determined from the monogenic signal by calculating the ratio of the two Riesz transforms.
Let us now demonstrate how $g_{U}^{(+)}(\x)$ is related
to the UQFT in $\bm{e}_{\nu}$.
It may be shown, again using the monogenic extension of a sinusoid given
in
\cite{FelsbergSommer}, that
$
G_{U,\bm{e}_{\nu}}^{(+)}(\f)=G_U(f),\;f>0,
$
so that the spectral support of the UQFT of the signal in $\bm{e}_{\nu}$
is one-sided, and limited to the positive frequencies.  Thus an equivalent
to equation (\ref{HT}) can be derived
that corresponds to an alternative to $G_Q^{(++)}(\q)=0.$ The result is to be expected,
as the plane waves are really a 1-D feature embedded in a 2-D space. If the signal 
corresponds to
unidirectional variation only, the monogenic signal is extracting the orientation,
and lives in a half-plane, whilst the anti-monogenic signal lives in the
second half-plane. 
Calculating the partial HT in the correct
direction will yield an equivalent representation to using the monogenic
signal, but the advantage of using
the latter is that the direction is {\em
retrieved} directly
from the observed image without having to perform HTs at all
values of $\theta.$ The assumption of global unidirectionality
is very strict, and constructing a local description of the signal will be
essential for the performance of the analysis method.

\section{Wavelet analysis}
Wavelets provide a method of constructing a local decomposition of an observed
image, that for many
classes of images constitute of decomposition coefficients mainly of negligible
magnitude. Furthermore, due to the definition of the wavelet coefficients,
it is also frequently easier to interpret the
non-negligible coefficients,
rather than trying to disentangle the full behaviour
of an image in the original spatial domain. The interpretation of the coefficients
relies on the choice of mother wavelet.
To be
able to satisfy the reconstruction formula of equation (\ref{reconny5}) for
any arbitrary function $g(\x)\in L^2({\mathbb{R}}^2),$ the mother
wavelet function, denoted $\psi(\x),$ must satisfy the two conditions of
\begin{equation}
0 < c_{\psi} < \infty , \ {\rm where} \ c_{\psi} = (2 \pi)^{2} \int_{E_{2}}
\frac{ |\Psi(\f)|^{2} }{f^2}\ d^{2}\f,\quad \int_{E_2}\left| \psi(\x )\right|^2\;d^2\x=1.
\label{admissibility}
\end{equation}
Thus the mother wavelet is constrained to be both a spatially
local
function, and to be mainly supported over some range of frequencies not including
the
origin, {\em i.e.} the function is oscillatory and zero-mean. Without loss of generality,
center the mother wavelet function in space to $\x={\mathbf 0},$ and assume that
$\Psi(\f)$ is maximum at $\f=\f_0=f_0\left[\cos(\phi_0)\quad \sin(\phi_0)\right]^T.$
Of course the possible structure $\psi(\x)$ may exhibit whilst still satisfying
equation (\ref{admissibility}) is very variable, and even more so in 2-D
compared
to 1-D. For
example, the wavelets developed in 
\cite{Antoine1,Antoine3,Antoine4,GonnetTorresani} are
complex and describe
directional plane waves, as do the wavelets of \cite{Selesnick2005}, whilst a great
deal of applications that rely on perfect reconstruction use discrete filter
bank
separable wavelets
\cite{mallat}, and a third possible option corresponds to using radial
functions, as described in Antoine {\em et al.} \cite{AMVA}, Metikas \& Olhede
\cite{Metikas}
and also by Farge
\cite{Farge}. Visually, in a contour plot description of the wavelets,
highly directional plane wavelets have the appearance of local maxima corresponding to
evenly spaced lines with
a given normal, separable wavelets have local maxima and minima in a checkerboard patter, and isotropic
wavelets correspond to local maxima appearing in repeating ring-shapes with approximately an even spacing
between the rings. The `optimal' choice of wavelet naturally depends on the
application in question. 

This article aims to provide classes of mother wavelet functions so that the
wavelet coefficients calculated are hyperanalytic
signals in the index $\be.$ Additionally
energy
in the observed image should be separated to the correct decomposition coefficient,
and this corresponds to using
wavelets that are essentially supported over a limited region of space and
spatial frequency. 
To satisfy the latter restriction wavelets that
are well-localised \cite{AMVA,Metikas}, are often used.
Local structure will be represented
in terms of interrelationships between
coefficients in a vector-valued representation at each local decomposition
coefficient index point. To construct a four-vector
(quaternionic) description of local structure, we use wavelets
of the form:
\begin{equation}
\label{quatwaveletdef}
\psi(\x)=\psi^{(r)}(\x)+\bm{i}\psi^{(i)}(\x)+\bm{j}\psi^{(j)}(\x)+\bm{k}\psi^{(k)}(\x),
\end{equation}
where each $\psi^{(s)}\in L^2({\mathbb{R}}^2)$ for $s=r,\;i,\;j,\;k.$ Thus
for
any $\x\in{\mathbb{R}}^2$
it follows that
$\psi(\x)\in {\mathbb{H}}.$ 

In
general we construct $\psi(\x)$ from $\psi^{(r)}(\x)$ and will denote the construction
operator by $\psi(\x)={\cal W}\psi^{(r)}(\x).$ 
If the observed image can locally be described as a plane wave, as is the
case in equation (\ref{planeoscillation}),
then the local
monogenic description of the signal is appropriate, and
the CWT
will be constructed to correspond
to a local UQFT ({\em cf} equation (\ref{dirFourier}).)
If on the other hand,
local separable variations in an appropriate frame of reference
are present like the signal modelled by equation (\ref{separablexxx}),
then a local version of the 
QFT, {\em cf} equation (\ref{qftdefin}), should be used for analysis. 

Naturally the whole family of wavelet coefficients
constructed from the family of wavelet functions must enjoy suitable properties.
\begin{definition}{Members of the quaternionic wavelet family.}\\
An arbitrary member of the quaternionic wavelet family is defined by
\begin{equation}
\psi_{\bm{\xi}}( \x ) =  U_{\bm{\xi}} \ \psi(\x ) 
= \psi^{(r)}_{\bm{\xi}}(\x) + \bm{i}  \psi^{(i)}_{\bm{\xi}}(\x )
 + \bm{j} \psi^{(j)}_{\bm{\xi}}(\x) + \bm{k} \psi^{(k)}_{\bm{\xi}}(\x).
\label{TSRanalyticwavelet}
\end{equation}
\end{definition}
The construction corresponds to individually translating, scaling and rotating
the three components of the quaternionic function. To define any member of the quaternionic wavelet family, that is $\psi_{\bm{\xi}}(\x),$
each of the real wavelet functions
that combine to make up the quaternionic
wavelet is scaled, rotated and shifted so that $\psi_{\bm{\xi}}( \x ) =  U_{\bm{\xi}}\  {\cal W}\ \psi^{(r)}(\x ) .$ Note that in general $\psi_{\bm{\xi}}( \x ) \neq   {\cal W}\ U_{\bm{\xi}} \ \psi^{(r)}(\x ) ,$ and some care must
be taken when noting the properties of the members of the wavelet family.
Also the properties of the CWT coefficients must be considered with some
care, as in general
neither is
$w_{{\cal W} \psi^{(r)}}(\bm{\xi};g) = {\cal W} w_{ \psi^{(r)}}(\bm{\xi};g),$
where ${\cal W}$ is implemented in $\be.$
It is very important
that the wavelet
coefficients are
interpretable, as they may be considered as a local projection of the image,
and they, rather than any other given quantity, will be the basis for analysis
of a given image. We shall discuss the properties enjoyed by the wavelets as well as the CWT coefficients in detail.

\subsection{Hypercomplexing wavelets}
The QFT decomposes an image into sinusoids, where the parity of variations
in both $x_1$ and $x_2$ can be considered separately via the $\bm{i},$ $\bm{j}$
and $\bm{k}$ components of the QFT. 
From equation (\ref{qftdefin}) we can note that the construction of the QFT is highly non-trivial,
as $g(\x)$ is left-multiplied by $\exp(-2\pi \bm{i} x_1 q_1)$ and right-multiplied by
$\exp(-2\pi \bm{j} x_2 q_2),$ so if a quaternionic $g(\x)$ is decomposed
in the QFT coefficients, then the transform cannot be described in terms
of an inner product of $g(\x)$ with
another function. If $g(\x)\in {\mathbb{R}}$ then
$
G_Q (\q)=\int_{-\infty}^{\infty} \int_{-\infty}^{\infty} g(\x)
\exp(-2\pi \bm{i} x_1 q_1)\exp(-2\pi \bm{j} x_2 q_2)\\ \;d^2\x,$ and
the QFT may thus be expressed as an inner product between 
$\left[\exp(-2\pi \bm{i} x_1 q_1)\exp(-2\pi \bm{j} x_2 q_2)\right]^*$ and
$g(\x).$ Note that 
\begin{eqnarray}
 \left(e^{-2\pi \bm{i} q_1 x_1}e^{-2\pi \bm{j} q_2 x_2}\right)^*
&=& \vartheta^{(ee)}(\x)+\bm{i}\vartheta^{(oe)}(\x)+\bm{j}\vartheta^{(eo)}(\x)
-\bm{k}\vartheta^{(oo)}(\x)
\label{hyperbasis},
\end{eqnarray}
thus defining the functions $\left\{\vartheta^{(s)}(\x)\right\},$ $s=ee,\;oe,\;eo,\;oo.$
$\vartheta^{(ee)}(\x)=\cos(2\pi q_1 x_1)\cos(2\pi q_2 x_2)$ is an
even function in both $x_1$ and $x_2,$ $\vartheta^{(oe)}(\x)={\cal H}_1\left\{\vartheta^{(ee)}\right\}(\x),$
$\vartheta^{(eo)}(\x)={\cal H}_2\left\{\vartheta^{(ee)}\right\}(\x),$ whilst
$\vartheta^{(oo)}(\x)={\cal H}_1\left\{{\cal H}_2\left\{\vartheta^{(ee)}\right\}\right\}(\x).$
To duplicate the structure of $\left\{\vartheta^{(\cdot \cdot)}(\x)
\right\},$ constructing a local QFT, a real-valued mother wavelet function
$\psi^{(ee)}(\x)$ is first chosen
to be even in both $x_1$ and $x_2,$ but {\em not} necessarily separable.
We
then define
the {\em hypercomplexing wavelet}, by adding three extra functions to this
function, where each extra function multiplies a unit quaternion, and where
the functions are odd instead of even in either (or possibly both) of the two spatial directions. 

\begin{definition}{Hypercomplexing wavelets.}\\
The hypercomplexing wavelet function $\psi^{++}(\cdot)$ is formed from a
2-D real valued wavelet function $\psi^{(ee)}(\x)\in L^2({\mathbb{R}}),$
satisfying the admissibility condition of (\ref{admissibility}) and additionally the condition of:
\begin{equation}
\label{hyperadmiss}
\psi^{(ee)}(x_1,x_2)=\psi^{(ee)}(\pm x_1,\pm x_2)=\psi^{(ee)}(x_2,x_1),
\end{equation}
{\em i.e.} the real wavelet is even in both $x_1$ and $x_2$ and symmetric
across
the indices.
For such $\psi^{(ee)}(\x)$ let:
$\psi^{(oe)}(\x)={\cal H}_1\left\{\psi^{(ee)}\right\}(\x),$
$\psi^{(eo)}(\x)={\cal H}_2\left\{\psi^{(ee)}\right\}(\x),$
and
$\psi^{(oo)}(\x)={\cal H}_2\left\{{\cal H}_1\left\{\psi^{(ee)}\right\}\right\}(\x),$
and form the {\em hypercomplexing} mother wavelet via
\begin{eqnarray}
\psi^{++}(\x)
&=&\widetilde{{\cal HC}}\psi^{(ee)}\left(\x\right) =\psi^{(ee)}\left(\x\right)
+\bm{i}\psi^{(oe)}\left(\x\right)
+\bm{j}\psi^{(eo)}\left(\x\right)
-\bm{k}\psi^{(oo)}\left(\x\right).
\label{number}
\end{eqnarray}
\end{definition}
The hypercomplexing wavelet is defined by 
a set of three functions constructed from the partial HTs
of the even/even function in the two spatial directions, but note that ${\cal{HC}}
\neq \widetilde{{\cal{HC}}},$ (where ${\cal{HC}}$ is given by equation (\ref{hypercdef}).)
The Fourier
transforms of the functions forming $\psi^{++}(\x)$ are given by:
\begin{eqnarray}
&& \nonumber
\Psi^{(oe)}(\f)=(-\bm{j}) {\mathrm{sgn}}\left(f_1\right)\Psi^{(ee)}(\f),\quad
\Psi^{(eo)}(\f)=(-\bm{j}) {\mathrm{sgn}}\left(f_2\right)\Psi^{(ee)}(\f),\quad\\
&&\Psi^{(oo)}(\f)=- {\mathrm{sgn}}\left(f_1\right){\mathrm{sgn}}\left(f_2\right)\Psi^{(ee)}(\f).
\label{hypercc}
\end{eqnarray}
$\Psi^{(ee)}(\f)\in{\mathbb{R}},$ as $\psi^{(ee)}(\x)$ is even in both $x_1$ and
$x_2,$ and thus $\Psi_Q^{(ee)}(\q)\equiv \Psi^{(ee)}(\q).$ 
\begin{eqnarray}
\nonumber
\Psi^{++ *}(\f)&=&\left(1-\bm{k}{\mathrm{sgn}}\left(f_1\right)
+{\mathrm{sgn}}\left(f_2\right)+
\bm{k}{\mathrm{sgn}}\left(f_1\right){\mathrm{sgn}}\left(f_2\right)
\right)^*\Psi^{(ee)}(\f)\\
&=&\left(1+\bm{k}{\mathrm{sgn}}\left(f_1\right)
+{\mathrm{sgn}}\left(f_2\right)-
\bm{k}{\mathrm{sgn}}\left(f_1\right){\mathrm{sgn}}\left(f_2\right)
\right)\Psi^{(ee)}(\f).
\end{eqnarray}
The Fourier transform of the convolution of $\psi^{++*}(\x)$ with  $g(\x)\in
L^2\left({\mathbb{R}}^2\right),$ taking values in ${\mathbb{R}}$ is given by 
\[{\cal F}\left\{g \ast\ast \psi^{++ *}\right\}(\f)=
G(\f) \left(1+\bm{k}{\mathrm{sgn}}\left(f_1\right)
+{\mathrm{sgn}}\left(f_2\right)-
\bm{k}{\mathrm{sgn}}\left(f_1\right){\mathrm{sgn}}\left(f_2\right)
\right)\Psi^{(ee)}(\f).\]

\begin{proposition}{Admissibility of the hypercomplexing wavelets.}\\
The quaternionic function $\psi^{++}(\x),$ as well as any of its
four component functions, all satisfy the admissibility condition.
\end{proposition}
\begin{proof}
As $\psi^{(ee)}(x)$ is an admissible wavelet function, it follows all its
partial HTs are admissible, {\em i.e.} square integrable and satisfying
the admissibility condition given by equation (\ref{admissibility}).
\end{proof}

\begin{theorem}{Lack of hypercomplexity of $\psi^{++}(\x).$}\\
The function $\psi^{++}(\x)$ defined by equation (\ref{number}) is {\em not}
hypercomplex.
\end{theorem}
\begin{proof}
By direct calculation we note that the QFT of $\psi^{++}(\x)$ takes the form
\begin{eqnarray}
{\cal F}\left\{\psi^{++}\right\}(\q)&=&\left\{\begin{array}{lcr}
2\Psi^{(ee)}\left(\q\right)
&{\mathrm{if}}
& q_1>0\\\nonumber
2{\mathrm{sgn}}(q_2)\Psi^{(ee)}\left(\q\right) &{\mathrm{if}}
& q_1<0\end{array}\right.\\
&\neq & 0 \quad {\mathrm{if}}\quad q_1<0\quad {\mathrm{or}} \quad q_2<0 .
\label{qftofhyper}
\end{eqnarray}
Thus the QFT of $\psi^{++}(\x)$ is {\em not} wholly supported on positive quaternionic
frequencies only, and $\psi^{++}(\x)$ is {\em not} a hypercomplex signal. \end{proof}
The hypercomplexing wavelets are therefore denoted $\psi^{++}(\cdot),$
rather than
$\psi^{(++)}(\cdot),$ but the wavelet coefficients they construct are $w^{(++)}_{\psi}(\bm{\xi};g),$
as they are hypercomplex signals in $\be.$
Note that
Chan {\em et al.} \cite{Chan} define discrete quaternionic wavelet filters
by combining two real-valued discrete wavelet functions that are approximately
HT pairs, in a discrete version of equation (\ref{number}). The objective of the quaternion-valued analysis is to form representations
of signals that are hyperanalytic signals rather than using hyperanalytic
signals in constructing the decomposition, and as we shall see, the hyperanalytizing
wavelets are useful decomposition filters.
%

\begin{proposition}{FT and QFT of the members
of the hypercomplexing wavelet family.}\\
Any member of the hypercomplexing wavelet family defined by substituting
in the form of the mother wavelet into equation (\ref{TSRanalyticwavelet})
has a FT and a QFT given by
\begin{eqnarray}
\nonumber
\Psi^{++}_{\bm{\xi}}(\f)&=&\left|a\right|
\left[1+{\mathrm{sgn}}\left(\left[\bm{r}_{-\theta}\f\right]_2\right)
-\bm{k}{\mathrm{sgn}}\left(\left[\bm{r}_{-\theta}\f\right]_1\right)
\left[1-
{\mathrm{sgn}}\left(\left[\bm{r}_{-\theta}\f\right]_2\right)\right]
\right]\Psi^{(ee)}(a\bm{r}_{-\theta}\f)e^{-2\pi\bm{j} \q^T \be}
\\
\nonumber
\Psi^{++}_{\bm{\xi},Q}(\q) &=& \left|a\right| e^{-2\pi\bm{i} q_1 b_1}
{\cal F}_{\cal Q}\left\{\psi^{++}_{(1,\theta,\bm{0})}\right\}
(a\q)e^{-2\pi\bm{j} q_2 b_2}\\
\nonumber
&=& \left|a\right| e^{-2\pi\bm{i} q_1 b_1}
\frac{1-\bm{k}}{2}\Psi^{(ee)}(a\bm{r}_{-\theta}\q)\left(1+{\mathrm{sgn}}\left(
\left[\bm{r}_{-\theta}\q\right]_{1} \right)+{\mathrm{sgn}}\left(
\left[\bm{r}_{-\theta}\q
\right]_{2} \right)\right.\\
\nonumber
&&\left.-{\mathrm{sgn}}\left(
\left[\bm{r}_{-\theta}\q
\right]_{1} \right){\mathrm{sgn}}\left(
\left[\bm{r}_{-\theta}\q
\right]_{2} \right)\right)
+\frac{1+\bm{k}}{2}\Psi^{(ee)}(a\bm{r}_{-\theta}\bm{J}_{\pi/2}\q)\left(1+{\mathrm{sgn}}\left(
\left[\bm{r}_{\theta}\q\right]_{1} \right)\right.\\
&&\left.+{\mathrm{sgn}}\left(
\left[\bm{r}_{\theta}\q
\right]_{2} \right)-{\mathrm{sgn}}\left(
\left[\bm{r}_{\theta}\q
\right]_{1} \right){\mathrm{sgn}}\left(
\left[\bm{r}_{\theta}\q
\right]_{2} \right)\right)e^{-2\pi\bm{j} q_2 b_2}
\label{hypercompqft}
\end{eqnarray}
respectively. 
\end{proposition}
\begin{proof}
The Fourier transform follows by direct calculation, whilst substituting
in the form of the mother wavelet into equation (\ref{qftofhyper}) yields
the expression for the QFT.
\end{proof}
Given the mother hypercomplexing wavelet was not a hypercomplex signal it
is not surprising that neither will an arbitrary member of the family be
a hypercomplex signal, as we may note from equation (\ref{hypercompqft}).
We will build hypercomplexing wavelets from either 
$\psi^{(ee)}(\x)=\psi^{(ee)}_S(\x),$ the {\em separable} mother wavelet corresponding
to a tensor product
of two 1-D functions, or the {\em isotropic} real mother wavelet $\psi^{(ee)}(\x)=\psi^{(ee)}_I(x),$
that is a function of $x=\sqrt{x_1^2+x_2^2}$. The reasoning behind
these choices, is that we replace $\vartheta^{(ee)}(\x)$ by functions
that are {\em even} in both spatial directions, and that are chosen to satisfy equation
(\ref{hyperadmiss}). Either of $\psi^{++}_I(\x)$ and $\psi^{++}_S(\x)$
will enable us
to decompose the observed image into local $\theta$-hypercomplex signals,
for any given
values of $a$ and $\theta.$

\subsubsection{Separable hypercomplexing wavelets \label{sepwavs}}
The separable hypercomplexing wavelets are constructed from a single 1-D
even function, denoted $\psi^{(e)}_{1D}(x_1).$ 
$\psi^{++}_S(\x)$
is constructed to filter the signal to contributions corresponding
to the same period in both spatial directions -- by augmenting the doubly
even function with the three extra wavelet functions we construct directional
wavelets that are well localised in both direction and scale as well as spatial
position. Using $\psi^{++}_S(\x)$ we separate components according to their
directional localisation,
and the CWT needs to be calculate across the full range of $\theta\in\left[0,\pi/2\right).$
First define the doubly
even mother wavelet that will serve as the real component of the quaternionic
function $\psi^{++}(\x),$ by $\psi^{(ee)}(\x)=\psi^{(e)}_{1D}(x_1)\psi^{(e)}_{1D}(x_2).$
We define the 1-D function
$\psi^{(o)}_{1D}(x)={\cal H}\left\{\psi^{(e)}_{1D}\right\}(x).$ We may then observe
that using equation (\ref{hypercc}) in this instance we find:
\[\psi^{(oe)}_S(\x)=\psi^{(o)}_{1D}(x_1)\psi^{(e)}_{1D}(x_2),\quad
\psi^{(oe)}_S(\x)=\psi_{1D}^{(o)}(x_1)\psi^{(e)}_{1D}(x_2),\quad
\psi^{(oo)}_S(\x)=\psi^{(o)}_{1D}(x_1)\psi^{(o)}_{1D}(x_2).
\]
Assume $f_0=\arg_{f>0} \max \Psi_{1D}^{(e)}(f)$. The wavelet functions may
be considered a filtering the signal to frequencies near $f_1,\;f_2\approx
\pm f_0.$ 
Note that in 1-D analytic and anti-analytic wavelets are defined by combining
$\psi^{(e)}_{1D}(x_1)$ and $\psi^{(o)}_{1D}(x_1)$ into a complex-valued wavelet
function by:
\begin{equation}
\label{analdef}
\psi^{(\pm)}_{1D}(x_1)=\psi^{(e)}_{1D}\left(x_1\right)
\pm\bm{j}\psi^{(o)}_{1D}\left(x_1\right)\quad
\Psi^{(\pm)}_{1D}(f_1)=\Psi^{(e)}_{1D}\left(f_1\right)
\pm\bm{j}\Psi^{(o)}_{1D}\left(f_1\right)\in{\mathbb{R}}.
\end{equation}
Thus the FT,
and the QFT of the members of the separable hypercomplexing wavelet family
are given by:
\begin{eqnarray}
\Psi^{++}_{S,\bm{\xi}}(\f)&=&\left|a\right|
\left[\Psi^{(e)}\left(a\left[\bm{r}_{-\theta}\f\right]_1\right)
\Psi^{(+)}\left(a\left[\bm{r}_{-\theta}\f\right]_2\right)+
\bm{i}\Psi^{(o)}\left(a\left[\bm{r}_{-\theta}\f\right]_1\right)
\Psi^{(-)}\left(a\left[\bm{r}_{-\theta}\f\right]_2\right)
\right]e^{-2\pi\bm{j}\f^T \be},
\end{eqnarray}
\begin{eqnarray}
\Psi^{++}_{S,\bm{\xi};Q}(\q)&=&\left|a\right| e^{-2\pi\bm{i} q_1 b_1}
\frac{1-\bm{k}}{2}\left(\Psi^{(+)}(a\left[\bm{r}_{-\theta}\q\right]_1)
\Psi^{(e)}(a\left[\bm{r}_{-\theta}\q\right]_2)+
\Psi^{(-)}(a\left[\bm{r}_{-\theta}\q\right]_1)
\Psi^{(o)}(a\left[\bm{r}_{-\theta}\q\right]_2)\right)
\nonumber\\ \nonumber
&&+\frac{1+\bm{k}}{2}\left(\Psi^{(+)}(a\left[\bm{r}_{-\theta}\bm{J}_{\pi/2}\q\right]_1)
\Psi^{(e)}(a\left[\bm{r}_{-\theta}\bm{J}_{\pi/2}\q\right]_2)\right.\\
&&\left.+
\Psi^{(-)}(a\left[\bm{r}_{-\theta}\bm{J}_{\pi/2}\q\right]_1)
\Psi^{(o)}(a\left[\bm{r}_{-\theta}\bm{J}_{\pi/2}\q\right]_2)
\right)e^{-2\pi\bm{j} q_2 b_2},
\nonumber
\end{eqnarray}
respectively. The interpretation of $\Psi^{++}_{S,\bm{\xi}}(\f)$ is that
the hypercomplexing wavelets are naturally treated in a rotated frame of reference,
whilst from the form of the QFT it is apparent that
there is no choice of rotation that will make the wavelet $\theta$-hypercomplex.
However, as the mother failed to be hypercomplex, this is not surprising.
\begin{figure*}[t]
\centerline{
\includegraphics[height=1.55in,width=1.55in]{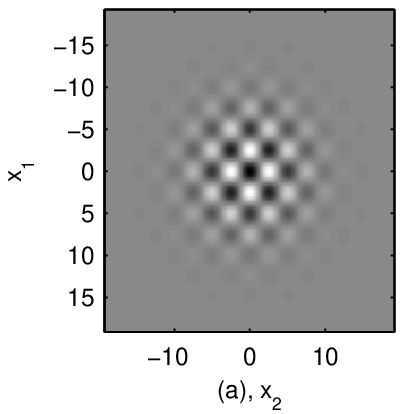}
\includegraphics[height=1.55in,width=1.55in]{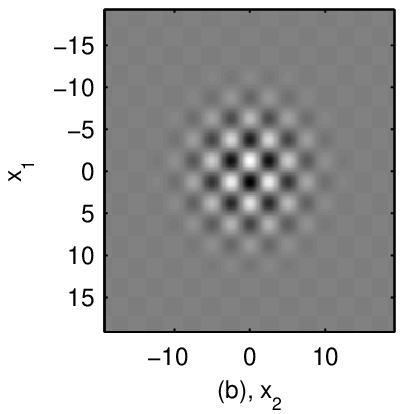}
\includegraphics[height=1.55in,width=1.55in]{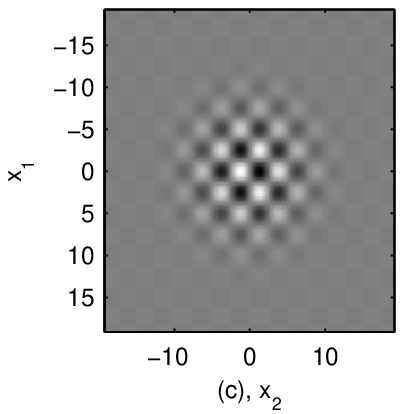}
\includegraphics[height=1.55in,width=1.55in]{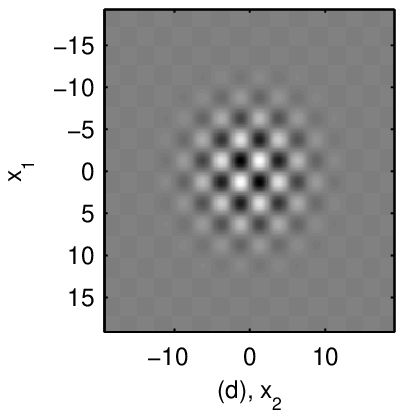}}
\caption{\label{figure1} 
The real part ((a)), first, second and third imaginary
parts ((b), (c) and (d)) of the separable hypercomplexing Morse wavelets,
with
$n_1=n_2=0,$ $\beta=9$ and $\gamma=4.$
}
\end{figure*}

Using the separable hypercomplexing wavelets each of the two
analytic wavelet functions may be given a polar representation, {\em i.e.}\\
$\psi^{(+)}_{1D}(x_1)=\left|\psi^{(+)}_{1D}(x_1)
\right|e^{2\pi \bm{j} \phi_{\psi,1D}^{(+)}(x)},
$
and thus it transpires that
$\psi^{++}_S(\x)=\left|\psi^{(+)}_{1D}(x_1)\right|\left|\psi^{(+)}_{1D}(x_2)\right|
e^{2\pi \bm{i} \phi_{\psi,1D}^{(+)}(x_1)}\\e^{2\pi \bm{k}0} e^{2\pi \bm{j} \phi_{\psi,1D}^{(+)}(x_2)}$ or $\psi^{++}_S(\x)=\left|\psi^{++}_S(\x)\right|
e^{2\pi \bm{i} \alpha_{\psi}(\x)}$ $e^{2\pi \bm{k}\gamma_{\psi}(\x)} e^{2\pi \bm{j} \beta_{\psi}(\x)},$
with $\left|\psi^{++}_S(\x)\right|=\left|\psi^{(+)}_{1D}(x_1)\right|\left|\psi^{(+)}_{1D}(x_2)\right|,$
$\alpha_{\psi}(\x)=\phi_{\psi,1D}^{(+)}(x_1),$ $\beta_{\psi}(\x)=\phi_{\psi,1D}^{(+)}(x_2),$
and
$\gamma_{\psi}(\x)=0.$
The separable hypercomplexing wavelets are representing local structure by finding the local structural information in $x_1$
and $x_2$ separately, and representing this in terms of a magnitude and two
phases. $\phi_{\psi,1D}^{(+)}(x)$ is the local frequency
content of the 1-D function, $\psi^{(+)}_{1D}(x_1),$ and in the two separate
spatial directions
the same frequency behaviour is considered. This is obvious from from Figure
\ref{figure1} where the wavelets correspond to two separable oscillations
local to $\x={\mathbf{0}},$ with the same local period in both spatial directions.
 
\subsubsection{Isotropic hypercomplexing wavelets \label{isowavs}}
Basing the hypercomplexing wavelet on an isotropic
mother wavelet, $\psi^{(ee)}_I(\x),$ 
wavelet coefficients are
produced by filtering the signal to retain contributions with the same
period. Using the extra components 
potentially separable structures can be characterised that have the same
period present at a given
spatial location.
Figure \ref{figure2} shows isotropic hypercomplexing wavelets, and how the
spatial localisation is radial, whilst each added component is odd in $x_1$
and $x_2,$ respectively. 
When an isotropic mother wavelet is used then we obtain that:
\begin{equation}
\psi_I^{(oe)}(\x)={\cal H}_1\left\{\psi_I^{(ee)}\right\}(\x),\quad
\psi_I^{(eo)}(\x)={\cal H}_2\left\{\psi_I^{(ee)}\right\}(\x),\quad
\psi_I^{(oo)}(\x)={\cal H}_T\left\{\psi_I^{(ee)}\right\}(\x).
\end{equation}
These wavelets are able to determine any preference for variation in the
direction $\chi=\tan^{-1}\left(x_2/x_1\right)=0,$ $\chi=\pi/2$ and $\chi=\pi/4,$ respectively. This is
clear from the example of an isotropic hypercomplexing wavelet given
in Figure \ref{figure2}.

\subsubsection{Hypercomplex wavelets \label{dirhyper}}
Following work on discrete wavelet filters given in \cite{Selesnick2005}[p.~138] a real directional wavelet
based on the 1-D analytic continuous wavelet
filters of section \ref{sepwavs} will be constructed via first defining the
$\pi/4$ directional wavelet function:
\begin{equation}
\label{dirryyy}
\psi^{(\pi/4)}_D(\x)=\psi^{(e)}_{1D}(x_1)\psi^{(e)}_{1D}(x_2)-
\psi^{(o)}_{1D}(x_1)\psi^{(o)}_{1D}(x_2).
\end{equation}
$\psi_D^{(\pi/4)}(\x)$ is local in frequency to $\f=\pm\sqrt{2}f_0\left[\cos(\pi/4)\; \sin(\pi/4)\right]^T.$ Formally, the above wavelet, does not satisfy (\ref{hyperadmiss}), and cannot
determine existing separable 
structure as it only considers structure along the diagonal direction and
not variation perpendicular to this direction.
To align the wavelet with
the axis of observation we define
\begin{equation}
\label{dirryyy2}
\psi_D(\x)=\psi^{(\pi/4)}_D(\bm{r}_{\pi/4}\x)=
\psi^{(e)}_{1D}\left(\frac{x_1-x_2}{\sqrt{2}}\right)\psi^{(e)}_{1D}\left(\frac{x_1+x_2}{\sqrt{2}}\right)-
\psi^{(o)}_{1D}\left(\frac{x_1-x_2}{\sqrt{2}}\right)\psi^{(o)}_{1D}\left(\frac{x_1+x_2}{\sqrt{2}}\right).
\end{equation}
Thus $\Psi_D(\x)$ will be local in frequency to $\f=\pm \sqrt{2}f_0\left[1\; 0\right]^T.$
The function is extended to a quaternionic object by calculating its
total and partial HTs:
\begin{eqnarray*}
{\cal H}_1\left\{\psi_D\right\}(\x)&=&
\psi^{(o)}_{1D}\left(\frac{x_1-x_2}{\sqrt{2}}\right)\psi^{(e)}_{1D}\left(\frac{x_1+x_2}{\sqrt{2}}\right)-
\psi^{(e)}_{1D}\left(\frac{x_1-x_2}{\sqrt{2}}\right)\psi^{(o)}_{1D}\left(\frac{x_1+x_2}{\sqrt{2}}\right)=
\psi_{D,2}(\x), \\
{\cal H}_2\left\{\psi_D\right\}(\x)&=&0, \quad
{\cal H}_T\left\{\psi_D\right\}(\x)=0.
\end{eqnarray*}
Formally by direct calculation:
$
\psi^{(++)}_D(\x)=\psi_D(\x)+\bm{i}\psi_{D,2}(\x)
.$ $
\psi^{(++)}_D(\x)$
is a directional quaternionic wavelet function, naturally taking the
form of a wavelet suggested by \cite{Selesnick2005}. This wavelet is formed by a similar set of operations 
to the hypercomplexing wavelet but we have
given it the notation $\psi^{(++)}_D(\x),$ as formally, it is also a hypercomplex
signal. As the third and fourth real-valued components of the quaternionic
objects are equivalently zero, and it is the non-zero fourth component, which
prevents the hypercomplexing wavelet from being a hypercomplex signal. $\psi_D^{(++)}(\x)$
can be treated as a complex-valued wavelet
(see footnote \ref{footie}), but
is not equivalent to the complex directional cone wavelets of Antoine {\em et
al.} \cite{Antoine3}. The two components of $\psi^{(++)}_D(\x)$ are plotted
in Figure \ref{figure4}.

\subsection{Monogenic wavelets}
The hypercomplexing wavelets were constructed so they filter the image
as a local version of the QFT in orientation $\theta,$ and their usage will
(as we shall see)
produce wavelet
coefficients that are $\theta$-hypercomplex signals.  
The monogenic signal is locally
identifying a dominant orientation of variations
and constructing, by adding
the Riesz transforms of the components to the signal, a polar representation
of such variation. The monogenic signal
of a globally unidirectional signal could be represented as an Inverse UQFT
(IUQFT) in $\bm{e}_{\nu},$ constructed only
from the positive frequencies of the UQFT in $\bm{e}_{\nu},$
and had an interpretation from
equation (\ref{dirFourier}). The direction of a unidirectional signal could
be determined directly by
the ratio of the Riesz transform components. A scale-localised version
of this construction would correspond to filtering the image in
space and spatial frequency, thus isolating structure with a given period
and spatial location, and then constructing the  Riesz transforms of the local component. The resulting quaternionic function could be represented as a quaternionic
object. The route to constructing the scale localised monogenic signal in
one step will be using quaternionic mother wavelet functions.

\begin{definition}{Monogenic wavelets.}\\
The monogenic extension of any given real-valued mother wavelet $\psi^{(r)}(\x)$
is given by
\begin{equation}
\psi^{(+)}(\x ) = {\cal M}\psi^{(r)}(\x)=\psi^{(r)}(\x) + {\cal R}\psi^{(r)}(\x) = \psi^{(r)}
(\x ) + \left[ \bm{i} \ \psi^{(1)}_R(\x )
+  \bm{j} \ \psi^{(2)}_R(\x ) \right],
\label{analyticwavelet}
\end{equation}
where $\psi^{(1)}_R = {\cal R}_{1} \ \psi^{(r)}$, $\psi^{(2)}_R = {\cal R}_{2} \ \psi^{(r)}.$ 
\end{definition}
Note that the terminology `monogenic wavelet' has been utilized before in
a different context: see work by Cerejeiras {\em et al.} \cite{cerejeiras}.
The monogenic wavelets
of \cite{cerejeiras}
are defined on the unit ball and based on M{\"o}bius transformations, and
are distinct from our wavelets.

\begin{lemma}{Properties of the wavelet Riesz components.}\\ 
Each component of the Riesz transform of a wavelet is also a wavelet.
\end{lemma}
\begin{proof}
As 
\begin{eqnarray*}
\int_{E_{2}} \psi^{(s)2}_R(\x)\;d^2\x=\int_{E_{2}} \left|\Psi^{(s)}_R(\f)\right|^2\;d^2\f
\le\int_{E_{2}} \left|\Psi^{(e)}(\f)\right|^2\;d^2\f<\infty,\;s=1,2,
\end{eqnarray*}
square integrability follows as the wavelets were constructed from a real
mother
wavelet.
Relations (\ref{riesssy5}), (\ref{weirdeqn}), hold of course, for $\psi^{(+)}$ and its components
$\psi^{(r)}$, $\psi^{(1)}_R$, and $\psi^{(2)}_R$. Denote by $c_{\psi^{(s)}}$
the value of equation (\ref{admissibility}) with $\psi^{(s)}_R.$
Equation (\ref{riesssy5}) in combination with the fact that $\cos^2{\phi}$ and
$\sin^{2}{\phi}$ are always between 0 and 1, imply that $0 < c_{\psi^{(1)}} \leq c_{\psi^{(r)}} < \infty $ and $ 0 
< c_{ \psi^{(2)} } \leq c_{\psi^{(r)}} < \infty $. Consequently $\psi^{(1)}_R$ and $\psi^{(2)}_R$ satisfy the admissibility condition, (\ref{admissibility}), and can be considered as wavelets. 
\end{proof}
The lemma implies that by the monogenic extension of the real wavelet we
have constructed two additional well-behaved functions and we may then consider
the properties of the monogenic extension as a whole.

\begin{proposition}{Admissibility of the monogenic mother wavelet.}\\
The monogenic extension of a real wavelet is also a wavelet.
\end{proposition}
\begin{proof}
It is obvious from (\ref{analyticwavelet}) and (\ref{weirdeqn}) that 
$|\Psi^{(+)}(\f)|^{2} = |\Psi^{(r)}(\f) |^{2} + 
|\Psi^{(1)}_R(\f)|^{2}  + 
|\Psi^{(2)}_R(\f)|^{2} + 2   \ \sin{\phi} \  |\Psi^{(r)}(\f)|^{2}$. Integrating over all frequencies as in
(\ref{admissibility}) and making use of the Hermitian property of the FT of a real image
 we find that 
$ c_{\psi^{(+)}} = c_{\psi^{(r)}} + c_{\psi^{(1)}} + c_{\psi^{(2)}}. $
Therefore, since each of $\psi^{(r)}$, $\psi^{(1)}_R$, and $\psi^{(2)}_R$
satisfy the admissibility condition, $\psi^{(+)}$ also satisfies the admissibility
condition and can be considered as a monogenic wavelet. 
Moreover, $c_{\psi^{(1)}} + c_{\psi^{(2)}} = c_{\psi^{(r)}}$, so
$c_{\psi^{(+)}} = 2 \ c_{\psi^{(r)}}$ also holds.
\end{proof}
Hence the monogenic extension of the wavelet is a valid analysis wavelet
and we may consider its further properties.


The monogenic wavelets may be represented in polar form
via its modulus, orientation, and phase respectively given by:
\begin{eqnarray}
\left| \psi^{(+)}_{\bm{\xi}}\left(\x\right) \right|
 &=& \sqrt{ [\psi^{(r)}_{\bm{\xi}}(\x)]^2 +  [\psi^{(1)}_{R,\bm{\xi}}(\x)]^2 +  
 [\psi^{(2)}_{R,\bm{\xi}}(\x)]^2},
\label{TSRmodulus}\\
\nonumber
{\bm{e}}_{\nu_{\bm{\xi}}}(\x)&=&  \left[ \bm{i} \ \frac{  
\psi^{(1)}_{R,\bm{\xi}}(\x) }{\sqrt{[\psi^{(1)}_{R,\bm{\xi}}(\x)]^2 + 
[\psi^{(2)}_{R,\bm{\xi}}(\x)]^2}}  + \bm{j} \ \frac{ \psi^{(2)}_{R,\bm{\xi}}(\x)
}{[\psi^{(1)}_{R,\bm{\xi}}(\x)]^2 +  [\psi^{(2)}_{R,\bm{\xi}}(\x)]^2}  \right]\\
&=&\bm{i}\cos(\nu_{{\xi}}(\x))+\bm{j}\sin(\nu_{{\xi}}(\x)),
\label{TSRunit}
\end{eqnarray}
\begin{eqnarray}
\phi_{{\bm{\xi}}}(\x) &= & \tan^{-1}\left(
{\rm sgn}\left[  \psi^{(2)}_{R,\bm{\xi}}\left(\x\right) \right]
\frac{\sqrt{[\psi^{(1)}_{R,\bm{\xi}}(\x)]^2 +
[\psi^{(2)}_{R,\bm{\xi}}(\x)]^2}}{ \psi_{\bm{\xi}}^{(r)}(\x)}\right).
\label{TSRphase}
\end{eqnarray}
The orientation is thus given an angular representation. The wavelet
is a function with spatial energy given by
the modulus, has a local zero-crossing structure given 
by $\phi_{\psi_{\xi}}(\x),$
and is showing a local orientation preference to angle $\nu_{\psi_{\xi}}(\x).$
The CWT is decomposing the image in terms of localised
oscillations with a local period determined from 
$\phi_{\psi_{\xi}}(\x).$ For a more
thorough discussion of the interpretation of the polar representation of
$\psi_{\bm{\xi}}^{(+)}(\x),$ see \cite{Spie}.

\begin{proposition}{FT and QFT of the members monogenic wavelet family.}\\
The FT and QFT of the translated, dilated, and rotated monogenic wavelet is
\begin{eqnarray}
&& \Psi^{(+)} _{\bm{\xi}}( \f ) =  {\cal F} \ \psi^{(+)}_{\xi}(\x) = \Psi^{(r)}_{\bm{\xi}}(\f) +    \  \left[  \bm{i} \   
\Psi^{(1)}_{R,\bm{\xi}}(\f)
 + \bm{j} \  \Psi^{(2)}_{R,\bm{\xi}}(\f)  \right] \nonumber \\ 
&& = \left\{ 1 +  \ \left[ - \bm{k} \cos{(\phi -  \theta)} +  \sin{(\phi -  \theta)} \right]
\right\} \ \Psi_{\bm{\xi}}^{(r)}(\f).
\label{FTanalyticwavelet}
\end{eqnarray}
\begin{eqnarray}
\nonumber
\Psi^{(+)}_{\bm{\xi};Q}( \q) &=&  {\cal F}_{{\cal Q}} \ \psi^{(+)}_{\bm{\xi}}(\x)
 = \Psi^{(r)}_{\bm{\xi};Q}(\q)  +  \ \left[ \bm{i} \ \Psi^{(1)}_{R,\bm{\xi};Q}( \q)  +
\Psi^{(2)}_{R,\bm{\xi}; Q}( \q) \ \bm{j} \right]  \nonumber \\
\nonumber
&=&\frac{a}{2}e^{-2\bm{i}\pi q_1 b_1}
\left[\left[1+\cos(\eta-\theta)+\sin(\eta-\theta)
\right]\left[\Psi^{(r)}_Q(\bm{r}_{-\theta}\q)-\bm{k}
\Psi^{(r)}_Q(\bm{J}_{\pi/2+\theta/2}\q)\right]\right.
\\&&\left.+\left[1+\cos(\eta+\theta)+\sin(\eta+\theta)
\right]\left[\Psi^{(r)}_Q(\bm{r}_{\theta}\q)+\bm{k}
\Psi^{(r)}_Q(\bm{J}_{\pi/2-\theta/2}\q)\right]\right]
e^{-2\bm{j}\pi q_2 b_2}.
\label{QFTanalyticwavelet}
\end{eqnarray}
\end{proposition}
\begin{proof}
See section \ref{ftqftmon}.
\end{proof}
As noted by Holschneider \cite{Holschneider} the entire family of wavelets
constructed from a 1-D analytic mother wavelet are analytic signals. 
Analogously, for the monogenic wavelets, we may state the result:

\begin{theorem}{Monogenicity of the members of the wavelet family.}\\
The members of the family of equation (\ref{TSRanalyticwavelet}) generated from the monogenic
wavelet $\psi^{(+)}(\cdot)$ through $U_{\bm{\xi}}$ are also $\theta$-monogenic
signals, this
should be contrasted with the $\theta$-hypercomplexing wavelets, that are
not $\theta$-hypercomplex signals.
\end{theorem}
\begin{proof}
This is not a trivial result as we noted in equation (\ref{notcommuating}),
and corresponds to more than just a simple 2-D extension of the 1-D result.
An arbitrary member of the wavelet family may be recast as $U_{\bm{\xi}}
\  \psi^{(r)}(\x) + \bm{i} \ U_{\bm{\xi}} \ {\cal R}_{1} \ \psi^{(r)}(\x)
+ 
\bm{j} \ U_{\bm{\xi}} \ {\cal R}_{2} \ \psi^{(r)}(\x)  $ as the Riesz kernels
commute with dilations and translations, $[{\cal D}_{a}, 
{\cal R}_{j}] = 0$, $[{\cal T}_{\be}, {\cal R}_{j}] = 0.$ Of course $R_{\theta}\psi^{(r)}(\x)$
is a $\theta$-monogenic signal/function from equation (\ref{thetaanalytic}), and thus
from
proposition \ref{thetamonny2}, we may note the result.
\end{proof}
Denote by $\psi^{(e)}(x)$ an isotropic real mother wavelet function.
\subsubsection{Isotropic monogenic wavelets}
\begin{proposition}{The wavelet family of the isotropic monogenic mother
wavelet.}\\
If  $\psi^{(r)}(\x)=\psi^{(e)}(x)$, then any member of the monogenic wavelet
family in the Fourier and Quaternionic Fourier
domain have a simplified form with $\bm{\xi}_0=\left[a,0;\be\right]^T$ by:
\begin{eqnarray}
\nonumber
\psi^{(+)}_{\bm{\xi}}(\x) &=& \psi^{(e)}_{\bm{\xi}_0}(\x) +\bm{i}\left[\cos(\theta)\psi^{(1)}_{R,\bm{\xi}_0}(\x)
+\sin(\theta)\psi^{(2)}_{R,\bm{\xi}_0}(\x)\right]+\bm{j}\left[-\sin(\theta)\psi^{(1)}_{R,\bm{\xi}_0}(\x)
+\cos(\theta)\psi^{(2)}_{\bm{\xi}_0}(\x)\right]\\
\label{isomon1}
\Psi^{(+)}_{\bm{\xi}}( \f ) 
&=& \left\{ 1 + \ \left[ - \bm{k} \cos{(\phi -  \theta)} +  \sin{(\phi -  \theta)} \right]
\right\} \ a \Psi^{(e)}(af) \ e^{-j 2\pi \f \be}, \\
\Psi^{(+) }_{\bm{\xi};Q}( \q )                                           &=& 
                                             e^{-2\pi\bm{i}q_1b_1}
                                             a\Psi^{(e)}(aq)
\left[1+\cos(\theta)\left(\cos(\eta)+\sin(\eta)\right)
+\bm{k}\sin(\theta)\left(-\sin(\eta)+\cos(\eta)\right)
\right]
e^{-2\pi\bm{j}q_2b_2}.
\nonumber
\end{eqnarray}
\end{proposition}
\begin{proof}
See section \ref{isofamily}.
\end{proof}
As an example of the isotropic monogenic wavelet, see Figure \ref{figure3}.
The isotropic structure of the magnitude is clear, whilst variations over
the $x_1$ and $x_2$ directions may respectively be determined.

\begin{proposition}{$\theta=\pi$ Equivalence to the anti-monogenic wavelet.}\\
When the monogenic mother has a radially symmetric real part, then the $\bm{\xi}_{\pi}$
monogenic wavelet coincides with the $\bm{\xi}$ anti-monogenic wavelet, and
with $\bm{\xi}_{\pi}=\left[a,\theta+\pi;\be\right]^T,$
\begin{equation}
\label{equivalentpirot}
\psi^{(+)}_{\bm{\xi}_{\pi}}(\x) = \psi^{(-)}_{\bm{\xi}}(\x).
\end{equation}
\end{proposition}
\begin{proof}
This we can see by replacing $\theta$ by
$\theta + \pi$ in  (\ref{isomon1}).
\end{proof}

\subsubsection{Directional monogenic wavelets \label{dirmon}}
The isotropic monogenic wavelets will localise an image to a particular scale
at each spatial position, and then represent the structure present at that
scale and position by a plane wave. In some applications it may be suitable to localise
in both scale and orientation, as several components are present at each
position and
at the same scale in different directions, and we start by constructing a
quaternionic wavelet starting with a real directional
wavelet
defined by equation (\ref{dirryyy2}). The Riesz transforms are found in
Appendix \ref{dirrictly}, and we obtain:
\begin{eqnarray*}
{\cal R}_1\left\{\psi_D\right\}(\x)&=&
4\int_{0}^{\infty} \int_0^{\infty}\frac{\sqrt{2} g_1}{g}\left[\Psi^{(e)}\left(g_1\right)
\Psi^{(e)}\left(g_2\right)\right] 
\sin(2\pi x_1(\frac{g_2+g_1}{\sqrt{2}}))
\;d^2\g=\psi_{D,3}(\x)\\
{\cal R}_2\left\{\psi_D\right\}(\x)&=& 0.
\end{eqnarray*}
The complex directional wavelet (see footnote \ref{footie} on treating $\psi_D^{(+)}(\x)$
as complex-valued) is formed by:
$
\psi^{(+)}_D(\x)=\psi_D(\x)+\bm{i}
\psi_{D,3}(\x).$ $\psi^{(+)}_D(\x)$ is 
a directional monogenic signal. $\psi^{(+)}_D(\x)$'s
directional structure may clearly be observed from a contour plot of its
two components in Figure \ref{figure4}.
$\psi^{(+)}_D(\x)$ similar to the complex
directional wavelet of \cite{Selesnick2005}, but is not equivalent to
$\psi^{(++)}_D(\x),$ defined in section \ref{dirhyper}. When the Riesz component is constructed from $\psi_D(\x),$ the spatial frequency
modulation to $\Psi_D(\f)$ is $f_1/f$
rather than ${\mathrm{sgn}}(f_1).$ Of course $\Psi_D\left(\f\right)$ is mainly
limited in frequency to $\f=\left[\pm \sqrt{2}f_0,\quad 0 \right]^T,$ and
for $\f=\left[\pm \sqrt{2}f_0+\delta f_1,\quad \delta f_2 \right]^T,$ $f_1/f=\pm
1+\delta f_1/f_0^2+o(\delta f_i^2),$ and so up to a small corrective error
term $\psi_{D,2}(\x)=\psi_{D,3}(\x).$ This is verified by the visual similarity
of Figures (\ref{figure4}) (b) and (c).
The advantage of using a second component
constructed from the Riesz transforms, {\em i.e.} using $\psi_{D,3}(\x)$ rather than using $\psi_{D,2}(\x),$ the
hypercomplex wavelet second component, is that the quaternionic mother wavelet
when rotated satisfies the relationships given by equation (\ref{notcommuating}).
The complex function is then
a bona-fide monogenic function coupled with a single quaternionic conjugate,
an anti-monogenic
function.

\subsection{Quaternionic Morse wavelets}
The 1-D Morse wavelets were introduced
by Daubechies and Paul \cite{DaubechiesPaul1988,Daub} and applied in signal
processing by Bayram and Baraniuk,
as well as Olhede \& Walden \cite{Bayram,Olhede}.
%
\begin{figure*}[t]
\centerline{
\includegraphics[height=1.55in,width=1.55in]{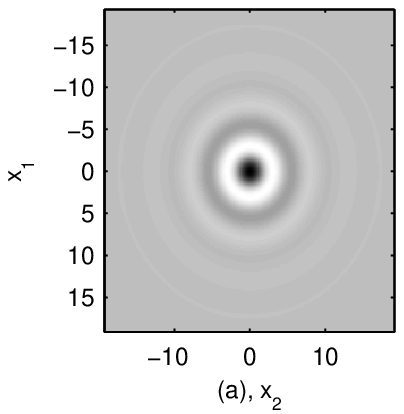}
\includegraphics[height=1.55in,width=1.55in]{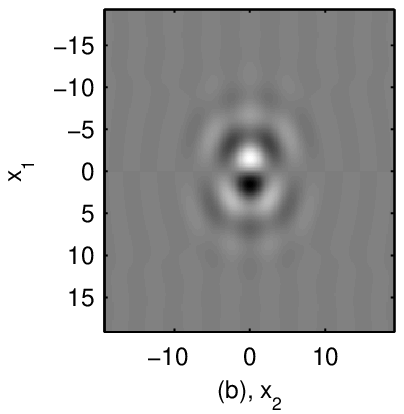}
\includegraphics[height=1.55in,width=1.55in]{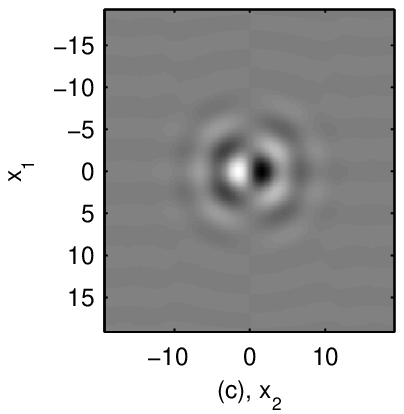}
\includegraphics[height=1.55in,width=1.55in]{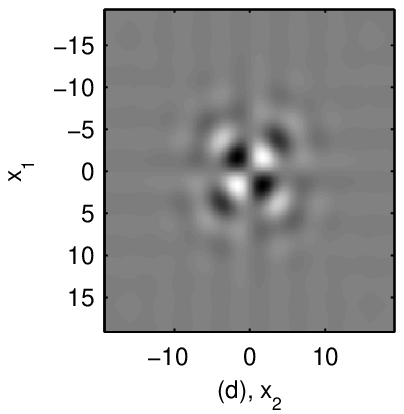}
}
\caption{\label{figure2} 
The real part (a), first, second and third imaginary
parts, (b), (c) and (d), of the isotropic hypercomplexing Morse wavelets, with
$n=0,$ $l=9$ and $m=4.$
}
\end{figure*}
%
In 1-D, according to \cite{Olhede}, the analytic $[\beta,\;\gamma,\;n]$ Morse wavelet is $\psi_{1D,n;\beta,\gamma}^{(+)}(\cdot)=\psi_{1D,n;\beta,\gamma}
^{(e)}(\cdot)
+\bm{j}\psi_{1D,n;\beta,\gamma}^{(o)}(\cdot)
,$ and is defined in
the frequency domain ({\em cf} \cite[p.~2666]{Olhede}) by
\begin{equation}
\Psi_{1D,n;\beta,\gamma}^{(+)}(f)=\sqrt{2}A_{n;\beta,\gamma}(2\pi f)^{\beta}e^{-(2\pi f)^{\gamma}}L_n^c \left[2(2\pi f)^\gamma\right],\;f>0. 
\end{equation}
For $f<0$, $\Psi_{1D,n;\beta,\gamma}^{(+)}(f)$ vanishes. We use the notation 
$c=(2\beta+1)/\gamma-1$, where $\beta,\gamma, \in {\mathbb{R}}^{+}$. Moreover, $
A_{n;\beta,\gamma} = \sqrt{\pi \gamma \ 2^{c+1} \ \Gamma(1+n)/ \Gamma(n+c+1)}$ is a normalisation constant and 
${\rm L}^{c}_{n}(x)$ is the generalised Laguerre polynomial which can be written in a series form:
\[{\rm L}^{c}_{n}(x) = \sum_{r=0}^{n} \ (-1)^{r} \ \\ \frac{\Gamma(1+n+c)}{\Gamma(1+c+r)
\ \Gamma(1+n-r)} \ \frac{x^{r}}{r!}.\]  
Metikas and Olhede \cite{Metikas} defined isotropic real Morse wavelets,
augmented into quaternionic functions, using the Riesz transforms. The work in \cite{Metikas} determined specific families of
isotropic wavelets that could be considered optimally local to a given radial
domain, as well as performed estimation of local image characteristics using what in this article would be
denoted isotropic monogenic wavelets. 
The $n$-th order real 
isotropic Morse wavelet is a radially symmetric function,
$\psi_{n;l,m}^{(e)}(\cdot)$, where $l=\beta+1/2,$ and $m=\gamma$. Its FT
and QFT are given by, \cite{Metikas}, 
\begin{eqnarray}
\Psi_{n;l,m}^{(e)}(f) &\equiv & 
 \Psi_{n;l,m,Q}^{(e)}(f) = A_{n;l,m}^{\prime} \ (2\pi f)^{l} \ e^{-(2\pi f)^{m}}
 \
 {\rm L}^{c^{\prime}}_{n}[ 2 \ (2\pi f)^{m} ], 
 \label{FTQFTmorse}
 \end{eqnarray}
where $c^{\prime}=(2l+2)/m-1,$ $m\ge 1,$ $l>0,$ $l\ge m/2-1,$ and
$A_{n;l,m}^{\prime} = \sqrt{ m \ 2^{c^{\prime}+1} \ \Gamma(1+n)/ \Gamma(n+c^{\prime}+1)}.$
The first two requirements are
necessary to make $\psi_{n;l,m}^{(e)}(f)$ a wavelet, and the second constraint
ensures the wavelet is related to a suitable localization operator \cite{Metikas}.
The 2-D generalization of the analytic Morse wavelets is approached by the
following method:
we first combine the even 1-D Morse wavelets in a suitable
tensor product, or we use isotropic real wavelet functions as $\psi^{(r)}(\x).$ We then, to mimic the analytic wavelets, augment the real mother wavelet
by two or three wavelet functions and form a quaternionic wavelet function,
by ${\cal W}\psi^{(r)}(\x).$


\begin{definition}{Separable hypercomplexing Morse wavelets.}\\
The $\beta,$ and $\gamma,$ separable hypercomplexing Morse wavelets are for any 
$(n_1,n_2)\in{\mathbb{N}}^2$ defined by
\begin{eqnarray}
\label{hypermorse}
\psi_{S,\bm{n};\beta,\gamma}^{++}(\x)&=& \psi_{1D,n_1;\beta,\gamma}^{(e)}(x_1)
\psi_{1D,n_2;\beta,\gamma}^{(e)}(x_2)
+\bm{i}\psi_{1D,n_1;\beta,\gamma}^{(o)}(x_1)\psi_{1D,n_2;\beta,\gamma}^{(e)}(x_2)\\
&&
\nonumber
+\bm{j}\psi_{1D,n_1;\beta,\gamma}^{(e)}(x_1)\psi_{1D,n_2;\beta,\gamma}^{(o)}(x_2)
-\bm{k}\psi_{1D,n_1;\beta,\gamma}^{(o)}(x_1)\psi_{1D,n_2;\beta,\gamma}^{(o)}(x_2).
\end{eqnarray}
\end{definition}
For a plot of the real-valued components of a given
separable quaternionic hypercomplexing Morse wavelet see Figures \ref{figure1}
(a), (b), (c) and (d). The separable
hypercomplexing Morse wavelets exhibit the chequer-board patterns of separable
wavelets, and can, once the CWT has been calculated, be used
to find local separable structure in an image at a given orientation.

\begin{definition}{Isotropic hypercomplexing Morse wavelets.}\\
The $l,$ and $m,$ isotropic hypercomplexing  Morse wavelets are for 
$n\in {\mathbb{N}},$ defined by
\begin{eqnarray}
\label{hypermorse2}
\psi_{I,\bm{n};l,m}^{++}(\x)&=& \psi_{n;l,m}^{(e)}(\x)
+\bm{i}{\cal H}_1\left\{\psi_{n;l,m}^{(e)}\right\}(\x)
+\bm{j}{\cal H}_2\left\{\psi_{n;l,m}^{(e)}\right\}(\x)
-\bm{k}{\cal H}_2{\cal H}_1\left\{\psi_{n;l,m}^{(e)}\right\}(\x).
\end{eqnarray}
\end{definition}
For a plot of the hypercomplexing isotropic Morse wavelets see Figure \ref{figure2}.
It is clear that these wavelets have an isotropic support but recognize variation
associated with the $x_1,$ $x_2$ and diagonal directions.

\begin{definition}{Isotropic monogenic Morse wavelets.}\\
The $l,$ and $m,$ isotropic monogenic Morse wavelets are defined in the frequency domain for $n\in{\mathbb{N}}$
by
\begin{eqnarray}
\label{FTmorsemonogenic}
\Psi^{(+)}_{I,n;l,m}(\f)=\left(1-\bm{k}\cos(\phi)+\sin(\phi)\right)\Psi^{(e)}_{n;l,m}(f).
\end{eqnarray} 
\end{definition}
\begin{figure*}[t]
\centerline{
\includegraphics[height=1.75in,width=1.75in]{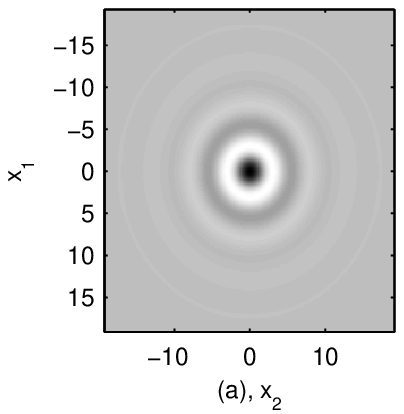}
\includegraphics[height=1.75in,width=1.75in]{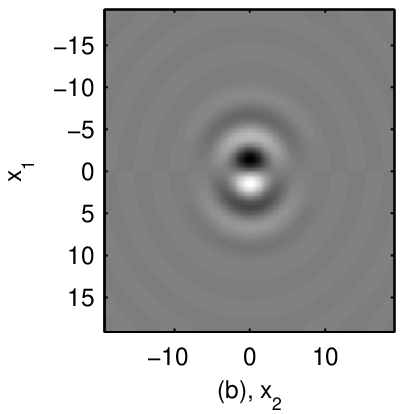}
\includegraphics[height=1.75in,width=1.75in]{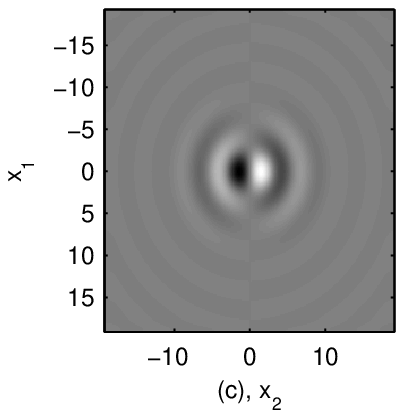}
}
\caption{\label{figure3} The real part (a), and two imaginary parts, (b)
and (c), of the isotropic monogenic Morse  wavelets, with
$n=0,$ $l=9$ and $m=4.$
}
\end{figure*}
By using the inverse Hankel transform the real part of the
wavelets may be retrieved in the spatial domain:
\begin{equation}
\psi_{I,n;l,m}^{(e)}(\x) =  A_{n;l,m}^{\prime} \ 2 \pi \ \left[ \int_{0}^{\infty} df \ f \  (2\pi f)^{l} \ e^{-(2\pi f)^{m}}  \  {\rm L}^{c}_{k}[ 2 \ (2 \pi f)^{m} ] \ J_{0}(2 \pi f x)  \right].
\label{morse}
\end{equation}
The monogenic extension of the $n$-th order mother Morse wavelet is found
by 
\begin{equation}
\psi_{I,n;l,m}^{(+)}(\x) = \psi_{n;l,m}^{(e)}(\x)  + \ ( \bm{i} \ \cos{\chi} + \bm{j} \ \sin{\chi } ) \ 2 \pi \ \left[ \int_{0}^{\infty} df \ f \ \Psi_{n;l,m}^{(e)}(\f)
\ J_{1}(2 \pi f x)  \right].
\label{morsemonogenic}
\end{equation}  
Exact forms for these functions for special choices of parameters can be
determined analytically (see Appendix \ref{spatialmorse}). Figure \ref{figure3}
displays a plot of the isotropic monogenic functions: it is clear that the wavelets have an isotropic support but recognize variation
associated with the two spatial directions. Note that the Poisson scale space
introduced by Felsberg and Sommer \cite{FelsbergSommer2} is not a special
case of a isotropic monogenic Morse wavelet as their analysis function is not zero at $\f={\mathbf{0}},$
and would require setting $l=0.$

\begin{figure*}[t]
\centerline{
\includegraphics[height=1.75in,width=1.75in]{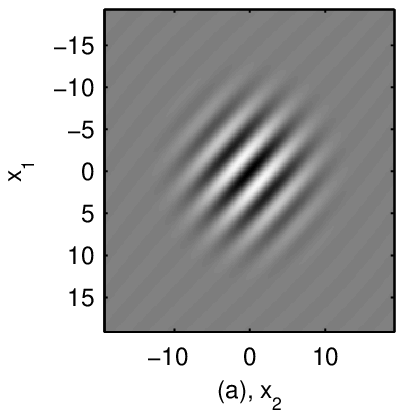}
\includegraphics[height=1.75in,width=1.75in]{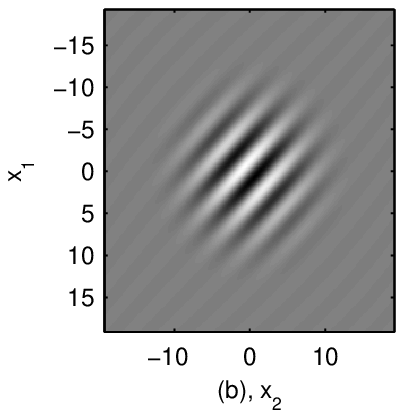}
\includegraphics[height=1.75in,width=1.75in]{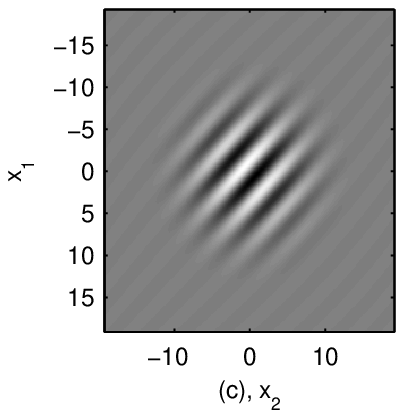}
}
\caption{\label{figure4} The real part $\psi_{D,0;9,4}(\x)$ ((a)), and first imaginary parts ((b) $\psi_{D,3, 0;9,4}(\x)$
and (c) $\psi_{D,2, 0;9,4}(\x)$) of the directional monogenic Morse wavelet and the hypercomplex Morse
wavelet, with
$n=0,$ $\beta=9$ and $\gamma=4.$
}
\end{figure*}

\begin{definition}{Directional monogenic Morse wavelets.}\\
The $\beta,$ and $\gamma,$ directional monogenic Morse wavelets are for 
$n\in{\mathbb{N}}$ defined by
\begin{eqnarray}
\label{dirmonmorse}
\psi_{D,n;\beta,\gamma}^{(+)}(\x)&=& R_{\pi/4}\left\{\psi_{1D,n;\beta,\gamma}^{(e)}
(x_1)\psi_{1D,n;\beta,\gamma}^{(e)}(x_2)
-\psi_{1D,n;\beta,\gamma}^{(o)}(x_1)\psi_{1D,n;\beta,\gamma}^{(o)}(x_2)\right\}
\\
&&
\nonumber
+\bm{i}{\cal{R}}_1\left\{R_{\pi/4}\left\{\psi_{1D,n;\beta,\gamma}^{(e)}(x_1)
\psi_{1D,n;\beta,\gamma}^{(e)}(x_2)
-\psi_{1D,n;\beta,\gamma}^{(o)}(x_1)\psi_{1D,n;\beta,\gamma}^{(o)}(x_2)\right\}\right\}.
\end{eqnarray}
\end{definition}
As expected
$\psi_{D,0;9,3}^{(+)}(\x)$ corresponds to a complex directional wavelet,
and is given by Figures \ref{figure4} (a) and (b), where the other quaternionic
components are identically zero.

\begin{definition}{Hypercomplex Morse wavelets.}\\
The $\beta,$ and $\gamma,$ hypercomplex Morse wavelets are for any integer
valued
$n$ defined by
\begin{eqnarray}
\label{quatbarmorse}
\psi_{D,n;\beta,\gamma}^{(++)}(\x)&=&R_{\pi/4}\left\{\psi_{1D,n;\beta,\gamma}^{(e)}
(x_1)\psi_{1D,n;\beta,\gamma}^{(e)}(x_2)
-\psi_{1D,n;\beta,\gamma}^{(o)}(x_1)\psi_{1D,n;\beta,\gamma}^{(o)}(x_2)\right\}
\\
&&
\nonumber
+\bm{i}\left[R_{\pi/4}\left\{\psi_{1D,n;\beta,\gamma}^{(o)}
(x_1)\psi_{1D,n;\beta,\gamma}^{(e)}(x_2)
+\psi_{1D,n;\beta,\gamma}^{(e)}(x_1)\psi_{1D,n;\beta,\gamma}^{(o)}(x_2)\right\}
\right].
\end{eqnarray}
\end{definition}
$\psi_{D,0;9,3}^{(++)}(\x),$ is thus calculated from
$\psi_{1D,0;9,3}^{(e)}(x_1)$ and $\psi_{1D,0;9,3}^{(o)}(x_1)$ in equation (\ref{dirryyy2}), and plotted in terms of its real and imaginary component
see Figures \ref{figure4} (a) and (c). The real component
is similar to $\psi_{D,0;9,3}^{(+)}(\x),$ but recall that the two imaginary components
are not equal, {\em i.e.} $\psi_{D,2,0;9,3}(\x)\neq \psi_{D,3,0;9,3}(\x).$
Directional wavelets have 
been previously constructed from wavelets related to the 1-D Morse wavelets,
see work by Antoine {\em et al}: \cite{Antoine1,Antoine3,Antoine4},
but these are distinct from $\psi_{D,n;\beta,\gamma}^{(+)}(\x),$ as well
as $\psi_{D,n;\beta,\gamma}^{(++)}(\x).$
The wavelets defined by Antoine {\em et al.} \cite{Antoine1}[p.~319]
are based on the generalized 1-D
Morse wavelets with $\beta=l=m$ and $\gamma=1,$ and are referred to as 2-D Cauchy
wavelets. In their construction a covariant change of axes has been implemented to
make the wavelets compactly supported on a cone in the frequency domain,
rather than by rotating a function supported in the $45^{\circ},$ or $\pi/4$ line. 

\subsection{Choice of mother wavelet function}
The Morse wavelets may be considered as optimally localised over a given
region of space and spatial frequency space, see for example Olhede \& Walden
\cite{Olhede}
and Metikas \& Olhede \cite{Metikas}. 
Note that the hyperanalyticizing 2-D extensions of the
real mother wavelet functions
are localised to the same region as the real wavelet functions. In order to choose an analysis mother wavelet, we either take a real mother wavelet that is isotropic or one which is anisotropic. 
If the directional choice of localisation is necessary to disentangle local
structure, we use a combinations of the tensor products of
1-D wavelets. However, if the isotropic choice of localisation is sufficient
to separate components at given scale and spatial positions, then we base the quaternionic
filters on the isotropic Morse wavelet. 

For a pre-specified spatial and spatial frequency region of either a directional or isotropic form
we choose $\beta$ and $\gamma,$ or equivalently $l$ and $m,$ to calculate a suitable family
of even or isotropic wavelets that are optimally local to a given region
(for a more thorough discussion see \cite{Olhede} or \cite{Metikas}). Thus
$\beta,$ $\gamma,$ or $l$ and $m,$ should be chosen depending on how
the image should be localised in space and spatial frequency.
As a second stage a suitable quaternionic extension of the real-valued mother
wavelet is chosen. The choice of quaternionic extension will depend on what
assumptions can be made regarding the underlying structure of the observed
local variation.
As
may be recalled from Figures \ref{figure1}, \ref{figure2}, \ref{figure3} and
\ref{figure4}, each quaternionic extension represents very different local
oscillatory structures associated with the same period. The separable hypercomplexing wavelets are supported in two
orthogonal directions, each real component of the quaternionic wavelet capable
of extracting even/odd and diagonal structure. The separable filters resemble the quaternionic Gabor
decomposition suggested by von B{\"u}low and Sommer \cite{Buelow1998}, and
the wavelet coefficients represent structure separable in a rotated coordinate
frame.
The directional hyperanalyticizing wavelets are constructed from a tensor product
of 1-D wavelets, but where the directional wavelets are a sum of
1-D
even and odd function tensor products, as was suggested by Selesnick {\em et
al} \cite{Selesnick2005}. Local plane waves, isolated in a given direction
can be represented well by such a decomposition.
The isotropic wavelets mix energy isotropically, {\em i.e.} localise
an image isotropically. Both the hypercomplex and monogenic wavelets
are built from the idea that the real wavelet isolates a single
component by its
scale and spatial localisation -- subsequently the isolated components will
be
represented by their plane-wave (monogenic), or separable wave (hypercomplex)
polar representation. 
This method resembles a Hilbert
spectrum representation of a function in 1-D \cite{Olhede3}. Once the wavelet
transform coefficients have been calculated at a single orientation, if the
assumption of single isolated structure holds, then no further calculations
are necessary. This represents a substantial reduction in computational time.
Once the CWT coefficients have been calculated, the local structure of the analysed
image is represented in terms of the quaternionic transform coefficients, and their polar representation.

\section{Quaternionic wavelet coefficients}
Having constructed various quaternionic mother wavelets, each suitable
for the analysis of local structure, we outline the properties enjoyed
by the coefficients. Equation (\ref{wavtrans}) defined the CWT
coefficients of image $g(\x)$ analysed with a wavelet family constructed
from wavelet $\psi(\x).$ If $\psi(\x)$ is quaternionic then equation
(\ref{wavtrans}) takes the form:
\begin{eqnarray}
\label{CWT}
w_{\psi}(\bm{\xi};g)
&=&\int_{-\infty}^{\infty}  \int_{-\infty}^{\infty} g(\x ) \ \psi^{*}_{\bm{\xi}}(\x)\;d^{2}
\x\\
\nonumber
&=&\int_{-\infty}^{\infty}  \int_{-\infty}^{\infty}
g(\x ) \left[\psi^{(r)}_{\bm{\xi}}(\x)-\bm{i}
\psi^{(i)}_{\bm{\xi}}(\x)-\bm{j}\psi^{(j)}_{\bm{\xi}}(\x)
-\bm{k}\psi^{(k)}_{\bm{\xi}}(\x)\right]d^2{\x}\\
&=&w_{\psi}^{(r)}(\bm{\xi};g)-\bm{i}w_{\psi}^{(i)}(\bm{\xi};g)
-\bm{j}w_{\psi}^{(j)}(\bm{\xi};g)-\bm{k}w_{\psi}^{(k)}(\bm{\xi};g),
\label{linearity}
\end{eqnarray}
thus defining $w_{\psi}^{(u)}(\bm{\xi};g),\;\;u=r,\;i,\;j,\;k.$ 
The local structure at $\bm{\xi}$ is given a
four-vector valued representation, via some phase/s representations or relative
magnitudes of the $w_{\psi}^{(u)}(\bm{\xi};g).$  `Local
energy' is determined from the local magnitude, $\left|w_{\psi}(\bm{\xi};g)\right|,$ of the image $g(\x).$ 
Define $\bm{\zeta}=
\left[a,\theta,\f_b\right]^T,$ where $\f_b$ is the Fourier variable
after the Fourier transform in $\bm{b}$ has been implemented, and 
$\bm{\kappa}=
\left[a,\theta,\bm{q}_b\right]^T,$ where $\bm{q}_b$ is the quaternion Fourier variable
of $\bm{b}.$ 

\begin{definition}{Scalogram of $g(\x).$}\\
The local energy of $g(\x)$ at $\bm{\xi}$ is defined
by the scalogram and is given by
\begin{equation}
\label{scaly}
S_{\psi}(\bm{\xi};g)=\left|w_{\psi}(\bm{\xi};g)\right|^2=
w_{\psi}^{(r)2}(\bm{\xi};g)+w_{\psi}^{(i)2}(\bm{\xi};g)
+w_{\psi}^{(j)2}(\bm{\xi};g)+w_{\psi}^{(k)2}(\bm{\xi};g).
\end{equation}
\end{definition}
Irrespective of what quaternionic mother wavelet function is used, the magnitude square
of the coefficients is interpreted as local signal presence. The interpretation
is appropriate as the four wavelet functions $\psi^{(u)}(\x)$ with $u=r,\;i,\;j,\;k,$
are chosen to be local to $\x=\be$ and scale $a,$ where the localisation
of $\psi^{(u)}(\x)$ with $u=i,\;j,\;k,$ is similar to that of $\psi^{(r)}(\x)$
from our previous comments in sections \ref{propmon} and \ref{prophyp}.
Note that $-\psi^{(oo)}(\x) $ enjoys the same localization as $\psi^{(oo)}(\x).$

The relative
relationships between the other components of the quaternionic wavelet will
depend on the choice of mother wavelet function used to decompose the image,
and so the phase representation of structure will vary with the mother wavelets
used.
To obtain coefficients with interpretable phase relations we calculate 
$w_{\psi}\left(\bm{\xi};g\right),$ so that it corresponds to a $\theta$-hyperanalytic
signal in index $\be$ for any fixed value of $\theta$ and $a.$ 
The signal $g(\cdot)$ could also be represented in
terms of some weights attached with functions $\psi_{\bm{\xi}}(\cdot),$
that are hyperanalytic in $\x$ for any fixed $\bm{\xi}_{0}.$ In the latter case equation
(\ref{reconny5})
states that the signal may be reconstructed in terms of the weighted and rotated hyperanalytic
functions $\psi_{\bm{\xi}}(\cdot).$ In 1-D these two perspectives
coincide, as analysis of a signal in terms of a analytic wavelet yields both
analytic wavelet coefficients and a reconstruction in terms of the analytic
wavelet functions. In 2-D the problem is more complicated, due
to the non-zero $\bm{k}$ component in the hypercomplexing wavelet.
We introduce the additional notation of $\bm{J}_{\pi/2}\bm{\kappa}\equiv \left[a,\theta,\bm{J}_{\pi/2}\q_{b}\right]^T.$ 
\begin{proposition}{Forms of the FT and QFT of wavelet coefficients with real wavelet.}\\
The FT and QFT of $w_{\psi}^{(r)}(\bm{\xi};g)$ for a real signal $g(\x)$
with a real-valued wavelet $\psi^{(r)}(\cdot)$ respectively are
\begin{eqnarray}
W_{\psi}^{(r)}(\bm{\zeta};g)&=& 
G(\f_b ) \ a \Psi^{(r)*}(a \bm{r}_{-\theta} \f_b),
\label{FTCWT}\\
W_{\psi;Q}^{(r)}(\bm{\kappa};g)&=&
\frac{1-\bm{k}}{2} \  G(\q_b ) \ a \Psi^{(r)*}(a \bm{r}_{-\theta} \q_b) + \frac{1+\bm{k}}{2} \ G(\bm{J}_{\pi/2} \q_b ) \ 
a \Psi^{(r)*}(a \bm{r}_{-\theta}\bm{J}_{\pi/2} \q_b)\nonumber \\
&=& \frac{1-\bm{k}}{2}W_{\psi}^{(r)}(\bm{\kappa};g)+
\frac{1+\bm{k}}{2}W_{\psi}^{(r)}(\bm{J}_{\pi/2}\bm{\kappa};g).
\label{QFTCWT}
\end{eqnarray}
\end{proposition}
\begin{proof}
Equations (\ref{FTCWT}) and (\ref{QFTCWT}) follow by direct calculation, and are only stated so that a comparison
can be made with the quaternionic wavelet coefficients.
\end{proof}
\begin{proposition}{Forms of the FT and QFT of wavelet coefficients with quaternionic wavelet.}\\
The FT and QFT of $w_{\psi}(\bm{\xi};g)$ for a real signal $g(\x)$
with a quaternion-valued wavelet $\psi(\cdot)$ respectively are, 
\begin{eqnarray}
\nonumber
W_{\psi}(\bm{\zeta};g)&=& G(\f_b)\left[ 
a \Psi^{(r)*}(a \bm{r}_{-\theta} \f_b) - \bm{i} a \Psi^{(i)*}(a \bm{r}_{-\theta} \f_b) - \bm{j}  a \Psi^{(j)*}(a \bm{r}_{-\theta} \f_b) - \bm{k}  a \Psi^{(k)*}(a \bm{r}_{-\theta} \f_b) \right] 
\label{FTCWT22}\\
\nonumber
W_{\psi;Q}(\bm{\kappa};g)&=& \frac{1}{2}W_{\psi}(\bm{\kappa};g)-
\frac{\bm{k}}{2}W_{\psi^{(r)}-i\psi^{(i)}+j\psi^{(j)}-k\psi^{(k)}}(\bm{\kappa};g)
+\frac{1}{2}W_{\psi}(\bm{J}_{\pi/2}\bm{\kappa};g)\\
&&+\frac{\bm{k}}{2}W_{\psi^{(r)}-i\psi^{(i)}+j\psi^{(j)}-k\psi^{(k)}}(\bm{J}_{\pi/2}\bm{\kappa};g)
\label{QFTCWT22}\\
\label{hahah}
&\neq & 
\frac{1-\bm{k}}{2}W_{\psi}(\bm{\kappa};g)+
\frac{1+\bm{k}}{2}W_{\psi}(\bm{J}_{\pi/2}\bm{\kappa};g).
\end{eqnarray}
\end{proposition}
\begin{proof}
See appendix \ref{nykonstig}.
\end{proof}
Equations (\ref{QFTCWT22}) and (\ref{hahah}) illustrate the new properties afforded by the
quaternionic decomposition, in comparison to a real-valued decomposition of $g(\bm{x})$. The $\bm{i}$ and $\bm{k}$ coefficients will enable
us to characterise structure in the $x_1$ axis, and we are not constrained
to average the coefficients over $x_1$ and $-x_1,$ {\em cf} equation (\ref{QFTCWT}). With a well-chosen
quaternionic 
mother wavelet this will give
us tools to represent oscillatory structures.

\subsection{Stability to changes of axes \label{changax}}
The axes of observation are not necessarily aligned with the local axes of
variation in the image.
The CWT is covariant with respect to the
transformation $\x \mapsto {\bm{r}}_{-\theta^{\prime}}\x-\x^{\prime},$ {\em
i.e.}
if we observe
$g(\x)=g_2({\bm{r}}_{-\theta^{\prime}}(\x-\x^{\prime})),$ rather than $g_2(\x),$
then using equation (2.18) of \cite{Antoine3}[p.~262] on each of the real
valued weights associated with $1,$ $\bm{i},$ $\bm{j},$ and $\bm{k},$ we
may formally note that:
\begin{equation}
\label{equivarrrys}
w_{\psi}(\bm{\xi};g(\x)=
w_{\psi}(a,\theta-\theta^{\prime},\bm{r}_{-\theta^{\prime}}\be-\x^{\prime};g_2(\x)).\end{equation}
Having noted this equivariance of the quaternionic CWT, one
might assume that no discussion needs to be provided of changes of axes of
observation.
However, a discretization of the calculation of the CWT, implies
that we calculate the CWT coefficients at a sampled subset of all values
of the locality index
$\bm{\xi}.$ Small misalignments in space
between the wavelet function and the object under observation 
may
cause the CWT coefficient at a given value of $\bm{\xi}$ to
be small even if the signal has a large contribution at $\x=\be$ and
$\f=a^{-1}\bm{r}_{-\theta} \f_0.$ 
The down-sampling inherent
in most filter bank implementations of the CWT exacerbates
the spatial initialization problem. In 1-D time shift variance has been considered in great detail
see for example \cite{Coifman,Fernandes2,gopinath2003}, and
small misalignment in time may be considered in terms of
shifts in phase.
In particular, \cite{gopinath2003}[p.~1794] discusses the relationship between
small spatial shifts and phase shifts.
It is important that, at any given $\bm{\xi}$, 
the transform of structure that corresponds to the same space
and spatial frequency locality should not correspond to very different $|w_{\psi}(\bm{\xi};g(\x))|$ due to
small misalignments in space between the image and the wavelet. This becomes
equivalent to requiring that the magnitude of the CWT does not change under
phase shifts of the signal.

In 1-D it was shown in
equation (6.11) and the subsequent discussion in \cite{Olhede1} that,
if a real signal of sinusoidal form is phase shifted by a constant angle $\theta_{s}$, the analytic CWT of this signal will correspond
to a complex phase times the CWT of the original signal. This ensures CWT
magnitude
invariance to phase shifts, and produces a stable CWT.
We shall see that the hyperanalytic wavelet coefficients
exhibit similar structure, but the definition of a phase-shifted signal in 2-D will
vary
with the choice of local structure. Depending on the local structure of the
image a phase
shift will be defined as a single shift in a given direction, or as two shifts
in perpendicular directions.

\subsection{The hypercomplex CWT}
The hypercomplex signal of a real image $g(\x)$ represents a real-valued
image
in terms of a four-vector representing local structure, where the local variation
in the image is considered
in terms of the $x_1$ and $x_2$ axes separately, as well as joint variation
in the two axes.
When the $\theta$-hypercomplex
signal is formed, variations in a rotated frame of reference
are considered, and this naturally introduces $\widetilde{g}_{-\theta}(\x),$
 rather than $g(\x)$ as the object under analysis.
If there are separable oscillations present in $g(\x)$ at a spatial position
$\x,$ in direction along and perpendicular to the axes rotated by $\theta,$ then using the polar representation of the $\theta$-hypercomplex
signal given in equation (\ref{polarhyper})
will allow for the characterisation of the local structure. The local structure
is, in this instance, given in terms of the amplitude
or energy associated with the variation,
the period in the first axis' oscillation, given by 
$\widetilde{\alpha}_{-\theta,\theta}(\x),$
as well as the period in the second axis' oscillation, given by 
$\widetilde{\beta}_{-\theta,\theta}(\x).$
If there truly is a completely separable structure present,  
then $\widetilde{\gamma}_{-\theta,\theta}(\x)$
will confirm this assumption by taking the value of zero. If not, a non-zero
value of $\widetilde{\gamma}_{-\theta,\theta}(\x)$ will characterise the non-separability of
the oscillation. By recalculating the decomposition at different values of
$\theta \in [0,\pi/2]$, the orientation may be determined from the analysed
image, as the value of $\theta$ that yields $\widetilde{\gamma}_{-\theta,\theta}(\x)=0.$
Unfortunately most observed images do not correspond to the simple structure
of separable variation at any given $\x,$ but variation exists at many different spatial scales and orientations
at the same spatial
position. For such images some localisation needs to be
implemented when forming the hypercomplex representation of the image variation. For an image corresponding to variation at many different scales, the full structure of the image
at a given spatial position will not be well represented in terms of an amplitude,
and two periods. Once different structures associated with
given variations have been
isolated in the image, the structure of each component may be given a hypercomplex
representation. Under the assumption of local separability the
hypercomplex CWT, or HCWT, will be calculated. 

To isolate the individual component we chose to combine the hypercomplex representation of an image with the CWT. 
The wavelet will localise the image in spatial position, scale and
potentially orientation, and is assumed to be sufficient
to separate the different components. In theory, to construct the HCWT,
the CWT may be combined with the hypercomplex representation
in many different orders. We could start by calculating the hypercomplex
signal annihilating all anti-hypercomplex components in equation (\ref{hypersplit}),
and then localise the signal in position, scale and orientation with
a real valued wavelet, with the aim that
once localised the CWT would isolate a single component that
would naturally be given in polar form. We could equivalently argue that
we should first project the image in position, scale and orientation using
a real wavelet to isolate
a single component and then find the hypercomplex representation of the localised in $\be$ 
signal. It is most convenient to use a single wavelet decomposition
that produces quaternionic coefficients in one calculation. Fortunately, a single
quaternionic wavelet which is not itself hypercomplex can be chosen so that both operations can be
done in a single step. The three procedures give the same
filtered image, and this defines a local
$\theta$-hypercomplex signal with interpretable
magnitude and phase structure.
To demonstrate the properties of the coefficients it is most easy to work in the
Fourier and quaternionic Fourier domains.

\begin{proposition}{FT of HCWT.}\\
The CWT of the real signal $g(\x)$ using the hypercomplexing
wavelet
$\psi^{++}(\cdot)$ has a Fourier transform given by
\begin{equation}
W_{\psi}^{(++)}({\bm{\zeta}};g)=\left[1+{\mathrm{sgn}}
\left(\left[\bm{r}_{-\theta}\f_b\right]_2\right)\right]
\left[1-\bm{k}{\mathrm{sgn}}\left(\left[\bm{r}_{-\theta}\f_b\right]_1\right)
\right]W_{\psi}^{(ee)}(\bm{\zeta};g)
\label{FThypercomplex}
\end{equation}
\end{proposition}
\begin{proof}
The result follows by direct calculation.
\end{proof}
The hypercomplexing wavelet is filtering the signal in the rotated coordinate
system, constructing a signal locally odd in $(\bm{r}_{-\theta}\x)_1$ as the $\bm{i}$ component,
a signal locally odd in the $(\bm{r}_{-\theta}\x)_2$ as the $\bm{j}$ component and a signal locally
odd in $(\bm{r}_{-\theta}\x)_1$ and $(\bm{r}_{-\theta}\x)_2$ as the $\bm{k}$
component. 

\begin{proposition}{QFT of HCWT \label{hfgffg}.}\\
The CWT of the real signal $g(\x)$ using the hypercomplexing wavelet
 $\psi^{++}(\cdot)$ has QFT given by
\begin{eqnarray}
\nonumber
W_{Q, \psi}^{(++)}({\bm{\kappa}};g) &=& \frac{1-\bm{k}}{2}
\left(1+{\mathrm{sgn}}\left(\left[\bm{r}_{-\theta}\q_b\right]_1\right)\right)
\left(1+{\mathrm{sgn}}\left(\left[\bm{r}_{-\theta}\q_b\right]_2\right)\right)
W^{(ee)}_{\psi}(\bm{\kappa};g)
 \nonumber \\
&&+ \frac{1+\bm{k}}{2}
\left(1+{\mathrm{sgn}}\left(\left[\bm{r}_{\theta}\q_b\right]_1\right)\right)
\left(1+{\mathrm{sgn}}\left(\left[\bm{r}_{\theta}\q_b\right]_2\right)\right)
W^{(ee)}_{\psi}(\bm{J}_{\pi/2}\bm{\kappa};g).
\label{QFThypercomplex}
\end{eqnarray}
\end{proposition}
\begin{proof}
The proof of this calculation follows by direct calculation.
\end{proof}
Comparing this equation with equation (\ref{thetahypercquat}) we observe
that this is a $\theta$-hypercomplex function, {\em i.e.} hypercomplex in a rotated
frame of reference in terms of
\begin{eqnarray}
H(\bm{r}_{-\theta}\q)&=&a \Psi^{(ee)*}(a\bm{r}_{-\theta}
\q)G(\q),\;
H(\q_1)=a \Psi^{(ee)}(a\q_1)G(\bm{r}_{\theta}\q_1)\equiv
a \Psi^{(ee)}(a\q_1)\widetilde{G}_{-\theta}(\q_1).
\end{eqnarray}
The hypercomplex CWT simultaneously constructs a $\theta$-hypercomplex signal as localising the signal in frequency. $h(\be)$ corresponds to the CWT of 
$\widetilde{g}_{-\theta}(\x),$ calculated at rotation $\theta=0,$ and we may also interpret $w_{\psi}^{(++)}(\bm{\xi};g)$ as the $\theta$-hypercomplex extension of the real 
CWT with $\theta=0$ of the signal $\widetilde{g}_{-\theta}(\x),$
{\em i.e.} calculating the real CWT of the real signal at $\bm{\kappa}_0=
\left[a,0,\be\right]^T,$ and then extending
it to a $\theta$-hypercomplex signal.
By equation (\ref{thetahypercquat}) we note:
\begin{eqnarray*}
\left(W^{(ee)}_{Q,\psi}(\bm{\kappa}_0;\widetilde{g}_{-\theta})\right)^{(++)}_{\theta}&=&
\left(1+{\mathrm{sgn}}\left(\left[\bm{r}_{-\theta}\q_b\right]_1\right)\right)
\left(1+{\mathrm{sgn}}\left(\left[\bm{r}_{-\theta}\q_b\right]_2\right)\right)
\left(\frac{1-\bm{k}}{2}\right)W^{(ee)}_{\psi}(\bm{\kappa};g)\\
&&+\left(1+{\mathrm{sgn}}\left(\left[\bm{r}_{\theta}\q_b\right]_1\right)\right)
\left(1+{\mathrm{sgn}}\left(\left[\bm{r}_{\theta}\q_b\right]_2\right)\right)
\left(\frac{1+\bm{k}}{2}\right)W^{(ee)}_{\psi}(\bm{J}_{\pi/2}\bm{\kappa};g).
\end{eqnarray*}
Thus  $w^{(++)}_{\psi}(\bm{\xi};g)$ is a $\theta$-hypercomplex
signal and corresponds to constructing the $\theta$-hypercomplex
extension of $w^{(ee)}_{\psi}(\bm{\xi}_0;\widetilde{g}_{-\theta}(\x)).$
This gives
the interpretation of the coefficients as the appropriate hyperanalytic extension
of a signal local to $\bm{\xi}_0 = [a, 0, \bm{\beta}]$ existing, and naturally described, in direction $\theta.$ 

\begin{proposition}{The FT of the HCWT
of the $\theta$-hypercomplex signal.}\\
The Fourier transform of the hypercomplex CWT of the $\theta$
hypercomplex signal takes the form:
\begin{equation}
W_{\psi}^{(++)}(\bm{\zeta};\widetilde{g}_{-\theta,\theta}^{(++)})=
2\left[1+{\mathrm{sgn}}\left(\left[\bm{r}_{-\theta}\f_b\right]_2\right)
\right]\left[1
-\bm{k}
{\mathrm{sgn}}\left(\left[\bm{r}_{-\theta}\f_b\right]_1\right)
\right]
W_{\psi}^{(ee)}(\bm{\zeta},g).
\end{equation}
\end{proposition}
\begin{proof}
See appendix \ref{konstig}.
\end{proof}
This results is intermediary to proving that the order of operations of constructing
the hypercomplex signal and calculating the CWT is not
important. The appendix gives the form of $W_{\psi}^{(++)}(\bm{\zeta};
\widetilde{g}_{-\theta,\theta}^{(\mu_1 \mu_2)}),$ for $\mu_1$ and $\mu_2$ taking
any of the values $\pm.$ For the 1-D signal, only the analytic wavelet
transform of the analytic signal was non-null - the analytic wavelet annihilated
the anti-analytic component. This is not the case for a hypercomplex signal,
as in not all three of the anti-hypercomplex components are annihilated.

\begin{proposition}{The FT of the real even CWT
of the $\theta$-hypercomplex signal.}\\
The Fourier transform of the hypercomplex CWT of the $\theta$-hypercomplex signal takes the form:
\begin{equation}
W^{(e e)}_{\psi}(\bm{\zeta};\widetilde{g}_{-\theta,\theta}^{(++)}) = 
\left(1+{\mathrm{sgn}}
\left(\left[\bm{r}_{-\theta}\f_b\right]_2\right)
\right)\left(1-\bm{k}{\mathrm{sgn}}
\left(\left[\bm{r}_{-\theta}\f_b\right]_1\right)
\right)W^{(e e)}_{\psi}(\bm{\zeta};g).
\end{equation}
\end{proposition}
\begin{proof}
See section \ref{konstig2}.
\end{proof}
We interpret $w^{(e e)}_{\psi}(\bm{\xi};\widetilde{g}_{-\theta,\theta}^{(++)})$ as the
localised $\theta$-hypercomplex signal of the signal observed in a different
axes: this should have a polar representation
that is representative of the signal structure.

\begin{theorem}{Construction of the local $\theta$-hypercomplex signal.}\\
By calculating the hypercomplex CWT of a real signal $g(\x),$
a local $\theta$-hypercomplex signal is constructed.
\end{theorem}
\begin{proof}
See appendix \ref{localhyp} and note that 
\begin{eqnarray}
W_{\psi}^{(++)}(\bm{\zeta};g) &=& \left[W_{\psi}^{(ee)}(\bm{\zeta}_0;\widetilde{g}_{-\theta})
\right]_{\theta}^{(++)} =
\frac{1}{2}W^{(+ +)}_{\psi}(\bm{\zeta};\widetilde{g}^{(++)}_{-\theta,\theta})=
W^{(e e)}_{\psi}(\bm{\zeta};\widetilde{g}^{(++)}_{-\theta,\theta}),\nonumber \\
W_{\psi,Q}^{(++)}(\bm{\kappa};g) &=& \left[W_{\psi,Q}^{(ee)}(\bm{\kappa}_0;\widetilde{g}_{-\theta})
\right]_{\theta}^{(++)}=
\frac{1}{2}W^{(+ +)}_{\psi,Q}(\bm{\kappa};\widetilde{g}^{(++)}_{-\theta,\theta})=
W^{(e e)}_{\psi,Q}(\bm{\kappa};\widetilde{g}^{(++)}_{-\theta,\theta}), \nonumber \\
w_{\psi}^{(++)}(\bm{\xi};g) &=& \left[w_{\psi}^{(ee)}(\bm{\xi}_0;\widetilde{g}_{-\theta})
\right]_{\theta}^{(++)}=\frac{1}{2}
w^{(+ +)}_{\psi}(\bm{\xi};\widetilde{g}^{(++)}_{-\theta,\theta})=
w^{(e e)}_{\psi}(\bm{\xi};\widetilde{g}^{(++)}_{-\theta,\theta}).
\label{localhyper}
\end{eqnarray}
Once one the above three sets of equations in (\ref{localhyper}) is
proved, the 
rest follow through Fourier, Inverse Fourier, Quaternionic Fourier, and Inverse
Quaternionic Fourier transforms. Furthermore the equations shed light on
the interpretation of $w_{\psi}^{(++)}(\bm{\xi};g).$
\end{proof}

$\left[w_{\psi}^{(ee)}(\bm{\xi}_0;\widetilde{g}_{-\theta} )\right]_{\theta}^{(++)}$ is a $\theta$-hypercomplex signal
by proposition
\ref{hfgffg}, thus its polar representation is meaningful, as it corresponds
to a real signal that has been localised in scale, and extended to a $\theta$-hypercomplex
signal.
Furthermore $w_{\psi}^{(ee)}(\bm{\xi};\widetilde{g}^{(++)}_{-\theta,\theta})$ can be
viewed as the local contribution of the $\theta$-hypercomplex signal, and
is therefore also meaningful if there is signal presence at that point.
These two quantities are both equal to half 
$w^{(+ +)}_{\psi}(\bm{\xi};\widetilde{g}^{(++)}_{-\theta,\theta}),$
the hypercomplex CWT of the $\theta$-hypercomplex signal.
Thus $w_{\psi}^{(++)}(\bm{\xi};g)$
has an interpretation in terms of $\widetilde{g}_{-\theta,\theta}(\x),$ the
signal naturally represented in the rotated coordinate system, but observed
in another coordinate system.
$w_{\psi}^{(++)}(\bm{\xi};g)$
can thus be used 
to locally at $\bm{\xi}$ represent a real-valued signal in
terms of phase and amplitude. Unlike von B{\"u}low and Sommer, we need {\em not} 
assume that a single component is present at $\x = \bm{b},$ and it 
suffices to assume that once $g(\x)$ has been localised to $(a,\theta)$
and $\x=\be$, then a single component is present. 
 Then $w^{(+ +)}_{\psi}(\bm{\xi};g )$ can be represented in polar form
using equation (\ref{polarhyper}). Furthermore, if we assume that 
the signal corresponds to a locally separable oscillation, the following theorem specifies
its representation.

\begin{theorem}{The HCWT of a separable oscillatory signal.}\\
The hypercomplex CWT of an oscillatory signal modelled
with $\x^{\prime}=\bm{r}_{-\theta}\x$ as $g(\x)=a_g(\x)$ \\$\cos\left(2\pi
\varphi_{g,1}(x^{\prime}_1)\right)
\cos\left(2\pi
\varphi_{g,2}(x^{\prime}_2)\right)$ where $a_g(\x),$ $\varphi_{g,1}(\cdot)$ and $\varphi_{g,2}(\cdot)$ are assumed to be slowly varying,  
defining $\varphi_{g,u}^{\prime}(x_1)=\frac{\partial}
{\partial x_1}\varphi_{g,u}(x_1),$ for $u=1,\;2,$ and additionally assuming $\varphi_{g,1}^{\prime}(x_1)
-\varphi_{g,2}^{\prime}(x_2)=C,$ a constant for all $x_1$ and $x_2$ such
that $a_g(\x)$ is non-negligible, is given
with $\be^{\prime}=\bm{r}_{-\theta}\be$ by
\begin{eqnarray}
\nonumber
w_{\psi}^{(++)}(\bm{\xi};g)&=& a_g(\be) 
a\Psi^{(ee)}\left(a f_{g,1}(\be,\theta),
af_{g,2}(\be,\theta)\right)
\left[\cos\left(2\pi \varphi_{g,1}(b_1^{\prime})) \right)
\cos\left(2\pi \varphi_{g,2}(b_2^{\prime}) \right)\right.\\
\nonumber
&&+\bm{i}\sin\left(2\pi \varphi_{g,1}(b_1^{\prime}) \right)
\cos\left(2\pi \varphi_{g,2}(b_2^{\prime}) \right)
+\bm{j}\cos\left(2\pi \varphi_{g,1}(b_1^{\prime}) \right)
\sin\left(2\pi \varphi_{g,2}(b_2^{\prime}) \right)
\\
\nonumber
&&\left.+\bm{k}\sin\left(2\pi \varphi_{g,1}(b_1^{\prime}) \right)
\sin\left(2\pi \varphi_{g,2}(b_2^{\prime}) \right)\right]+o(1)\end{eqnarray}
\begin{eqnarray}
w_{\psi}^{(++)}(\bm{\xi};g)&=& a_g(\be) a
\Psi^{(ee)}\left(af_{g,1}(\be,\theta),
af_{g,2}(\be,\theta)\right)
e^{2\pi\bm{i}\varphi_{g,1}(b_1^{\prime}) }e^{2\pi\bm{j} \varphi_{g,2}(b_2^{\prime})}
+o(1).
\label{wavridge1}
\end{eqnarray}
where $\bm{f}_g(\x,\theta)=\left[\varphi_{g,1}^{\prime}(x^{\prime}_1),\;\;
\varphi_{g,2}^{\prime}(x^{\prime}_2)\right]^T.$ Note that the constraints imposed on $g(\x)$ includes
a local plane wave
structure as then one of the two phase functions will itself be constant,
and separable structure with different local periods, as this will be modelled
via two different local plane waves separated to different values of $a$ and $\theta$.
\end{theorem} 
\begin{proof} The proof follows by direct calculation, if some care is taken
with the conjugation and appropriate assumptions are made regarding the variability
of the amplitude and phase functions. It is necessary to use a hypercomplexing
wavelet for the formula to follow.
\end{proof}
Let $\f_{0}=\arg_{f_1>0} \max \left|\Psi^{(ee)}(\f)\right|.$ We note that
a simplified description of the separable oscillatory signal $g(\x)$ can be determined
from the ridge of the CWT {\em cf} \cite{GonnetTorresani}[p.~394] via analysis for values of $\bm{\xi}$ in
a subspace given by
$\left\{\bm{\xi}:\;\bm{f}_g(\be,\theta)=a^{-1}\f_{0}\right\}$ such that the measure of the degree of separability of $w^{(+ +)}_{\psi}( \bm{\xi};g)$ is zero.
Direct calculation shows that, when the CWT is calculated at another
angle than that of the variation in $g(\x)$, the magnitude of the wavelet
transform is less than that of equation (\ref{wavridge1}). 
The ridges, and thus $\bm{f}_g(\be,\theta),$ can be determined from
maxima in 
$\left|w_{\psi}^{(++)}(\bm{\xi};g)\right|,$ when the latter magnitude is
calculated for fixed $a$ and $\be$ varying $\theta,$
but to verify
the chosen form of the separable orientation is appropriate we additionally
calculate the polar representation of $
w_{\psi}^{(++)}(\bm{\xi};g)$ and determine that indeed $\widetilde{\gamma}_{-\theta,\theta}(\be)=0,$
where the latter is the $\bm{k}$ angle in the polar representation of $
w_{\psi}^{(++)}(\bm{\xi};g),$ see equation (\ref{polarhyper}). Once $\theta^{\prime},$
has been determined
we calculate $\widetilde{\alpha}_{-\theta^{\prime},\theta^{\prime}}(\be)$ and
$\widetilde{\beta}_{-\theta^{\prime},\theta^{\prime}}(\be),$
these two functions then characterising
the structure of the locally separable oscillation, where
$\left|w_{\psi}^{(++)}(\bm{\xi};g)\right|$ corresponds to the magnitude of the oscillation. 
 Signals
corresponding to aggregations of several separable oscillatory components
can be locally analysed
as long as the wavelets are sufficiently concentrated in space and spatial
frequency. For further notes on separating distinct components see 
\cite{GonnetTorresani}[p.~395], where such methods {\em mutatis mutandis}
can be applied in this context. The main focus in this article is not the
development of additional wavelet ridge methods, but the expression in equation
(\ref{wavridge1}) shows the structural representation possible with the transform
coefficients. 

%
\subsection{Phase shifts in both orthogonal axes}
We define phase shifts in both directions of variation for separable signals.
The shifted signal depends
on the axes of the direction of variation specified via $\theta,$ something
that
will not be the case for the single phase-shift defined via the monogenic structure.

\begin{definition}{Phase shift for a separable structure of orientation $\theta.$}\\
For a real signal that is separable when viewed in the correct axes of observation,
{\em i.e.} $g(\x)$ satisfies equation (\ref{sepintheta}), the signal 
phase shifted by $\bm{\theta}_s=\begin{pmatrix}\theta_{s,1} & \theta_{s,2}\end{pmatrix}^T$ is defined by:
\begin{eqnarray}
\nonumber
\Lambda_{\bm{\theta}_{s},\theta}^{2D} g(\x) &\equiv&  
\left|\widetilde{g}^{(++)}_{-\theta,\theta}(\x)\right| \ \cos{( 2\pi\widetilde{\alpha}_{-\theta,\theta}(\x)
- \theta_{s,1} )} \cos{( 2\pi\widetilde{\beta}_{-\theta,\theta}(\x)
- \theta_{s,2} )}\nonumber \\
&=&\Re\left\{|\widetilde{g}^{(++)}_{-\theta,\theta}(\x)| e^{( \bm{i}
\left(2\pi\widetilde{\alpha}_{-\theta,\theta}(\x)
- \theta_{s,1}\right) )} e^{( \bm{j}\left(2\pi\widetilde{\beta}_{-\theta,\theta}(\x)
- \theta_{s,2} \right))}
\right\},
\label{phaseq1}
\end{eqnarray}
where  $\Lambda_{\bm{\theta}_{s},\theta}^{2D}$ is denoted the phase shift
operator. $\Lambda_{\bm{\theta}_{s},\theta}^{2D}$  shifts the
two separable oscillations in cycle by $\theta_{s,1}$ and $\theta_{s,2},$
respectively. 
\end{definition}

\begin{proposition}{Hypercomplex signal of the phase-shifted separable signal.}\\
For a real signal $g(\x)$ satisfying equation (\ref{sepintheta}) the 
$\theta$-hypercomplex extension of the 
phase shifted real image in direction $\theta$ is given by:
\begin{eqnarray}
\nonumber
\widetilde{(\Lambda_{\bm{\theta}_{s},\theta} g )}_{-\theta, \theta}^{(++)}(\x) &=& 
\nonumber
|\widetilde{g}^{(++)}_{-\theta,\theta}(\x)| e^{( \bm{i}\left[2\pi\widetilde{\alpha}_{-\theta,\theta}(\x)
- \theta_{s,1} )\right]}  e^{ \bm{j}\left[2\pi\widetilde{\beta}_{-\theta,\theta}(\x)
- \theta_{s,2} )\right]}
\\
&=& e^{ - \bm{i} \theta_{s,1} }\  \widetilde{g}^{(+)}_{-\theta,\theta}(\x)e^{
- \bm{j} \theta_{s,2} }.
\label{ShiftedCWThyper}
\end{eqnarray}
\end{proposition}
\begin{proof}
See section \ref{2dphasehyper}.
\end{proof}

\begin{theorem}{HCWT of a phase-shifted separable signal.}\\
The CWT of a phase-shifted signal, where $g(\x)$ satisfies equation (\ref{sepintheta}),
using a
hypercomplex wavelet 
is given by
\begin{equation}
\label{phassyhyp}
w_{\psi}^{(++)}(\bm{\xi};\Lambda_{\bm{\theta}_{s},\theta}^{2D} g)=e^{-\bm{i}\theta_{s,1}}
w_{\psi}^{(++)}(\bm{\xi};g)e^{-\bm{j}\theta_{s,2}}
\end{equation}
\end{theorem}
\begin{proof}
See appendix \ref{phasehyperinv}.
\end{proof}

\begin{corollary}{Magnitude invariance of the HCWT of a separable signal under
phase-shifts.}\\
The magnitude of the CWT of a phase-shifted separable signal is equivalent
to the CWT
of the original signal.
\end{corollary}
\begin{proof}
From equation (\ref{phassyhyp}) we may note that $
\Lambda_{\bm{\theta}_{s},\theta}^{2D} w^{(++)}(\bm{\xi}; g)=e^{\bm{i}\theta_{s,1}}
w_{\psi}^{(++)}(\bm{\xi};g)e^{\bm{j}\theta_{s,2}},$ and thus it follows that
\begin{equation}
\left|w_{\psi}^{(++)}(\bm{\xi};\Lambda_{\bm{\theta}_{s},\theta}^{2D} g)\right|^2
=\left|w_{\psi}^{(++)}(\bm{\xi};g)\right|^2.
\end{equation}
\end{proof}
It thus follows that {\em if} the original image was separable then the magnitude of the CWT coefficients calculated at the correct orientation
is stable with respect to phase-shifts. Given the dependence of the analysis
on the orientation of the signal it is not surprising to find that the analysis
strongly depends on the local orientation of the image. Our representation of the same image {\em should}
change substantially if we assume separability in different directions. We may deduce that for signals that
can be described as locally separable that small
changes in the spatial alignment between the wavelets and the signal will
not cause migration of energy over frequency bands. The result is important,
as if discrete filter banks are constructed to implement the CWT
then $\psi^{++}(\x)$ phase-shift invariant in magnitude implies that
the sampled wavelet coefficients will not suffer from substantial spatial
shift variance,
but there is still orientation dependence.

\subsection{The monogenic CWT}
The interpretability of the hypercomplex decomposition coefficients depends
on the
assumption of locally separable variation present in the image.
To be able to retrieve
the structural representation of the signal the appropriate axes for analysis
have to
be identified; this in general requires the
calculation of the HCWT at the full set of $\theta\in\left[0,\pi\right],$
unless it transpires that the locally preferred orientation is already known.
If at a fixed spatial location $\x=\be$ and given set of frequencies there
was
variation in only one orientation, the structure could
be represented by the scale localised monogenic signal. 
If at any given spatial point $\x,$ $g(\x)$ corresponded to a
single component, then the local characteristics
of $g(\x)$ could be extracted from $g^{(\pm)}(\x),$ and the
polar representation of equation (\ref{polarform}) can be used directly to
represent the signal.
However, in general,
$g(\x)$ will correspond to a multiscale structure at the spatial point $\x,$
and thus a scale local representation using the
monogenic wavelets is necessary to implement in order to produce an interpretable polar representation.
The monogenic CWT, or MCWT, represents local
structure in terms of locally unidirectional variation, where the direction
of variation can be determined from the transform coefficients at a single
local point. In general, analysis with a
monogenic, but not necessarily isotropic monogenic mother wavelet, is implemented
{\em i.e.} a quaternionic mother wavelet of the form
$\psi(\cdot)=\psi^{(+)}(\cdot)=
\psi^{(r)}(\cdot)+\bm{i}\psi^{(1)}_R(\cdot)+\bm{j}\psi^{(2)}_R(\cdot)$ is
used. 
A non-isotropic mother wavelet is a suitable choice as it is not reasonable to assume that the scale localisation alone
will be sufficient to separate components present in the image. We use an isotropic mother wavelet, when, locally, there is only a single unidirectional oscillation
present.

The HCWT were constructed to ensure that the
coefficients rather than the wavelets had
interpretable polar representation. Unfortunately, we found that these two requirements
could not be simultaneously achieved. It will transpire for the monogenic wavelets
that even
if the wavelets are themselves $\theta$-monogenic, so are the wavelet coefficients,
unlike the previous case. This all is due to the zero
$\bm{k}$-component
in the quaternionic representation.
The coefficients are represented in terms of an amplitude,
an orientation and a local period. Again, to be able to establish the properties of the monogenic coefficients, the transform
in the Fourier domain is considered in detail. We show that the MCWT
annihilates the anti-monogenic
component, and this ensures that, equivalently, firstly the
monogenic signal may be constructed and then scale localised; or firstly
the
image can be scale localised and then 
the monogenic extension of the local signal constructed. Indeed both operations
can be implemented
in one step
using the MCWT,
and this then establishes the interpretation of the wavelet coefficients.

\begin{proposition}{FT \& QFT of the MCWT.}\\
The FT, and the QFT, of the MCWT of a real signal are 
given by:
\begin{eqnarray}
W_{\psi}^{(+)}(\bm{\zeta};g)&=&\left[1 -\bm{k}\cos{(\phi_b - \theta)}+\sin{(\phi_b - \theta)}\right]
W_{\psi}^{(r)}(\bm{\zeta};g)
\label{FTCWTaw} \\
W_{\psi;Q}^{(+)}(\bm{\kappa};g)&=&  W_{\psi;Q}^{(r)}(\bm{\kappa};g) + \cos{\theta} \left[ \cos{\nu}_b + \sin{\nu}_b \right] W_{\psi;Q}^{(r)}(\bm{\kappa};g)
 + \sin{\theta} \left[ \cos{\nu}_b - \sin{\nu}_b \right] \bm{i} W_{\psi;Q}^{(r)}(\bm{\kappa};g) \bm{j}. \nonumber
\end{eqnarray}
\end{proposition}
\begin{proof}
See section \ref{monowav}. 
\end{proof}
From these equation the frequency domain properties of the monogenic wavelet
coefficients can be determined. Note that from equations (\ref{FTCWTaw})
and (\ref{weirdeqn}) we
can immediately deduce that the MCWT coefficients are $\theta$-monogenic.
If the real mother wavelet is radially symmetric, $\psi^{(r)}(\x)=\psi^{(e)}(x),$ then the rotation has no important effect in the definition of CWT, (\ref{CWT}).
This assumption implies that 
$\Psi^{(e)}(f) \in {\mathbb{R}},$ and $a \Psi^{(e)*}(a \bm{r}_{-\theta } \f)
= a\Psi^{(e)}(af)$. 

\begin{corollary}{FT and QFT of the isotropic MCWT.}\\
The FT, and the QFT, of the MCWT based on a real isotropic mother wavelet of a real signal is
given by:
\begin{eqnarray}
W_{\psi}^{(e)}(\bm{\zeta};g) & = &  \Psi^{(e)}(af_b)G(\f_b),
\label{RFTCWT}\\
W_{\psi}^{(+)}(\bm{\zeta};g) &=&   \left[1 +  \ \left( - \bm{k} 
\ \cos(\phi_b-\theta) +  \sin(\phi_b-\theta) \right) \right] \ \Psi^{(e)}(af_b)G(\f_b) \nonumber \\ 
&=&\Psi^{(+)}(a\bm{r}_{-\theta}\f_b)G(\f_b)\neq \Psi^{(+)*}(a\bm{r}_{-\theta}\f_b)G(\f_b)
\nonumber\\
&=&  \left[1 +  \ \left( - \bm{k} 
\ \cos(\phi_b-\theta) +  \sin(\phi_b-\theta) \right) \right] \    W_{\psi}^{(e)}(\bm{\zeta};g)
\nonumber
\\
W_{\psi;Q}^{(+)}(\bm{\kappa};g) &=& \left\{ G_{Q}(\q_b) \left[1 + \cos{\theta} \left( \cos{\nu}_b + \sin{\nu}_b \right) \right]
 + \sin{\theta} \left[ \cos{\nu}_b - \sin{\nu}_b \right] \bm{i} G_{Q}(\q) \bm{j} \right\} a \Psi^{(e)}(aq_b).
\nonumber
\end{eqnarray} 
\end{corollary}  
These relations will enable us to consider the properties of the magnitude
of the wavelet coefficients under rotation, and we wish to derive certain
relations regarding the properties of the CWT in these circumstances.

\begin{proposition}{FT of the MCWT.}\\
The FT of the MCWT of a monogenic extension
of a real signal takes the form:
\begin{equation}
W_{\psi}^{(\pm)}(\bm{\zeta};g^{(+)})=
\left[\left(1\pm\cos(\theta) \right)G^{(+)}(\f_b)\pm\bm{k} \sin(\theta) G^{(-)}(\f_b)
\right]a\Psi^{(r)*}\left(a\bm{r}_{-\theta}\f_b \right).
\end{equation}
The anti-monogenic extension has an equivalent form with the $\pm$ on the
right hand side replaced by $\mp.$
\end{proposition}
\begin{proof}
See section \ref{splitting} and set $\theta^{\prime}=0.$
\end{proof}

\begin{proposition}{CWT of $\theta$-monogenic \& anti-monogenic decomposition components.}\\
The Fourier transform of the MCWT of the extended
$\theta$-monogenic and anti-monogenic extension
of a real signal take the forms:
\begin{eqnarray}
\nonumber
W_{\psi}^{(+)}(\bm{\zeta};\widetilde{g}_{-\theta,\theta}^{(+)})&=& \left\{  
2+2\sin(\phi_b-\theta)-2\bm{k}\cos(\phi_b-\theta)  \right\}  \nonumber \ 
W^{(r)}_{\psi}(\bm{\zeta};g) \nonumber \\
&& =2\widetilde{G}^{(+)}_{-\theta,\theta}(\f_b)
a\Psi^{(r)*}(a \bm{r}_{-\theta}\f_b).
\nonumber\\
&=&2W_{\psi}^{(r)}(\bm{\zeta};\widetilde{g}_{-\theta,\theta}^{(+)}),\quad
w_{\psi}^{(+)}(\bm{\xi};\widetilde{g}_{-\theta,\theta}^{(+)})=
2w_{\psi}^{(r)}(\bm{\xi};\widetilde{g}_{-\theta,\theta}^{(+)})\nonumber\\
W_{\psi}^{(+)}(\bm{\zeta};\widetilde{g}_{-\theta\theta}^{(-)})&=&0,\quad
w_{\psi}^{(+)}(\bm{\xi};\widetilde{g}_{-\theta\theta}^{(-)})=0
\label{FTCWTmnkl3}
\end{eqnarray}
\end{proposition}
\begin{proof}
See section \ref{splitting}.
\end{proof}\\
Hence the MCWT of a $\theta$-monogenic signal observed
in the rotated frame of reference, {\em i.e.} $\widetilde{g}_{-\theta,\theta}^{(+)}(\x)$
is twice that of the real  wavelet
transform of the $\theta$-monogenic signal $\widetilde{g}_{-\theta, \theta}(\x).$
The MCWT of $\widetilde{g}_{-\theta,\theta}^{(-)}(\x)$ is zero. Thus the monogenic
wavelet annihilates the anti-monogenic component of the image, in perfect
analogy with 1-D, but we may further generalise the result to consider any
combination of monogenic/anti-monogenic wavelet and monogenic/anti-monogenic
signal.

\begin{corollary}{Annihilation using Monogenic \& Anti-Monogenic Wavelets.}\\
The MCWT annihilated the $\theta$-anti-monogenic image
and the anti-monogenic CWT annihilated the $\theta$-monogenic image:
\begin{eqnarray}
W_{\psi}^{(\pm)}(\bm{\zeta};g)   &=& \frac{1}{2} \  
W_{\psi}^{(\pm)}(\bm{\zeta};\widetilde{g}^{(\pm)}_{-\theta,\theta}),\quad
w_{\psi}^{(\pm)}(\bm{\xi};g)   = \frac{1}{2} \  
w_{\psi}^{(\pm)}(\bm{\xi};\widetilde{g}^{(\pm)}_{-\theta,\theta}).
\label{phaseinvariance}
\end{eqnarray}
For any other angle than $\theta \pm l \pi, l \in {\cal Z}$, $ w_{\psi}^{(\pm)}(\bm{\xi};g)$ is a linear combination of $
w_{\psi}^{(\pm)}(\bm{\xi};g^{(+)}_{-\theta,\theta})$ and $w_{\psi}^{(\pm)}(\bm{\xi};g
g^{(-)}_{-\theta,\theta})$.
\end{corollary}
\begin{proof}
Comparing equation (\ref{splitting})
with equation (\ref{FTCWTaw}) the first result follows directly: furthermore
the annihilation follows {\em mutatis mutandis} from equation (\ref{FTCWTmnkl3p-}).
\end{proof}
Thus, when analysing a real signal, the MCWT can be
used
to annihilate the anti-monogenic components in the real image, which
is in complete analogy to the $d=1$ analytic/anti-analytic case, \cite{Olhede1, Olhede2}. See equations (6.8) to (6.10), and in particular the unnumbered equation over equation (6.10), of \cite{Olhede1}.
These show that the CWT of an $d=1$ real signal with respect to an analytic (anti-analytic) wavelet
is equal to one half of the CWT of the 
analytic (anti-analytic) extension of the real signal with respect to the
analytic  
(anti-analytic) wavelet.  
In other words equation (\ref{phaseinvariance}) is a generalisation of the equation
over equation (6.10) of \cite{Olhede1}. 
This implies that both in $d=1$ and $d=2$, the CWT of the real signal or
image has the same phase and one half of the modulus of the CWT of the 
analytically or monogenically extended signal or image.
\begin{corollary}{Annihilation of $\pi$-rotated component.}\\
If we assume $\psi^{(r)}=\psi^{(e)},$ 
then with the additional notation of $\bm{\zeta}_{\pi}=
\left[a,\theta+\pi,\bm{b}\right]^T,$ it follows that:
\begin{equation}
W_{\psi}^{(\pm)}(\bm{\zeta}_{\pi};g)   = \frac{1}{2} \  
W_{\psi}^{(\pm)}(\bm{\zeta}_{\pi};\widetilde{g}^{(\mp)}_{-\theta,\theta}),\quad
w_{\psi}^{(\pm)}(\bm{\xi}_{\pi};g)   = \frac{1}{2} \  
w_{\psi}^{(\pm)}(\bm{\xi}_{\pi};\widetilde{g}^{(\mp)}_{-\theta,\theta}).
\end{equation}
\end{corollary}
If an isotropic wavelet is used then the monogenic wavelet rotated by $\pi$
is identical to
the anti-monogenic wavelet, see equation (\ref{equivalentpirot}), and so
the result is an immediate consequence of the previous proposition.

\begin{theorem}{Construction of the local $\theta$-monogenic signal.}\\
The MCWT of a real image $g(\x)$ is equivalent to
the scale localisation of the $\theta$-monogenic extension of the signal
and the $\theta$-monogenic extension of the scale-localised version of the
image $g(\x):$ 
\begin{eqnarray}
\label{propi=ii}
w^{(+)}_{\psi}(\bm{\xi};g)&=&\left[w^{(r)}_{\psi}(\bm{\xi}_0;\widetilde{g}_{-\theta}
)\right]^{(+)}_{\theta}=w^{(r)}_{\psi}(\bm{\xi};\widetilde{g}_{-\theta,\theta}^{(+)})
\end{eqnarray}
\end{theorem}
\begin{proof}
See appendix \ref{ftthetamon}.
\end{proof}
Thus all three local descriptions, {\em i.e.} the
MCWT of a real signal,
the real CWT of a $\theta$-monogenic signal and
the monogenic extension of a real CWT of a real signal of
a real-valued signal may all be viewed as equivalent.
The left hand side of equation (\ref{propi=ii}) simply corresponds
to the MCWT of the real signal, whilst the middle
equation is the $\theta$-monogenic extension of the scale local signal 
$\widetilde{g}_{-\theta}(\x),$
which corresponds to adding some suitable components to the real part of
the CWT of $g(\x).$ The right hand side of equation (\ref{propi=ii})
corresponds to first forming the monogenic extension of the signal, and then
scale-localising the extension to make the phase description of the signal
interpretable. Hence the MCWT can, assuming a single component
has been retrieved at local index $\bm{\xi},$ be represented in polar form
using equation (\ref{polarform}) where the magnitude, phase and orientation are interpretable
in terms of a local univariate variations and the $\theta$-monogenic representation.
Furthermore if we may assume
the signal corresponds to a local plane wave the following theorem specifies
its representation.

\begin{theorem}{The MCWT of a plane oscillatory signal.}\\
The MCWT of a single component separable oscillatory signal modelled
by \\$g(\x)=a_g(\x)\cos\left(2\pi
\varphi_{g}(\x)\right)$ where $a_g(\x),$ as well as $\varphi_{g}(\x)$ are assumed to be slowly varying, is given
by:
\begin{eqnarray}
w_{\psi}^{(+)}(\bm{\xi},g)&=&a_g(\be)a\left|\Psi^{(r)}(a \bm{r}_{-\theta}\bm{f}_g(\be))\right|
e^{2\pi \bm{e}_{\bm{n}_g(\be,\theta)}
\left[\varphi_g(\be)-\varphi_{\psi}(a \bm{r}_{-\theta}\bm{f}_g(\be))\right]}+o(1),
\label{wavridge2}
\end{eqnarray}
where $\bm{f}_g(\x)=\nabla \varphi_{g}(\x)=f_g(\x)\left[\cos(\phi_g(\x))\quad
\sin(\phi_g(\x))
\right]^T,$\\ $\bm{n}_g(\x,\theta)=\left[\cos(\phi_g(\x)-\theta)\;\;
\sin(\phi_g(\x)-\theta)
\right]^T,$ $\Psi^{(r)}(\f)=\left|\Psi^{(r)}(\f)\right|
e^{-2\pi\bm{j} \varphi_{\psi}(\f)},$ and note that Hermitian symmetry imposes
$\left|\Psi^{(r)}(\f)\right|=
\left|\Psi^{(r)}(-\f)\right|$ as well as $ \varphi_{\psi}(\f)= -\varphi_{\psi}(-\f).$
\end{theorem}
\begin{proof}
This result was shown for the special case of the isotropic Morse wavelets
in Metikas \& Olhede \cite{Metikas}, and follows by direct calculation.
\end{proof}
Let $\f_{0}=\arg_{f_1>0} \max \left|\Psi^{(r)}(\f)\right|.$ We note that
a simplified description of the oscillatory signal $g(\x)$ can be determined
from the ridge of the CWT {\em cf} \cite{GonnetTorresani}[p.~394] via
the subspace of the locality index given by
$\left\{\bm{\xi}:\;\bm{f}_g(\be)=a^{-1}\bm{r}_{\theta}\f_{0}\right\}$.
Equation (\ref{wavridge2}) can be used to characterise the oscillation at
$\x=\be$ and if the signal would be more appropriately modelled as an aggregation
of of oscillatory signals as long as they are sufficiently separated, {\em
i.e.}
as long as the wavelet is sufficiently narrow in space and spatial frequency
to separate the different component. Further analysis of the signal such
as that of Olhede and Metikas \cite{Metikas,Isspit,Spie},
using isotropic monogenic wavelets,
can then be implemented for arbitrary monogenic wavelets. Using directional
monogenic wavelets will allow us to analyse a larger class of signal as it
is more reasonable to assume an image corresponds to a collection of plane-waves
if we localise in scale, position {\em and} orientation.

\subsection{Phase shifts for Unidirectional Variation}
Section \ref{changax} noted that the stability of the CWT coefficients
under given affine transformations of the argument was important. For the
HCWT we remarked that the coefficients exhibited
desirable properties for phase-shifted separable structure. With the monogenic
transform
we can obtain desirable properties for an arbitrary phase-shifted signal,
and in some instances obtain invariance under rotation. As the local
structure is described as a plane wave with orientation determined from the
signal under observation this enables us to derive the additional results.
Furthermore the phase-shift is not dependent on
the choice of axes: the monogenic wavelet is identifying the direction of
variation from the observed image, and then the phase-shift is given in terms
of this direction. This is not equivalent to the hypercomplex wavelet coefficient
structure.
\begin{proposition}{Local magnitude invariance under rotation.}\\
If the real mother wavelet is isotropic, then the MCWT coefficients
of
a real image $g(\x)$
at
$\bm{\xi}$ will have a magnitude invariant to the value of $\theta.$ 
\end{proposition}
\begin{proof}
See section \ref{Locinv}.
\end{proof}
Hence the magnitude of the wavelet coefficients does not depend on the rotation
angle, and we may deduce the result. 
This implies that we are not sensitive to loss of magnitude due to local
orientation misalignment. However, we would naturally like to also be able to retrieve
the local orientation, from the MCWT.
\begin{definition}{A directional signal.}\\
A real signal $g(\x)$ is considered as one-dimensionally directional
at frequency $\f$ if $\exists\phi_0\in\left[0,\pi\right)$ such that
\begin{equation}
G\left(\f\right)
= \frac{\widetilde{G}(f)}{2}
\left[\delta(\phi-\phi_0-\pi)+\delta(\phi-\phi_0)\right].
\label{dirsignaldef}
\end{equation}
\end{definition}
We shall in this article be analysing real images: as $g(\x)\in {\mathbb{R}}$
its FT is Hermitian
$
G^{*}(\f ) = G(- \f),$ and thus necessitates having both a delta distribution
component at $\phi_0$ and $\phi_0+\pi.$

\begin{theorem}{Directional selectivity of the MCWT.}\\
If the real mother wavelet is isotropic, and the analysed signal $g(\x)$
is directional over the support of the monogenic wavelet at $\bm{\xi}$ with
directionality $\phi_0,$ 
then with $\bm{\xi}^{\dagger}=\left[a,\phi_0,\bm{b}\right]^T$
\begin{eqnarray}
\left|w_{\psi}^{(1)}\left(\bm{\xi}^{\dagger};g \right)
\right|^2
&=&\left|w_{\psi}^{(1)}\left(\bm{\xi};g \right)
\right|^2+\left|w_{\psi}^{(2)}\left(\bm{\xi};g \right)
\right|^2,\quad
\phi_0-\theta=\tan^{-1}\left(\frac{w_{\psi}^{(2)}\left(\bm{\xi};g \right)}
{w_{\psi}^{(1)}\left(\bm{\xi};g \right)} \right)
\end{eqnarray}
\end{theorem}
\begin{proof}
See section \ref{Dirsel}.
\end{proof}
Thus with an isotropic mother wavelet, the magnitude of the Riesz component
wavelet coefficients is invariant to rotations, but is maximised in the first
component if we rotate the wavelet to align with the directionality of the
variations. Hence we may determine the local directionality of the signal
from the monogenic
isotropic wavelet coefficients.

\subsubsection{Phase-shifted plane wave signals \label{Phase Shift}}
For a plane wave present at a given point $\x$ we define the phase-shift
operation as follows.
\begin{definition}{Phase-shift of a plane wave.}\\
For a real signal $g(\x)=\widetilde{g}_{-\theta}(\bm{r}_{-\theta}\x)$ we note from equation
(\ref{decomposition}) that the image
may for any rotation angle $\theta$ be written in sinusoidal form and the
phase shifted by $\theta_s$ version of the signal is
defined via:
\begin{eqnarray}
\Lambda_{\theta_{s}} g(\x) &\equiv&  
\left|g^{(+)}(\x)\right| \ \cos{( 2\pi\phi(\x)
- \theta_{s} )} =\Re\left\{|g^{(+)}(\x)| \exp\left( (2\pi\phi (\x) - \theta_{s})
\bm{e}_{\widetilde{\nu}_{-\theta,\theta}}(\x) \right)
 \right\}
\nonumber
 \\
&=&\frac{1}{2}\left(e^{-\bm{e}_{\widetilde{\nu}_{-\theta,\theta}}(\x)\theta_s}
g^{(+)}_{-\theta,\theta}(\x) + e^{\bm{e}_{\widetilde{\nu}_{-\theta,\theta}}(\x)\theta_s}
g^{(-)}_{-\theta,\theta}(\x)\right)
,
\end{eqnarray}
where  $\Lambda_{\theta_{s}}$ is denoted the phase shift operator.
\end{definition}

\begin{proposition}{The $\theta$-monogenic phase-shifted signal.}\\
For a real signal $g(\x)$ the
$\theta$-monogenic extension of the 
phase shifted real image if $\widetilde{\nu}_{-\theta,\theta}(\x)$ is varying
sufficiently slowly, {\em i.e.} assuming $\widetilde{\nu}_{-\theta,\theta}(\x)=
\widetilde{\nu}_{-\theta,\theta}(\x_0)$ for all $\x$ such that $
|\widetilde{g}^{(++)}_{-\theta,\theta}(\x)|$ is of non-negligible magnitude,
takes the form:
\begin{eqnarray}
\nonumber
(\Lambda_{\theta_{s}} g)_{\theta}^{(+)}(\x) &=&  |\widetilde{g}^{(++)}_{-\theta,\theta}(\x)| \ \left[
\cos{( 2\pi\widetilde{\phi}_{-\theta,\theta}(\x) - \theta_{s} ) }   + 
\bm{e}_{\widetilde{\nu}_{-\theta,\theta}}(\x) \ \sin{( 2\pi\widetilde{\phi}_{-\theta,\theta}(\x) - \theta_{s} ) }     \right]   \nonumber
\\
&=& e^{  - \bm{e}_{\widetilde{\nu}_{-\theta,\theta}}(\x) \theta_{s} }\  
\widetilde{g}^{(+)}_{-\theta,\theta}(\x).
\label{ShiftedCWT}
\end{eqnarray}
\end{proposition} 
If $\bm{e}_{\widetilde{\nu}_{-\theta,\theta}}(\x)=\bm{e}_{\widetilde{\nu}_{-\theta,\theta}},$ {\em
i.e.} the orientation is constant across the width of the
image, then as the Fourier transform corresponds to a right-hand multiplication,
we may note
\begin{eqnarray}
{\mathcal{F}}\left\{(\Lambda_{\theta_{s}} g_{-\theta,\theta}^{(+)}(\x) \right\}&=&
e^{- \bm{e}_{\widetilde{\nu}_{-\theta,\theta}} \theta_{s} }\widetilde{G}^{(+)}_{-\theta,\theta}(\f).
\end{eqnarray}

\begin{theorem}{MCWT of a phase-shifted signal.}\\
The MCWT of a phase-shifted signal, where the orientation
of the signal is constant over the width of the wavelet
is given by
\begin{equation}
\label{phassy}
w_{\psi}^{(\pm)}(\bm{\xi};\Lambda_{\theta_{s}} g)=e^{\mp  
2\pi\bm{e}_{\widetilde{\nu}_{\theta}}(\bm{b})\theta_s}
w_{\psi}^{(\pm)}(\bm{\xi};g).
\end{equation}
\end{theorem}
\begin{proof}
See section \ref{CWTINV}.
\end{proof}
This result implies that the phase shift between two images that locally
correspond to plane waves may easily be determined by equation (\ref{phassy}).

\begin{corollary}{Magnitude invariance of MCWT under phase-shifts.}\\
The magnitude of the MCWT of a phase-shifted signal, when the local orientation
of the signal is stable over the width of the wavelet is equivalent to that
of the not phase-shifted signal.
\end{corollary}
\begin{proof}
From equation (\ref{phassy}) we may note that $
w_{\psi,\theta_s}^{(\pm)}(\bm{\xi}; g)=e^{2\pi\bm{e}_{\nu_{\theta}}(\bm{b})\theta_s}
w_{\psi}^{(\pm)}(\bm{\xi};g),$ and thus it follows that
\begin{equation}
\left|w_{\psi}^{(\pm)}(\bm{\xi};\Lambda_{\theta_{s}} g)\right|^2
=\left|w_{\psi}^{(\pm)}(\bm{\xi};g)\right|^2.
\end{equation}
\end{proof}
Hence the magnitude of the MCWT is invariant to shifts of phase.
This implies that there will be no migration of energy across scales due
to misalignment between the wavelet and the signal. This shows the stability
of the MCWT to misalignment between the observed signal
and the analysis filters.

\section{Conclusions}
In this paper we have defined and provided suitable 2-D extensions of 1-D
analytic
wavelet decompositions that we denoted hyperanalytic.
We
have constructed classes of 2-D mother wavelets,
as well as specific examples of
such functions that may be used to calculate hyperanalytic
decomposition coefficients,
and
discussed properties of the hyperanalytic wavelets and wavelet coefficients in detail.
We have stressed the importance of magnitude invariance of the coefficients
to phase-shifts of the image,
and have constructed wavelets yielding coefficients that may be thought of as `locally
hyperanalytic' to exhibit such invariance. 

We argued that the design of 2-D hyperanalytic decompositions in continuous
space must start by understanding 2-D
hyperanalytic signals. Hyperanalytic signals are constructed as the appropriate
limit of hyperanalytic functions. A hyperanalytic function 
is a set of functions satisfying any generalisation of the Cauchy-Riemann
system. When considering analysis in dimensions higher than 1-D there is
in general more than one such system, and so a choice of the hyperanalytic
system
must be made. We considered two different systems in 2-D, namely the monogenic
and hypercomplex systems. We carefully
distinguished
between hyperanalytic signals (limits of hyperanalytic functions) and hyperanalytic
functions in our discussion.

Once a choice of hyperanalytic has been made, the hyperanalytic signal had
to be
suitably merged into the CWT to yield a hyperanalytic
CWT set of coefficients.
Two additional problems of interest were then resolved. Firstly the
rotation operator
does not commute
with the convolutions that combine to construct a hyperanalytic
signal from a real-valued image. To
determine the structure of rotated hyperanalytic objects in a suitable
frame of reference, we defined the
notion
of $\theta$-hyperanalytic as
a function that is a hyperanalytic signal in a rotated frame. This definition
allows for the interpretation of the hyperanalytic wavelet coefficients in
terms of local variational structure in a given orientation.
Furthermore the monogenic and hypercomplex
wavelet transforms were
related to the local UQFT (a unidirectional oscillatory representation), as well as
the local QFT (a separable oscillatory representation), of an image. We demonstrated
that our choices of hyperanalytic
wavelet
decompositions thus correspond to
the representation of local image structure as either a plane wave (locally 1-D
structure) and a local version of the UQFT basis elements, or as
a combination of locally separable structures with a natural
axis of observation and a local version of the QFT basis elements.

Given a hyperanalytic signal is a three or four vector-valued object it admits
a representation as a quaternionic function. The local energy of the real
signal was described
by 
the magnitude of the quaternion, whilst its local structure was described
by some suitable phase function/s defined in terms of the quaternionic elements,
thus parameterising the local structure
of the image in terms of
local magnitude,
scales and orientation of variation. 
Different polar descriptions were chosen for each of the two $\theta$-hyperanalytic
constructions, as each hyperanalytic definition is suited to particular forms
of local signal variations, and the interpretation of the magnitude and phases
were discussed.

Secondly, as we demonstrated, a decomposition constructed using a hyperanalytic
wavelet does not in general, even for zero rotation of the mother wavelet,
yield wavelet
coefficients that are hyperanalytic extensions of real-valued wavelet coefficients.
In fact the wavelet coefficients need not be
the real wavelet decomposition of a signal that has been extended
to a hyperanalytic object. We introduced the concept of an hyperanalyticizing
wavelet,
as a function constructing hyperanalytic wavelet coefficients. In 1-D the
analytic wavelets
localise the signal
in scale
and time simultaneously as yielding a local interpretable description of the signal in
terms of the local amplitude and the phase of the signal at that scale. As
both the operation of localisation and construction of an analytic object
commute in 1-D (both are multiplications in terms of complex objects) this simultaneous operation is not problematic. This, we demonstrated,
is not the case in 2-D (quaternionic objects in general do not commute and,
in addition, the rotation and the construction of hyperanalytic images do
not commute). 
We showed
that
the hyperanalytic coefficients were equivalent to either 
constructing the real wavelet decomposition of the real image and then calculating
the $\theta$-hyperanalytic extension of the wavelet coefficients or to constructing
the $\theta$-hyperanalytic
signal and then constructing the real wavelet coefficients of this object.
This gave the interpretation of the hyperanalytic wavelet coefficients, as
 the $\theta$-hyperanalytic signal of the local projection or equivalently
the local projection of the $\theta$-hyperanalytic component.
To
be able to prove these properties we defined
the decomposition of any real-valued image into $\theta$-hyperanalytic
and $\theta$-anti-hyperanalytic
objects, in direct analogue with the analytic/anti-analytic decomposition
in 1-D
\cite{Papoulis}.

We demonstrated that the hyperanalyticizing
wavelet either eliminates, or converts, the anti-hyperanalytic objects into
hyperanalytic objects. The hyperanalytic decomposition allows
for the determination of the CWT coefficients of phase-shifted
images, where the elimination/conversion of the anti-components is pivotal
for easy parameterisation of the coefficients of the phase shifted signal in terms of the coefficients of the original signal. 
For both the monogenic and hypercomplex coefficients we retrieved magnitude invariance
of the wavelet coefficients under phase-shifts, even if in the hypercomplex
case some restrictions had to be placed on the observed signal. Magnitude invariance
to phase shift is an important property when discretizing the implementation:
the local representation will be sampled and so small shifts in space should
not change the format of the representation substantially.
In the monogenic case
we additionally demonstrated local invariance
of magnitude under rotations with certain classes of mother wavelets and
proved
directional selectivity of the quaternionic
coefficients. The stability and other properties of the transform depend on constructing wavelet
coefficients that are hyperanalytic: this justifies our careful development
of the hyperanalytizing wavelets.

In 1-D, analytic wavelets are also used for the
decomposition
and characterisation of non-stationary oscillatory signals, 
using wavelet ridge techniques \cite{Delprat,mallat}.
Usage of a localisation method is necessary
for the analysis of multi-component non-stationary oscillatory signals,
as for such signals the phase and amplitude description of the local structure
in the spatial domain is not informative \cite{Cohen}. 
In 2-D we showed that wavelet ridges for unidirectional and separable oscillatory
signals can be determined from the monogenic
and
hypercomplex wavelet coefficients. The advantage of using quaternionic wavelets
is that separable oscillations can easily be parameterised in both directions
simultaneously, and the analysis
of unidirectional oscillations is simplified as once components have been
separated in phase space the orientation of the variation can be determined
from the monogenic wavelet coefficients
without calculating the transform at all values of $\theta.$
Our work thus complements existing ridge analysis using complex wavelets
in 2-D.
Given all the components of the quaternionic wavelet can be argued to have
the same spatial and spatial frequency localisation only a few of the wavelet
coefficients of most deterministic image features will be of non-negligible magnitude, whilst most noise spreads across all coefficients.
Thus the localisation also
increases the signal-to-noise ratio of the signal in the coefficients that
contain most of the signal presence, and this further facilitates analysis
of noisy realisations.

Finally, even if general results are of interest, to be able to implement
the CWT, explicit examples of mother wavelets 
must be provided for analysis. Existing discrete hyperanalyticizing wavelets \cite{Chan}
and \cite{Selesnick2005}
were given an interpretation in terms of the general framework developed
in this paper.
Furthermore we provided families
of continuous wavelets that are hyperanalyticizing. We chose to base the
wavelets
on the 1-D and 2-D generalized Morse wavelets, that can be shown to enjoy
optimal localisation in phase space: we provided five classes of wavelets,
the hypercomplex Morse wavelet, the separable and isotropic hypercomplexing Morse
wavelets,
as well as the the directional
and isotropic monogenic Morse wavelets.
 
The study and development of analytic 1-D filters have 
lead to the synthesis of a wealth of signal processing methodology
and applications. Before the discrete implementation was developed, continuous
analytic wavelet decompositions using analytic wavelets were introduced and championed,
such as wavelet ridge methods, and local decompositions of multivariate time
series.
This paper demonstrates the potential for hyperanalytic (quaternionic) methods,
and the representation of local 2-D structure in terms of quaternionic objects.
It
is anticipated that methods based on quaternionic wavelet coefficients hold
many
novel challenges and potential future developments.

\section*{Acknowledgments}
SO \& GM would like to thank Dr Frederik Simons for valuable discussions
and SO would like to thank Professor Andrew Walden for introducing her to
this
research area. GM gratefully acknowledges financial support from EPSRC (UK).



\renewcommand{\theequation}{\thesection.\arabic{equation}}
\renewcommand{\theequation}{A-\arabic{equation}}
\setcounter{equation}{0}  

\appendices
\section{Hyperanalytic Properties}
\subsection{QFT of $g^{(++)}_{\theta}(\x)$ \label{qftofg++}}
Note that $G(\f)$ commutes with $\bm{j}$ as it is a function of $\bm{j}$
only. We obtain:
\begin{eqnarray}
{\cal F}_{\cal Q}\left\{g^{(++)}_{\theta}(\x)\right\}
&=&{\cal F}_{\cal Q}\left\{g(\bm{r}_{-\theta}\x)\right\}+
\bm{i}{\cal F}_{\cal Q}\left\{g_H^{(1)}(\bm{r}_{-\theta}\x)\right\}+
{\cal F}_{\cal Q}\left\{g_H^{(2)}(\bm{r}_{-\theta}\x)\right\}\bm{j}+
\bm{i}{\cal F}_{\cal Q}\left\{g_H^{(3)}(\bm{r}_{-\theta}\x)\right\}\bm{j}
\nonumber \\
&=&\frac{1-\bm{k}}{2}G(\bm{r}_{-\theta}\q)+
\frac{1+\bm{k}}{2}G(\bm{r}_{-\theta}\bm{J}_{\pi/2}\q)+
\bm{i}\left(\frac{1-\bm{k}}{2}G_H^{(1)}(\bm{r}_{-\theta}\q)\right.
\nonumber \\
&&\left.+
\frac{1+\bm{k}}{2}G_H^{(1)}(\bm{r}_{-\theta}\bm{J}_{\pi/2}\q)
\right)+
\left(\frac{1-\bm{k}}{2}G_H^{(2)}(\bm{r}_{-\theta}\q)+
\frac{1+\bm{k}}{2}G_H^{(2)}(\bm{r}_{-\theta}\bm{J}_{\pi/2}\q)
\right)\bm{j}
\nonumber \\
&&+\bm{i}\left(\frac{1-\bm{k}}{2}G_H^{(3)}(\bm{r}_{-\theta}\q)+
\frac{1+\bm{k}}{2}G_H^{(3)}(\bm{r}_{-\theta}\bm{J}_{\pi/2}\q)
\right)\bm{j}
\nonumber \\
&=&\frac{1-\bm{k}}{2}\left(1+{\mathrm{sgn}}\left(
\left[\bm{r}_{-\theta}\q\right]_{1} \right)\right)\left(1+{\mathrm{sgn}}\left(
\left[\bm{r}_{-\theta}\q
\right]_{2} \right)\right)G(\bm{r}_{-\theta}\q)\nonumber\\
&&+\frac{1+\bm{k}}{2}\left(1+{\mathrm{sgn}}\left(
\left[\bm{r}_{\theta}\q\right]_{1} \right)\right)\left(1+{\mathrm{sgn}}\left(
\left[\bm{r}_{\theta}\q
\right]_{2} \right)\right)G(\bm{r}_{-\theta}\bm{J}_{\pi/2}\q)
\label{qffft}
\end{eqnarray}
Furthermore, if we had started by defining $\x^{\prime}=\bm{r}_{-\theta}\x$
and
then found the QFT in $\x^{\prime},$ with QFT variable $\q^{\prime} $ we
would have obtained
$
{\cal F}_{\cal Q}\left\{g_{\theta}^{(++)}\right\}(\q^{\prime})=
\left(1+ {\mathrm{sgn}}\left(
\left[\q^{\prime}\right]_{1} \right)\right)\left(1+{\mathrm{sgn}}\left(
\left[\q^{\prime}
\right]_{2} \right)\right) G_{Q}(\q^{\prime}),
$
and so in the rotated frame of reference the QFT is supported on positive
quaternionic frequency only. 

\subsection{The $\theta$-Monogenic Signal \label{monoe}}
With $p(\x,y)$ and $q_s(\x,y)$ defined by equation (\ref{ddversionhilbert}) let us denote
by operator notation the construction of a system defined from $g(\x)$ satisfying
the Riesz system of equations:
${\cal P}g(\x,y)=\left(p(\cdot,y) \ast \ast g(\cdot)\right)(\x),$ and
${\cal Q}_{s} g(\x,y)
=\left( q_s(\cdot,y) \ast \ast g(\cdot)\right)(\x).$
When $u_{R,g}(\x,y)=
{\cal P}g(\x,y),$ and $v^{(s)}_{R,g}(\x,y)
={\cal Q}_{s}g(\x,y),$ we denote
$k_{R;g}^+(\x,y)=u_{R,g}(\x,y) + \bm{i} \ v^{(1)}_{R,g}(\x,y) + \bm{j} \ v^{(2)}_{R,g}(\x,y)
=
{\cal P} \ g(\x) + \bm{i} \ {\cal Q}_{1} \ g(\x) + \bm{j}
\ {\cal Q}_{2} \ g(\x),$ which for $y \rightarrow 0^{+}$ in the
upper half-space
reduces to the monogenic function $g^{(+)}(\x) =g(\x)
+ \bm{i} \ g^{(1)}_R(\x) + \bm{j} \ g^{(2)}_R(\x).$ Rotate the solution to
form,
$R_{\theta} \ {\cal P} \ g(\x,y) + \bm{i} \ R_{\theta} \ {\cal
Q}_{1} \ g(\x,y) + 
\bm{j} \ R_{\theta} \ {\cal Q}_{2} \ g(\x,y) $.  Rotations require
some
attention; they commute with the convolution with the Poisson kernel $[{\cal
R}_{\theta},
{\cal P}]= 0 $, but not with operation of the convolution with the conjugate Poisson kernels
$[R_{\theta},{\cal
Q}_{j}] \neq 0$, see \cite{SteinWeiss}. In particular,
\begin{equation}
 \left( \begin{array}{c}
                  {R}_{\theta} \  {\cal Q}_{1} \ g(\x,y) \\
                  {R}_{\theta}\ {\cal Q}_{2} \ g(\x,y)
                  \end{array} \right) = \bm{r}_{-\theta} \
 \left( \begin{array}{c}
      {\cal Q}_{1} \ { R}_{\theta} g(\x,y)  \\
      {\cal Q}_{2} \ { R}_{\theta} \ g(\x,y)
      \end{array} \right) = \bm{r}_{-\theta} \
 \left( \begin{array}{c}
     {\cal Q}_{1} \ g(\bm{r}_{-\theta} \x,y) \\
     {\cal Q}_{2} \ g(\bm{r}_{-\theta} \x,y) 
 \end{array} \right).
\label{poissonrotations}
\end{equation}
As noted by  \cite{SteinWeiss1}, the Riesz system 
is invariant under dilation, translation, and simultaneous rotation of $(x_{1},x_{2})$
and $({\cal Q}_{1}g ,{\cal Q}_{2} \ g)$  by the same angle.  In other words, (\ref{poissonrotations}) is also a solution of the Riesz system, but this time
defined in variable $\x^{\prime}=\bm{r}_{-\theta}\x.$
For $y \rightarrow 0^{+}$, (\ref{poissonrotations}) reduces to
\begin{eqnarray}
&&   R_{\theta} \ g(\x) + 
\bm{ i} \ \left[ \cos{\theta} \ {\cal R}_{1} 
\  R_{\theta} \ g(\x)  
 + \sin{\theta} \ {\cal R}_{2} \  R_{\theta} \ g(\x) \right] \nonumber \\
 &&  +\bm{ j} \  \ \left[ - \sin{\theta} \ {\cal R}_{1} 
\  R_{\theta} \ g(\x)  + \cos{\theta} \ {\cal R}_{2} \  R_{\theta} \ g(\x) \right].
\end{eqnarray}
This can be recast as $R_{\theta} \ g(\x) + \bm{i} \ R_{\theta} \ {\cal R}_{1} \ g(\x) + 
\bm{j} \ R_{\theta} \ {\cal R}_{2} \ g(\x),$ see \cite{SteinWeiss} for the
algebraic properties of the Riesz kernels, and in analogy to (\ref{poissonrotations}),
\begin{equation}
 \left( \begin{array}{c}
                  {\cal R}_{\theta} \  {\cal R}_{1} \ g(\x) \\
                  {\cal R}_{\theta}\ {\cal R}_{2} \ g(\x)
                  \end{array} \right) = \bm{r}_{-\theta} \
 \left( \begin{array}{c}
      {\cal R}_{1} \ {\cal R}_{\theta} \ g(\x)  \\
      {\cal R}_{2} \ {\cal R}_{\theta} \ g(\x)
      \end{array} \right) = \bm{r}_{-\theta} \
 \left( \begin{array}{c}
     {\cal R}_{1} \ g(\bm{r}_{-\theta} \x) \\
     {\cal R}_{2} \ g(\bm{r}_{-\theta} \x) 
 \end{array} \right).
\label{rietzrotations}
\end{equation}
Therefore, if a given function $g^{(+)}(\x)$ is a monogenic signal, {\em
i.e.}
the limit of a monogenic function, then the function
generated by rotating its axis by $\theta$ is also a monogenic signal, {\em i.e.} it is the limit
as $y\rightarrow 0^{+}$ of a system of functions, satisfying the Riesz system
in a rotated frame of reference. A key component in the functionality of
the set-up is that $\left|\bm{r}_{-\theta}\x\right|=\left|\x\right|.$ Hence
the $\theta$-monogenic signal is aptly named. {\em Mutatis mutandis} the
proofs extend to the anti-monogenic case.
\section{Properties of Quaternionic Wavelets}
\subsection{FT and QFT of the Monogenic Wavelets \label{ftqftmon}}
We note that
\begin{eqnarray}
\Psi_{\bm{\xi},Q}^{(+)}&=&a e^{-2i\pi \bm{i} q_1 b_1}
\left[{\cal F}_Q R_{\theta}\left\{\psi^{(r)}+\bm{i}\psi^{(1)}_R+
\psi^{(2)}_R\bm{j}
\right\}(a\q)\right]e^{-2\pi \bm{j} q_2 b_2},
\end{eqnarray}
as each real wavelet is transformed according to:
\begin{eqnarray}
\nonumber
\Psi_{\bm{\xi};Q}^{(r)}(\q)
&=&{\cal{F}}_Q U_{\bm{\xi}} \left\{\psi\right\}(\q)\\
\nonumber
&=&\int_{-\infty}^{\infty}\int_{-\infty}^{\infty} e^{-2\pi\bm{i}q_1x_1}\frac{1}{a}\psi^{(r)}\left(
\bm{r}_{-\theta}\left(\frac{\x-\be}{a}\right)\right)e^{-2\pi\bm{j}q_2x_2}
\;d^2\x\\
\nonumber
&=&\int_{-\infty}^{\infty}\int_{-\infty}^{\infty} e^{-2\pi\bm{i}(q_1b_1+aq_1y_1)}a\psi^{(r)}\left(
\bm{r}_{-\theta}\y \right)e^{-2\pi\bm{j}(q_2b_2+q_2ay_2)}
\;d^2\y\\
\label{localreppy}
&=&a e^{-2\pi\bm{i}q_1b_1}{\cal F}_Q\left\{\psi^{(r)}\left(
\bm{r}_{-\theta}\x \right)\right\}(a\q)e^{-2\pi\bm{j}q_2b_2}.
\end{eqnarray}
Note that from \cite{HahnSnopek} for $\psi(\x)\in{\mathbb{R}}$
\begin{equation*}
\Psi_{Q}(\q)=\frac{1-\bm{k}}{2}\Psi(\q)+\frac{1+\bm{k}}{2}\Psi((-q_1,q_2)),
\end{equation*}
and we shall use this property on numerous occasions.
We find with $\widetilde{\psi}^{(r)}_{\theta}(\x)=R_{\theta}\left\{\psi^{(r)}\right\}(\x)$
that 
\begin{eqnarray*}
{\cal F}_Q R_{\theta}\left\{\psi^{(r)}\right\}(\q)
={\cal F}_Q\left\{\widetilde{\psi}^{(r)}\right\}(\q)=\frac{1-\bm{k}}{2}\widetilde{\Psi}_{\theta}^{(r)}(\q)+
\frac{1+\bm{k}}{2}\widetilde{\Psi}_{\theta}^{(r)}(\bm{J}_{\pi/2}\q).
\end{eqnarray*}
We also note that
$\widetilde{\Psi}_{\theta}^{(r)}(\q)=\Psi^{(r)}(\bm{r}_{-\theta}\q).$
Thus we find
\begin{eqnarray}
\nonumber
{\cal F}_Q R_{\theta}\left\{\psi^{(r)}\right\}(\q)
&=&\frac{1-\bm{k}}{2}\Psi^{(r)}(\bm{r}_{-\theta}\q)+
\frac{1+\bm{k}}{2}\Psi^{(r)}(\bm{r}_{-\theta}\bm{J}_{\pi/2}\q)\\
\nonumber
&=&\frac{1-\bm{k}}{2}\left[\frac{1+\bm{k}}{2}\Psi_Q^{(r)}(\bm{r}_{-\theta}\q)+
\frac{1-\bm{k}}{2}\Psi_Q^{(r)}(\bm{J}_{\pi/2}\bm{r}_{-\theta}\q)
\right]\\
\nonumber
&&+\frac{1+\bm{k}}{2}\left[\frac{1+\bm{k}}{2}\Psi_Q^{(r)}(\bm{r}_{-\theta}\bm{J}_{\pi/2}\q)+
\frac{1-\bm{k}}{2}\Psi_Q^{(r)}(\bm{J}_{\pi/2}\bm{r}_{-\theta}\bm{J}_{\pi/2}\q)
\right]\\
\nonumber
&=&\frac{1}{2}\Psi_Q^{(r)}(\bm{r}_{-\theta}\q)-\frac{\bm{k}}{2}
\Psi_Q^{(r)}(\bm{J}_{\pi/2+\theta/2}\q)+\frac{\bm{k}}{2}\Psi_Q^{(r)}(\bm{J}_{\pi/2-\theta/2}\q)
\\
&&+\frac{1}{2}\Psi_Q^{(r)}(\bm{r}_{\theta}\q).
\label{lastref}
\end{eqnarray}
Similarly to $\Psi^{(r)}_{\bm{\xi},Q}(\q )$ we can now derive the QFTs of the members of the families produced from the first and second Riesz components of the monogenic mother wavelet, $\Psi^{(1)}_{R,\bm{\xi},Q}(\q )$
and $\Psi^{(2)}_{R,\bm{\xi},Q} (\q )$, as these two functions are
real, see also \cite{HahnSnopek}. As before, the crucial step is to calculate
the QFT of the rotated real wavelets, by combining equation (\ref{lastref})
with equation (\ref{FTRiesz}). This yields:
\begin{eqnarray}
{\cal F}_Q R_{\theta}\left\{\psi^{(1)}_R\right\}(\q)
&=&\left(-\frac{\bm{i}}{2}\right)\cos(\eta-\theta)
\left[\Psi_Q^{(r)}(\bm{r}_{-\theta}\q)-\bm{k}\Psi_Q^{(r)}(\bm{J}_{\pi/2+\theta/2}\q)
\right] \nonumber \\
&+& \left(-\frac{\bm{i}}{2}\right)\cos(\eta+\theta)
\left[\bm{k}\Psi_Q^{(r)}(\bm{J}_{\pi/2-\theta/2}\q)+\Psi_Q^{(r)}(\bm{r}_{\theta}\q)
\right],
\label{lastref1}
\end{eqnarray} 
\begin{eqnarray}
{\cal F}_Q R_{\theta}\left\{\psi^{(2)}_R\right\}(\q)
&=& \frac{1}{2}\left[
\Psi_Q^{(r)}(\bm{r}_{-\theta}\q)
-\bm{k}
\Psi_Q^{(r)}(\bm{J}_{\pi/2+\theta/2}\q)\right]
(-\bm{j})\sin(\eta-\theta)
 \nonumber  \\ 
&+& \frac{1}{2}
\left[\bm{k}\Psi_Q^{(r)}(\bm{J}_{\pi/2-\theta/2}\q)
+\Psi_Q^{(r)}(\bm{r}_{\theta}\q)\right]
(-\bm{j})\sin(\eta+\theta).
\label{lastref2}
\end{eqnarray} 
Putting equations (\ref{localreppy}), (\ref{lastref}), (\ref{lastref1}),
and (\ref{lastref2}) together, we get the QFT of a member
of the family generated through translations, dilations, and rotations from
the entire mother monogenic wavelet:
\begin{eqnarray}
\Psi_{\bm{\xi},Q}^{(+)}(\bm{q}) &=&\frac{a}{2} e^{-2\pi \bm{i} q_1 b_1}
\left[\left(1+\cos(\eta-\theta)+\sin(\eta-\theta)
\right)\left(
\Psi_Q^{(r)}(\bm{r}_{-\theta}\q)-\bm{k}\Psi_Q^{(r)}(\bm{J}_{\pi/2+\theta/2}\q)
\right)
\right. \nonumber \\
 &+&\left.
\left(1+\cos(\eta+\theta)+\sin(\eta+\theta)
\right)\left(\Psi_Q^{(r)}(\bm{r}_{\theta}\q)+
\bm{k}\Psi_Q^{(r)}(\bm{J}_{\pi/2-\theta/2}\q)
\right)\right]e^{-2\pi \bm{j} q_2 b_2}.
\label{equationtolabel}
\end{eqnarray}
 
\subsection{The Wavelet Family of the Isotropic Monogenic Mother
Wavelet \label{isofamily}}
From equation (\ref{FTanalyticwavelet}) we may note
that the FT of a translated, dilated, and rotated monogenic wavelet which
originates in a monogenic mother wavelet with radially symmetric real part
respectively is given by:
\begin{eqnarray}
\Psi^{(+)}_{\bm{\xi}}( \f ) &=&  {\cal F} \ \psi^{(+)}_{\bm{\xi}} (\x)
\nonumber \\ 
&=& 
\label{isoFourier}
\left\{ 1 + \ \left[ - \bm{k} \cos{(\phi -  \theta)} +  \sin{(\phi -  \theta)} \right]
\right\} \ a \Psi^{(e)}(af) \ e^{-\bm{j} 2\pi \f^T \be}, \\
\nonumber
\psi^{(+)}_{\bm{\xi}}( \f )&=&\int_{-\infty}^{\infty}\int_{-\infty}^{\infty} \left\{ 1 + \ \left[ - \bm{ij} \cos{(\phi -  \theta)} +  \sin{(\phi -  \theta)} \right]
\right\} \ a \Psi^{(e)}(af) \ e^{-\bm{j} 2\pi \f^T (\be-\x)}\;d^2\x\\
&=&\psi^{(e)}_{\bm{\xi}_0}(\x)+\bm{i}\left[\cos(\theta)\psi^{(1)}_{R,\bm{\xi}_0}(\x)
+\sin(\theta)\psi^{(2)}_{R,\bm{\xi}_0}(\x)\right]+
\bm{j}\left[-\sin(\theta)\psi^{(1)}_{R,\bm{\xi}_0}(\x)
+\cos(\theta)\psi^{(2)}_{R,\bm{\xi}_0}(\x)\right]. \nonumber \\
\label{quatisofourier}
\end{eqnarray}
This gives the stated results. For the QFT expression, using equation (\ref{QFTanalyticwavelet}),
we note
from equation (\ref{equationtolabel}), recalling that the QFT of an isotropic
real wavelet
is a real-valued object:
\begin{eqnarray*}
\Psi_{Q,\bm{\xi}}^{(+)}&=&\frac{a}{2} e^{-2\pi \bm{i} q_1 b_1}
\left[\left(1+\cos(\eta-\theta)+\sin(\eta-\theta)
\right)\left(
\Psi_Q^{(e)}(q)-\bm{k}\Psi_Q^{(e)}(q)
\right)
\right.\\
\nonumber
&&+\left.
\left(1+\cos(\eta+\theta)+\sin(\eta+\theta)
\right)\left(\Psi_Q^{(e)}(q)+
\bm{k}\Psi_Q^{(e)}(q)
\right)\right]e^{-2\pi \bm{j} q_2 b_2}\\
&=&a e^{-2\pi \bm{i} q_1 b_1}
\Psi_Q^{(e)}(q)\left(1+\cos(\theta)\left(
\cos(\eta)+\sin(\eta)\right)+\bm{k}\sin(\theta)\left(-
\sin(\eta)+\cos(\eta)\right)
\right)e^{-2\pi \bm{j} q_2 b_2}.
\end{eqnarray*}
This yields the desired result.
%
\subsection{Directional Monogenic Wavelets \label{dirrictly}}
We recall that 
\begin{equation}
\label{analyticparity}
\Psi_{1D}^{(o)}(f)=-\bm{j}\Psi_{1D}^{(e)}(f),\quad{\mathrm{if}} \quad f>0,
\quad
\Psi_{1D}^{(o)}(f)=\bm{j}\Psi_{1D}^{(e)}(f),\quad{\mathrm{if}} \quad f>0.
\end{equation}
Direct calculation yields the result
\begin{eqnarray}
{\cal R}_1\left\{\psi_{D}\right\}(\x)&=&
\int_{-\infty}^{\infty} \int_{-\infty}^{\infty} 
\frac{f_1}{\bm{j} f}\left[\Psi^{(e)}_{1D}\left(\frac{f_1-f_2}{\sqrt{2}}\right)
\Psi^{(e)}_{1D}\left(\frac{f_1+f_2}{\sqrt{2}}\right)-
\Psi^{(o)}_{1D}\left(\frac{f_1-f_2}{\sqrt{2}}\right)\Psi^{(o)}_{1D}\left(\frac{f_1+f_2}{\sqrt{2}}\right)\right]
 \nonumber \\ 
&&\times e^{2\bm{j}\pi \f^T \x}\;d^2\f
\nonumber \\
 &=&4\int^{\infty}_{0} \int_0^{\infty}\frac{g_1+g_2}{\sqrt{2}g}\left[\Psi^{(e)}_{1D}\left(g_1\right)
 \Psi^{(e)}_{1D}\left(g_2\right)\right]
 \sin(2\pi \g^T \bm{r}_{\pi/4} \x)\;d^2\g \nonumber \\
&=&4\int_{0}^{\infty} \int_0^{\infty}\frac{\sqrt{2}g_1}{g}\left[\Psi^{(e)}_{1D}\left(g_1\right)
\Psi^{(e)}_{1D}\left(g_2\right)\right] 
\sin(2\pi x_1(\frac{g_1+g_2}{\sqrt{2}}))
\;d^2\g.
\label{riessyresult}
\end{eqnarray}
Thus ${\cal R}_1\left\{\psi_D\right\}(\x)$ is constant in $x_2$ and
odd in $x_1,$ and a highly directional wavelet. 
We can in a similar fashion calculate the second Riesz transform to retrieve.
\begin{eqnarray}
\nonumber
{\cal R}_2\left\{\psi_D\right\}(\x)&=&\frac{1}{\sqrt{2}}
\left[\psi^{(\pi/4,1)}_D(\bm{r}_{\pi/4}\x)-
\psi^{(\pi/4,2)}_D(\bm{r}_{\pi/4}\x)\right] = 0.
\label{Riesz2dir}
\end{eqnarray}
Furthermore note that:
\begin{eqnarray*}
\nonumber
\psi_D(\x)&=&
\int_{-\infty}^{\infty} \int_{-\infty}^{\infty} 
\left[\Psi^{(e)}_{1D}\left(\frac{f_1-f_2}{\sqrt{2}}\right)
\Psi^{(e)}_{1D}\left(\frac{f_1+f_2}{\sqrt{2}}\right)-
\Psi^{(o)}_{1D}\left(\frac{f_1-f_2}{\sqrt{2}}\right)\Psi^{(o)}_{1D}\left(\frac{f_1+f_2}{\sqrt{2}}\right)\right]
\\\nonumber
&&\times e^{2\bm{j}\pi \f^T \x}\;d^2\f\\
\nonumber
 &=&4\int^{\infty}_{0} \int_0^{\infty}\left[\Psi^{(e)}_{1D}\left(g_1\right)
 \Psi^{(e)}_{1D}\left(g_2\right)\right]
 \cos(2\pi \g^T \bm{r}_{\pi/4} \x)\;d^2\g.
\label{Riesz1direven}
\end{eqnarray*}
\subsection{Exact form of the spatial representation of the isotropic Morse
Wavelets \label{spatialmorse}}
To be able to calculate the spatial domain representation of $\psi_{I}^{(+)}(\x),$
given by equation (\ref{morsemonogenic})
in the instance of the Morse wavelets, we make use of the following integrals, see \cite{McLachlan}, 
\[
\int_{0}^{2 \pi}  d\phi \ e^{\bm{j} z \cos{(\theta - \phi)}} \sin{m \phi} = 2 \pi \bm{j}^{m} J_{m}(z) \sin{m \theta}, \quad \int_{0}^{2 \pi} d\phi \ e^{\bm{j} z \cos{(\theta - \phi)}} \cos{m \phi} = 2 \pi \bm{j}^{m} J_{m}(z) \cos{m \theta}. \]
In special cases  $n=0$, $l>-2$, $m=1,2$ the analytic forms corresponding to equations (\ref{morse})
and (\ref{morsemonogenic}) can be found. 
Note that ${\rm L}^{c}_{0}(x) = 1$. 
For $m=1$ and $n=0$ we get
\begin{equation}
\nonumber
\psi_{0;l,1}^{(e)}(x) =  \frac{A_{0,l,1}^{\prime}}{2 \pi}  \frac{1}{[1 +  x^2]^{\frac{l+2}{2}}}  \ \Gamma(l+2) \  
P^{0}_{l+1}(\frac{1}{ \sqrt{1 + x^2} } )
\end{equation} 
where we have used equation 6.621.1 of \cite{GR} and $P^{\mu}_{\nu}$ are
the Legendre functions of degree 
$\nu$ and order $\mu$, see e.g \cite{AS}. The monogenic wavelet is
\begin{eqnarray*}
\psi_{0;l,1}^{(+)}(\x) &=& \psi_{0;l,1}^{(e)}(x) + \ \left( \bm{i} \
\cos{\chi} +  \bm{j} \ \sin{\chi} \right) \frac{A_{0,l,1}^{\prime}}{ 2 \pi} \nonumber \\
&\times &   \frac{1}{[1 + x^2]^{\frac{l+2}{2}}} 
\ 
\Gamma(l+3) \ P^{-1}_{l+1}(\frac{1}{ \sqrt{1 + x^2} } ).   
\end{eqnarray*}
For $m=2$ we get:
\begin{equation} 
\nonumber
\psi_{0;l,2}^{(e)}(x) =  \frac{A_{0,l,2}^{\prime}}{2\pi} \frac{1}{2} \Gamma(\frac{l+2}{2})\  _{1}F_{1}(\frac{l+2}{2};1;-\frac{x^2}{4})
\end{equation}
where we have used equation 6.631 of \cite{GR}. For $l=2$, $s=1/2$, this wavelet
 is identical with the mexican hat, see e.g. \cite{Antoine3}. The monogenic wavelet is:
\begin{eqnarray*}
\psi_{0;l,2}^{(+)}(\x) &=& \psi_{0;l,2}^{(e)}(x) +   \ 
\left( \bm{i} \ \cos{\chi} + \bm{j} \ \sin{\chi} \right) \frac{A_{0,l,2}^{\prime}}{2 \pi} \nonumber \\ 
&\times & \frac{x}{4} 
\Gamma(\frac{l+3}{2})\  _{1}F_{1}(\frac{l+3}{2};2;\frac{-x^2}{4}).
\end{eqnarray*}

\section{Properties of the Quaternionic Wavelet Coefficients}
\subsection{QFT of the wavelet coefficients with an arbitrary quaternionic wavelet
\label{nykonstig}}
The CWT with respect to a quaternionic wavelet is written in terms of its real components and calculate the QFT component
 by component as follows:
\begin{eqnarray}
{\cal F}_Q\left\{w_{\psi}(\bm{\xi};g)\right\}
&=&{\cal F}_Q\left\{w_{\psi}^{(r)}(\bm{\xi};g)\right\}-
\bm{i}{\cal F}_Q\left\{w_{\psi}^{(i)}(\bm{\xi};g)\right\}
-{\cal F}_Q\left\{w_{\psi}^{(j)}(\bm{\xi};g)\right\}\bm{j}-
\bm{i}{\cal F}_Q\left\{w_{\psi}^{(k)}(\bm{\xi};g)\right\}\bm{j}
\nonumber \\
 &=&\frac{1}{2}W_{\psi}(\bm{\kappa};g)-\frac{\bm{k}}{2}W_{\psi^{(r)}-i\psi^{(i)}+j\psi^{(j)}-k\psi^{(k)}}(\bm{\kappa};g)
+\frac{1}{2}W_{\psi}(\bm{J}_{\pi/2}\bm{\kappa};g) \nonumber \\
&+&\frac{\bm{k}}{2}W_{\psi^{(r)}-i\psi^{(i)}+j\psi^{(j)}-k\psi^{(k)}}(\bm{J}_{\pi/2}\bm{\kappa};g)
\label{qftwav1} \nonumber \\
&\neq & \frac{1-\bm{k}}{2}W_{\psi}(\bm{\kappa};g)+
\frac{1+\bm{k}}{2}W_{\psi}(\bm{J}_{\pi/2}\bm{\kappa};g).
\end{eqnarray}
\subsection{The Fourier Transform of the Hypercomplex wavelet transform of
the $\theta$-Hypercomplex signal \label{konstig}}
To establish the properties of the wavelet transform, we find the Fourier
decomposition of the real signal in terms of the $\theta$-hypercomplex and
anti-hypercomplex components.
We note that for $\mu_1=\pm,$ and $\mu_2=\pm,$
\begin{eqnarray*}
w^{(\mu_1 \mu_2)}_{\psi}(\bm{\xi};\widetilde{g}^{(++)}_{-\theta^{\prime},\theta^{\prime}})
&=&
w^{(\mu_1 \mu_2)}_{\psi}(\bm{\xi};\widetilde{g}_{-\theta^{\prime},\theta^{\prime};H})+\bm{i}
w^{(\mu_1 \mu_2)}_{\psi}(\bm{\xi};\widetilde{g}_{-\theta^{\prime},\theta^{\prime};H}^{(1)})+\bm{j}
w^{(\mu_1 \mu_2)}_{\psi}(\bm{\xi};\widetilde{g}_{-\theta^{\prime},\theta^{\prime};H}^{(2)})  \\
&+& \bm{k}
w^{(\mu_1 \mu_2)}_{\psi}(\bm{\xi};\widetilde{g}_{-\theta^{\prime},\theta^{\prime};H}^{(3)}).
\end{eqnarray*}
Of course for any real-valued function $h(\cdot)$ by the linearity of the
wavelet transform:
\begin{eqnarray*}
w^{(\mu_1 \mu_2)}_{\psi}(\bm{\xi};h)&=&
w^{(e e)}_{\psi}(\bm{\xi};h)-\mu_1 \bm{i}w^{(o e)}_{\psi}(\bm{\xi};h)
-\mu_2\bm{j}w^{(e o)}_{\psi}(\bm{\xi};h)+
\mu_1\mu_2\bm{k}w^{(o o)}_{\psi}(\bm{\xi};h).
\end{eqnarray*}
Thus
\begin{eqnarray}
\nonumber
W^{(++)}_{\psi}(\bm{\zeta};h)&=&
W^{(ee)}_{\psi}(\bm{\zeta};h)-\bm{i}W^{(ee)}_{\psi}(\bm{\zeta};h)\bm{j}
{\mathrm{sgn}}\left(\left[\bm{r}_{-\theta}\f_b\right]_1\right)
-\bm{j}W^{(ee)}_{\psi}(\bm{\zeta};h)\bm{j}
{\mathrm{sgn}}\left(\left[\bm{r}_{-\theta}\f_b\right]_2\right)\\
\nonumber
&&
+\bm{k}W^{(ee)}_{\psi}(\bm{\zeta};h)(-1)
{\mathrm{sgn}}\left(\left[\bm{r}_{-\theta}\f_b\right]_1\right)
{\mathrm{sgn}}\left(\left[\bm{r}_{-\theta}\f_b\right]_2\right)\\
\nonumber
&=&\left[1+{\mathrm{sgn}}\left(\left[\bm{r}_{-\theta}\f_b\right]_2\right)
-\bm{k}{\mathrm{sgn}}\left(\left[\bm{r}_{-\theta}\f_b\right]_1\right)
-\bm{k}{\mathrm{sgn}}\left(\left[\bm{r}_{-\theta}\f_b\right]_1\right)
{\mathrm{sgn}}\left(\left[\bm{r}_{-\theta}\f_b\right]_2\right)
\right]W^{(ee)}_{\psi}(\bm{\zeta};h)\nonumber\\
&=&\left[1+{\mathrm{sgn}}\left(\left[\bm{r}_{-\theta}\f_b\right]_2\right)\right]
\left[1-\bm{k}{\mathrm{sgn}}\left(\left[\bm{r}_{-\theta}\f_b\right]_1\right)\right]
W^{(ee)}_{\psi}(\bm{\zeta};h).
\end{eqnarray}
Thus the FT of the CWT of the hypercomplex image with respect to the hypercomplexing
wavelet is given by:
\begin{eqnarray}
W_{\psi}^{(++)}(\bm{\zeta};\widetilde{g}^{(\mu_1 \mu_2)}_{-\theta^{\prime},\theta^{\prime}})&=&
W_{\psi}^{(++)}(\bm{\zeta};\widetilde{g}_{-\theta^{\prime},\theta^{\prime}})
+\mu_1\bm{i}W_{\psi}^{(++)}(\bm{\zeta};\widetilde{g}^{(1)}_{-\theta^{\prime},\theta^{\prime};H})
+\mu_2\bm{j}W_{\psi}^{(++)}(\bm{\zeta};\widetilde{g}^{(2)}_{-\theta^{\prime},\theta^{\prime};H}) \nonumber \\
&+& \mu_1\mu_2\bm{k}W_{\psi}^{(++)}(\bm{\zeta};\widetilde{g}^{(3)}_{-\theta^{\prime},\theta^{\prime};H}). 
 \label{newequation}
 \end{eqnarray}
 If we consider signals that are naturally aligned with the wavelet
$\psi_{\bm{\xi}}^{++}(\cdot),$ i.e. take $\theta^{\prime}=\theta$, then 
 a tedious but straightforward calculation yields 
\begin{eqnarray}
W_{\psi}^{(++)}(\bm{\zeta};\widetilde{g}^{(++)}_{-\theta,\theta})
 \label{++ext}
&=&2\left[1+{\mathrm{sgn}}\left(\left[\bm{r}_{-\theta}\f_b\right]_2\right)
\right]\left[1
-\bm{k}
{\mathrm{sgn}}\left(\left[\bm{r}_{-\theta}\f_b\right]_1\right)
\right]
W_{\psi}^{(ee)}(\bm{\zeta},g).
\end{eqnarray}
 This quantity is non-zero in the first quadrant of the QFT domain,
when $\theta=0.$ The hypercomplex wavelet transforms of the other 
three single quadrant supported signals follow from similar computations and are
 given by 
\begin{eqnarray}
W_{\psi}^{(++)}(\bm{\zeta};\widetilde{g}^{(-+)}_{-\theta,\theta})
&=&2\left[1+{\mathrm{sgn}}\left(\left[\bm{r}_{-\theta}\f_b\right]_2\right)
\right]\left[1
+\bm{k}
{\mathrm{sgn}}\left(\left[\bm{r}_{-\theta}\f_b\right]_1\right)
\right]
W_{\psi}^{(ee)}(\bm{\zeta},g),
\label{-+ext}
\end{eqnarray}
\begin{eqnarray}
W_{\psi}^{(++)}(\bm{\zeta};g^{(+-)}_{-\theta,\theta})
 \label{+-ext}
&=&2\left[1+{\mathrm{sgn}}\left(\left[\bm{r}_{-\theta}\f_b\right]_2\right)
\right]\left[1
-\bm{k}
{\mathrm{sgn}}\left(\left[\bm{r}_{-\theta}\f_b\right]_1\right)
\right]
W_{\psi}^{(ee)}(\bm{\zeta},g),
\end{eqnarray}
\begin{eqnarray}
W_{\psi}^{(++)}(\bm{\zeta};g^{(--)}_{-\theta,\theta})
&=&-2\left[1+{\mathrm{sgn}}\left(\left[\bm{r}_{-\theta}\f_b\right]_2\right)
\right]\left[1
+\bm{k}
{\mathrm{sgn}}\left(\left[\bm{r}_{-\theta}\f_b\right]_1\right)
\right]
W_{\psi}^{(ee)}(\bm{\zeta},g).
\label{--ext}
\end{eqnarray}
Thus we have considered the hypercomplex wavelet transform of the $\theta$-monogenic
decomposition of the real-valued signal $g(\x),$ 
%
%
%
\subsection{The Fourier Transform of the Even Wavelet Transform
of the $\theta$-Hypercomplex Signal, and Coefficients Equivalences \label{konstig2}}
One possible strategy to construct an interpretable split into amplitude
and phase descriptions would be to first construct the $\theta$-hypercomplex
extension of $\widetilde{g}_{-\theta}(\x)$ and then take the real wavelet of
this object to project the $\theta$ hypercomplex extension into a small region
of space and spatial frequency. The object of this procedure would be to
extend the aggregate of sinusoids into a suitable quaternionic extension,
and then once the quaternionic extension has been projected, the projected
components are well represented in polar form. 
It follows that:
\begin{eqnarray}
\nonumber
w^{(e e)}_{\psi}(\bm{\xi};\widetilde{g}^{(++)}_{-\theta,\theta})
&=&
w^{(e e)}_{\psi}(\bm{\xi};\widetilde{g}_{-\theta,\theta})+\bm{i}
w^{(e e)}_{\psi}(\bm{\xi};\widetilde{g}^{(1)}_{-\theta,\theta;H})+\bm{j}
w^{(e e)}_{\psi}(\bm{\xi};\widetilde{g}^{(2)}_{-\theta,\theta;H})+\bm{k}
w^{(e e)}_{\psi}(\bm{\xi};\widetilde{g}^{(3)}_{-\theta,\theta;H})\\
W^{(e e)}_{\psi}(\bm{\zeta};\widetilde{g}^{(++)}_{-\theta,\theta})
&=&
\nonumber
W^{(e e)}_{\psi}(\bm{\zeta};g)+\bm{i}
W^{(e e)}_{\psi}(\bm{\zeta};g)(-\bm{j}){\mathrm{sgn}}
\left(\left[\bm{r}_{-\theta}\f_b\right]_1\right)+\bm{j}
W^{(e e)}_{\psi}(\bm{\xi};g)(-\bm{j}){\mathrm{sgn}}
\left(\left[\bm{r}_{-\theta}\f_b\right]_2\right)\\
\nonumber
&&+\bm{k}
W^{(e e)}_{\psi}(\bm{\xi};g)(-1){\mathrm{sgn}}
\left(\left[\bm{r}_{-\theta}\f_b\right]_1\right){\mathrm{sgn}}
\left(\left[\bm{r}_{-\theta}\f_b\right]_2\right)\\
&=&\left(1+{\mathrm{sgn}}
\left(\left[\bm{r}_{-\theta}\f_b\right]_2\right)
\right)\left(1-\bm{k}{\mathrm{sgn}}
\left(\left[\bm{r}_{-\theta}\f_b\right]_1\right)
\right)W^{(e e)}_{\psi}(\bm{\zeta};g).
\label{onerep}
\end{eqnarray}
We used the expression of the Fourier transform of a $\theta$-hypercomplex
signal given by equation (\ref{thetahyperft}).
\subsection{Construction of the Local $\theta$-Hypercomplex Signal
\label{localhyp}}
Firstly let us derive an expression for $\left[W^{(ee)}_{\psi}(\bm{\zeta}_0;\widetilde{g}_{-\theta})
\right]_{\theta}^{(++)},$ to show equivalence to the other expressions.
Note that:
\begin{equation}
w_{\psi}^{(ee)}(\bm{\xi}_0,\widetilde{g}_{-\theta})=
\int_{-\infty}^{\infty}
\int_{-\infty}^{\infty} 
G(\bm{r}_{\theta}\f_b)a\Psi^{(ee)}(a\f_b)e^{2\pi\bm{j}\f_b^T \be}\;d^2\f_b.
\end{equation}
Thus if we take the $\theta$-hypercomplex extension of 
$w_{\psi}^{(ee)}(\bm{\xi}_0,\widetilde{g}_{-\theta})$ we obtain the representation:
\begin{eqnarray}
\nonumber
\left[w_{\psi}^{(ee)}(\bm{\xi}_0,\widetilde{g}_{-\theta})\right]_{\theta}^{(++)}
&=&\left[w_{\psi}^{(ee)}(\bm{\xi}_0,\widetilde{g}_{-\theta})\right]_{\theta}
+\bm{i}\left[w_{\psi}^{(ee)}(\bm{\xi}_0,\widetilde{g}_{-\theta})\right]_{\theta;H}^{(1)}\\
&&+\bm{j}\left[w_{\psi}^{(ee)}(\bm{\xi}_0,\widetilde{g}_{-\theta})\right]_{\theta;H}^{(2)}
+\bm{k}\left[w_{\psi}^{(ee)}(\bm{\xi}_0,\widetilde{g}_{-\theta})\right]_{\theta;H}^{(3)},
\end{eqnarray}
We thus using equation 
(\ref{thetahyperft}) we retrieve the following in the Fourier domain:
\begin{eqnarray}
\nonumber
{\cal F}\left\{\left[w_{\psi}^{(ee)}(\bm{\xi}_0,\widetilde{g}_{-\theta})\right]_{\theta}^{(++)}
\right\}&=&\left[1+{\mathrm{sgn}}\left((\bm{r}_{-\theta}\f_b)_2\right)
\right]\left[1-\bm{k} {\mathrm{sgn}}\left((\bm{r}_{-\theta}\f_b)_1\right)
\right]
G(\f_b)a\Psi^{(ee)}(a\bm{r}_{-\theta}\f_b)\\
&=&W_{\psi}^{(++)}(\bm{\zeta};g),
\label{ftevenrealhyper}
\end{eqnarray}
and thus we may note from equations (\ref{FThypercomplex}), (\ref{onerep}) and
finally (\ref{ftevenrealhyper}) that 
\begin{equation}
\label{equivalences}
w_{\psi}^{(++)}(\bm{\xi};g)
\equiv
w^{(e e)}_{\psi}(\bm{\xi};\widetilde{g}^{(++)}_{-\theta,\theta})
\equiv \left[w_{\psi}^{(ee)}(\bm{\xi}_0,\widetilde{g}_{-\theta})\right]_{\theta}^{(++)}.
\end{equation}
%
%
Finally as the wavelet transform is linear we may note that for any fixed value of the rotation parameter
$\theta=\theta^{\prime},$ that by equation (\ref{hypersplit}): $g(\x)=\sum_{\pm}
\widetilde{g}^{(\pm \pm)}_{-\theta^{\prime},\theta^{\prime}}(\x),$
we retrieve:
\begin{eqnarray}
\nonumber
w^{(++)}_{\psi}(\bm{\xi};g)
&=&\frac{1}{4}\left[
w^{(++)}_{\psi}(\bm{\xi};\widetilde{g}^{(++)}_{-\theta^{\prime},\theta^{\prime}})+
w^{(++)}_{\psi}(\bm{\xi};\widetilde{g}^{(-+)}_{-\theta^{\prime},\theta^{\prime}})+
w^{(++)}_{\psi}(\bm{\xi};\widetilde{g}^{(+-)}_{-\theta^{\prime},\theta^{\prime}})+
w^{(++)}_{\psi}(\bm{\xi};\widetilde{g}^{(--)}_{-\theta^{\prime},\theta^{\prime}})
\right].
\label{4thingies}
\end{eqnarray}
Hence we find that using equations (\ref{++ext}), (\ref{-+ext}),
(\ref{+-ext}) and (\ref{--ext}) that the hypercomplex wavelet transform
of $g(\x)$ viewed at $\theta^{\prime}=\theta,$ via equation (\ref{hypersplit})
the wavelet transform of $g(\x)$ is decomposed into the transform of each
component, but that naturally sum to a hypercomplex signal:
\begin{eqnarray}
W_{\psi}^{(++)}(\bm{\zeta};g)
&=&\left[1+{\mathrm{sgn}}\left(\left[\bm{r}_{-\theta}\f_b\right]_2\right)
\right]\left[1
-\bm{k}
{\mathrm{sgn}}\left(\left[\bm{r}_{-\theta}\f_b\right]_1\right)
\right]
W_{\psi}^{(ee)}(\bm{\zeta},\widetilde{g}_{\theta}).
\label{wavtrans3}
\end{eqnarray}
Thus equation (\ref{wavtrans3}) naturally agrees with equation (\ref{FThypercomplex}) but
by carrying out the calculation we may note
that we annihilate the sum of the
anti-hypercomplex components corresponding to $x_1$ having reversed polarity,
whilst the anti-hypercomplex component corresponding to $x_2$ reversing sign
is projected/converted into a hypercomplex signal: thus the three anti-hypercomplex
components are not individually annihilated like in the 1-D case. Two of the anti-components thus cancel, whilst the third is projected
into a hypercomplex object. These results will of course be necessary to
derive the properties of the hypercomplex wavelet transform of the phase-shifted signal.
%

\subsection{$\theta$-Hypercomplex Extension of Phase-Shifted Signal \label{2dphasehyper}}
Note that the $\theta$-hypercomplex extension of the phase shifted signal has the FT:
\begin{eqnarray*}
{\cal F}\left\{\left(
\Lambda^{2D}_{\bm{\theta}_s,\theta}g\right)_{\theta}^{(++)}\right\}&=&
\left(1+{\mathrm{sgn}}\left(\left[
\bm{r}_{-\theta}\f\right]_{2}\right)\right)\left(1-\bm{k}{\mathrm{sgn}}\left(\left[
\bm{r}_{-\theta}\f\right]_{1}\right)\right){\cal F}\left\{\left(
\Lambda^{2D}_{\bm{\theta}_s,\theta}g\right)\right\}(\bm{r}_{-\theta}\f)\\
{\cal F}\left\{
\widetilde{\left(\Lambda^{2D}_{\bm{\theta}_s,\theta}g\right)}_{-\theta,\theta}^{(++)}\right\}&=&
\left(1+{\mathrm{sgn}}\left(\left[
\bm{r}_{-\theta}\f\right]_{2}\right)\right)\left(1-\bm{k}{\mathrm{sgn}}\left(\left[
\bm{r}_{-\theta}\f\right]_{1}\right)\right){\cal F}\left\{\left(
\Lambda^{2D}_{\bm{\theta}_s,\theta}g\right)\right\}(\f).
\end{eqnarray*}
The FT of the separable phase-shifted signal is
\begin{eqnarray*}
{\cal F}\left\{
\Lambda^{2D}_{\bm{\theta}_s,\theta}g\right\}(\f)&=&\int_{-\infty}^{\infty} 
\int_{-\infty}^{\infty}  \left|
\widetilde{g}_{-\theta}^{(++)}(\bm{r}_{-\theta}\x)\right|
\cos\left(2\pi\widetilde{\alpha}_{-\theta}(\bm{r}_{-\theta}\x)-
\theta_{s,1} \right)\cos\left(2\pi\widetilde{\beta}_{-\theta}(\bm{r}_{-\theta}\x)-
\theta_{s,2} \right)\\
&&e^{-2\pi \bm{j}(\bm{r}_{-\theta}\x)^T \bm{r}_{-\theta}\f}\;d^2\x\\
&=&\frac{1}{4}
\left[e^{-\bm{j}\theta_{s,1}}\widetilde{G}_{1,-\theta}^{(+)}(\left[\bm{r}_{-\theta}\f
\right]_1)+e^{\bm{j}\theta_{s,1}}\widetilde{G}_{1,-\theta}^{(-)}(\left[\bm{r}_{-\theta}\f
\right]_1)\right]\\
&&\left[e^{-\bm{j}\theta_{s,2}}\widetilde{G}_{2,-\theta}^{(+)}(\left[\bm{r}_{-\theta}\f
\right]_2)+e^{\bm{j}\theta_{s,2}}\widetilde{G}_{2,-\theta}^{(-)}(\left[\bm{r}_{-\theta}\f
\right]_2)\right]\\
&=&\frac{1}{4}
\left[e^{-\bm{j}\theta_{s,1}}(1+{\mathrm{sgn}}\left(\left[
\bm{r}_{-\theta}\f\right]_{1}\right))
+e^{\bm{j}\theta_{s,1}}(1-{\mathrm{sgn}}\left(\left[
\bm{r}_{-\theta}\f\right]_{1}\right))\right]\widetilde{G}_{-\theta}(\bm{r}_{-\theta}\f)
\\
&&\left[e^{-\bm{j}\theta_{s,2}}(1+{\mathrm{sgn}}\left(\left[
\bm{r}_{-\theta}\f\right]_{2}\right))
+e^{\bm{j}\theta_{s,2}}(1-{\mathrm{sgn}}\left(\left[
\bm{r}_{-\theta}\f\right]_{2}\right))
\right].
\end{eqnarray*}
After some algebra we then obtain:
\begin{eqnarray*}
{\cal F}\left\{\widetilde{\left(
\Lambda^{2D}_{\bm{\theta}_s,\theta}g\right)}_{-\theta,\theta}^{(++)}\right\}(\f)&=&
\left(1+{\mathrm{sgn}}\left(\left[
\bm{r}_{-\theta}\f\right]_{2}\right)\right)\left(1-\bm{k}{\mathrm{sgn}}\left(\left[
\bm{r}_{-\theta}\f\right]_{1}\right)\right)\\
&&\frac{1}{4}
\left[e^{-\bm{j}\theta_{s,1}}(1+{\mathrm{sgn}}\left(\left[
\bm{r}_{-\theta}\f\right]_{1}\right))
+e^{\bm{j}\theta_{s,1}}(1-{\mathrm{sgn}}\left(\left[
\bm{r}_{-\theta}\f\right]_{1}\right))\right]\widetilde{G}_{-\theta}(\bm{r}_{-\theta}\f)
\\
&&\left[e^{-\bm{j}\theta_{s,2}}(1+{\mathrm{sgn}}\left(\left[
\bm{r}_{-\theta}\f\right]_{2}\right))
+e^{\bm{j}\theta_{s,2}}(1-{\mathrm{sgn}}\left(\left[
\bm{r}_{-\theta}\f\right]_{2}\right))
\right]\\
 &=&e^{-\bm{i}\theta_{s,1}}\widetilde{G}_{-\theta,\theta}^{(++)}(\f)e^{-\bm{j}\theta_{s,2}}.
 \end{eqnarray*}
Hence we may deduce that, as would be expected for $g(\x)$ a separable signal
we determine that:
\begin{equation}
\label{phasenew}
\widetilde{\left(\Lambda^{2D}_{\bm{\theta}_s,\theta}g\right)}_{-\theta,\theta}^{(++)}(\x)
=e^{-\bm{i}\theta_{s,1}}
\widetilde{g}_{-\theta,\theta}^{(++)}(\x)e^{-\bm{j}\theta_{s,2}}.\end{equation}

\subsection{Hypercomplex CWT of a Phase-Shifted Signal \label{phasehyperinv}}
We note that the hypercomplex CWT of the observed signal $g(\x)$ is by equations
(\ref{equivalences}) and (\ref{phasenew}):
\begin{eqnarray*}
\nonumber
w_{\psi}^{(++)}(\bm{\xi};\Lambda^{2D}_{\bm{\theta}_s,\theta}g)
&=&w_{\psi}^{(ee)}(\bm{\xi};\widetilde{\left(
\Lambda^{2D}_{\bm{\theta}_s,\theta}g\right)}_{-\theta,\theta}^{(++)})\\
&=&e^{-\bm{i}\theta_{s,1}}
w_{\psi}^{(ee)}(\bm{\xi};\widetilde{g}_{-\theta,\theta}^{(++)}(\x))
e^{-\bm{j}\theta_{s,2}}\\
&=&e^{-\bm{i}\theta_{s,1}}
w_{\psi}^{(++)}(\bm{\xi};g)e^{-\bm{j}\theta_{s,2}},
\end{eqnarray*}
from equation (\ref{phasenew}). Thus the result follows.

\subsection{FT and QFT of the Monogenic Wavelet Transform \label{monowav}}
Recalling equations (\ref{FTRiesz}), (\ref{FTCWT22}) as well as (\ref{QFTCWT22})
 we obtain the FT  
\begin{eqnarray}
W_{\psi}^{(+)}(\bm{\zeta};g)&=&  
W_{\psi}^{(r)}(\bm{\zeta};g)-\bm{i}
W_{\psi}^{(1)}(\bm{\zeta};g)-\bm{j}
W_{\psi}^{(2)}(\bm{\zeta};g)
 \nonumber \\
&=&\left[1 -\bm{k}\cos{(\phi_b - \theta)}+\sin{(\phi_b - \theta)}\right]
W_{\psi}^{(r)}(\bm{\zeta};g).
\end{eqnarray}
The corresponding QFT is
\begin{eqnarray}
W_{\psi;Q}^{(+)}(\bm{\kappa};g)&=& W_{\psi;Q}^{(r)}(\bm{\kappa};g)-\bm{i}
W_{\psi;Q}^{(1)}(\bm{\kappa};g)-
W_{\psi;Q}^{(2)}(\bm{\kappa};g)\bm{j} \nonumber \\
&=& \frac{1-\bm{k}}{2} \left[ 1+ \cos{(\nu_b - \theta)} + \sin{(\nu_b-\theta)} \right] 
G(\q_b) a\Psi^{(r)*}\left(a\bm{r}_{-\theta} \q_b \right) \nonumber \\
&+& \frac{1+\bm{k}}{2} \left[ 1+ \cos{(\nu_b + \theta)} + \sin{(\nu_b + \theta)} \right] 
G(\bm{J}_{\pi/2} \q_b) a\Psi^{(r)*}\left( a \bm{r}_{-\theta} \bm{J}_{\pi/2}\q_b \right) \nonumber \\
&=& W_{\psi;Q}^{(r)}(\bm{\kappa};g) + \cos{\theta} \left[ \cos{\nu_b} + \sin{\nu_b} \right] W_{\psi;Q}^{(r)}(\bm{\kappa};g)
 + \sin{\theta} \left[ \cos{\nu_b} - \sin{\nu_b} \right] \bm{i} W_{\psi;Q}^{(r)}(\bm{\kappa};g) \bm{j}. \nonumber \\
\end{eqnarray}

\subsection{Wavelet transform of $\theta$-Monogenic \& Anti-Monogenic Decomposition Components \label{splitting}}
We find that the Fourier transform of $\widetilde{g}_{-\theta^{\prime},\theta^{\prime}}^{(+)}(\x)$
is given by:
\begin{eqnarray}
&& W_{\psi}^{(+)}(\bm{\zeta};
\widetilde{g}_{-\theta^{\prime},\theta^{\prime}}^{(+)}) = 
W_{\psi}^{(+)}(\bm{\zeta};g) +  \left[  \bm{i} \; 
W_{\psi}^{(+)}(\bm{\zeta};\widetilde{g}^{(1)}_{-\theta^{\prime},\theta^{\prime},R}) +   \bm{j} \;
W_{\psi}^{(+)}(\bm{\zeta};\widetilde{g}^{(2)}_{-\theta^{\prime},\theta^{\prime},R}) \right] \nonumber \\
&&=\left[1+\sin(\phi_b-\theta)+\cos(\phi_b-\theta)\cos(\phi_b-
\theta^{\prime})+\sin(\phi_b-
\theta^{\prime})+
\sin(\phi_b-\theta)\sin(\phi_b-
\theta^{\prime})
\right.\nonumber \\
&&+\left.\bm{k}\left(-\cos(\phi_b-\theta)-\cos(\phi_b-
\theta^{\prime})-\sin(\phi_b-\theta)\cos(\phi_b-
\theta^{\prime})
+\cos(\phi_b-\theta)\sin(\phi_b-
\theta^{\prime})
\right)
\right]W^{(r)}_{\psi}(\bm{\zeta};g). \nonumber \\
\label{complicatedexpression}
\end{eqnarray}
We remark that with $\theta=
\theta^{\prime}$ this then becomes
\begin{eqnarray}
\nonumber
 W_{\psi}^{(+)}(\bm{\zeta};\widetilde{g}_{-\theta,\theta}^{(+)})
&=&   \left\{  2+2\sin(\phi_b-\theta)-2\bm{k}\cos(\phi_b-\theta)  \right\}  \nonumber \ W^{(r)}_{\psi}(\bm{\zeta};g) \nonumber \\ 
&=&\left\{  2+2\sin(\phi_b-\theta)-2\bm{k}\cos(\phi_b-\theta)  \right\}
G(\f_b)a\Psi^{(r)*}(a \bm{r}_{-\theta}\f_b)\nonumber\\
&=&\left\{  2+2\sin(\phi_b-\theta)-2\bm{k}\cos(\phi_b-\theta)  \right\}
\widetilde{G}_{-\theta,\theta}(\f_b)a\Psi^{(r)*}(a \bm{r}_{-\theta}\f_b)\nonumber\\
& =&2\widetilde{G}^{(+)}_{-\theta,\theta}(\f_b)
a\Psi^{(r)*}(a \bm{r}_{-\theta}\f_b)=2
W_{\psi}^{(r)}(\bm{\zeta};\widetilde{g}_{-\theta,\theta}^{(+)}).
\label{FTCWTmnkl3p}
\end{eqnarray}
Furthermore {\em mutatis mutandis} with $\widetilde{g}_{-\theta,\theta}^{(-)}(\x)$
replacing $\widetilde{g}_{-\theta,\theta}^{(+)}(\x)$
we derive that
\begin{equation}
\label{FTCWTmnkl3p-}
W_{\psi}^{(+)}(\bm{\zeta};\widetilde{g}_{-\theta,\theta}^{(-)})=0.\end{equation}
This completes the proof. We can note that equation (
\ref{complicatedexpression}) establishes
the necessity of introducing the $\theta$-monogenic signal. The Monogenic
wavelet transform will not annihilate the monogenic signal for example, and
in the subsequent section this may produce problems.

\subsection{Construction of Local $\theta$-Monogenic Signal 
\label{ftthetamon}}
We note the form of the Fourier transform of $w_{\psi}^{(+)}(\bm{\xi};g)$
from equation (\ref{FTCWTaw}) and then 
from equation (\ref{weirdeqn}) we find that (recalling $\bm{\xi}_0=
\left[a,0,\be\right]^T$):
\begin{eqnarray*}
{\cal F}\left\{\left[w^{(r)}_{\psi}(\bm{\xi}_0;\widetilde{g}_{-\theta})\right]^{(+)}_{\theta}
\right\}&=&\left[1+\sin(\phi_b-\theta)-\bm{k}\cos(\phi_b-\theta)
\right]\widetilde{G}_{-\theta}\left(\bm{r}_{-\theta}\f_b\right)
a\Psi^{(r)*}\left(a\bm{r}_{-\theta}\f_b\right)
\\
&=&\left[1+\sin(\phi_b-\theta)-\bm{k}\cos(\phi_b-\theta)
\right]G(\f_b)a\Psi^{(r)*}\left(a\bm{r}_{-\theta}\f_b\right)\\
&=&W_{\psi}^{(+)}(\bm{\zeta};g)\\
\left[w^{(r)}_{\psi}(\bm{\xi}_0;\widetilde{g}_{-\theta})\right]^{(+)}_{\theta}
&=&w_{\psi}^{+}(\bm{\xi};g).
\end{eqnarray*}
Thus the monogenic wavelet transform coefficients correspond to the monogenic
extension of the real wavelet transform in the rotated frame of reference.
For the second equality
we start by using equation (\ref{decomposition}) and take
\[g(\x)=\frac{1}{2}\left(\widetilde{g}_{-\theta^{\prime},\theta^{\prime}}^{(+)}(\x)+
\widetilde{g}_{-\theta^{\prime},\theta^{\prime}}^{(-)}(\x)\right).\]
By the linearity of the wavelet transform we have that with $\theta^{\prime}=\theta$
using equations (\ref{FTCWTmnkl3p}) and (\ref{FTCWTmnkl3p-}):
\begin{eqnarray}
\nonumber
W_{\psi}^{(+)}(\bm{\zeta};g)&=&\frac{1}{2}\left[
W_{\psi}^{(+)}(\bm{\zeta};\widetilde{g}_{-\theta,\theta}^{(+)})+
W_{\psi}^{(+)}(\bm{\zeta};\widetilde{g}_{-\theta,\theta}^{(-)})
\right]\\
&=&\frac{1}{2}\left[2 W_{\psi}^{(r)}(\bm{\zeta};\widetilde{g}_{-\theta,\theta}^{(+)})
+0\right]= W_{\psi}^{(r)}(\bm{\zeta};\widetilde{g}_{-\theta,\theta}^{(+)}),
\end{eqnarray}
and thus the result follows. Hence the operation of constructing the monogenic
wavelet transform
can either be considered in light of the operation of 1) taking the monogenic wavelet
transform of a real signal, 2) finding the real wavelet transform of a $\theta$-monogenic
signal or as finding the $\theta$-monogenic image of the wavelet transform
in $\be$ of the rotated by $\theta$ signal.
%
\subsection{Local Invariance under Rotation \label{Locinv}}
If the monogenic wavelet is constructed from an isotropic wavelet,
then it follows that:
\begin{eqnarray}
w_{\psi}^{(e)}(\bm{\xi};g)
&=&\int_{-\infty}^{\infty}
\int_{-\infty}^{\infty} W_{\psi}^{(e)}(\bm{\zeta};g) e^{2\bm{j}\pi
\f_b^T \be}\;d^2\f_b=\int_{-\infty}^{\infty}
\int_{-\infty}^{\infty}G(\f_b)a \Psi\left(a\bm{r}_{-\theta}\f_b\right)e^{2\bm{j}\pi
\f_b^T \be}\;d^2\f_b
\nonumber \\ 
&\overset{(1)}{=}&
\int_{-\infty}^{\infty}
\int_{-\infty}^{\infty}G(\f_b)a \Psi\left(a\f_b\right)e^{2\bm{j}\pi
\f_b^T \be}\;d^2\f_b=\int_{-\infty}^{\infty}\int_{-\infty}^{\infty} 
W_{\psi}^{(e)}(\bm{\zeta}_0;g) e^{2\bm{j}\pi
\f_b^T \be}\;d^2\f_b
\nonumber \\ 
&=&w_{\psi}^{(e)}(\bm{\xi}_0;g),
\label{reltheta} 
\end{eqnarray}
where we have used that the real mother wavelet is isotropic in (1). In a similar fashion, we calculate the for the 
 CWTs with respect to the first and second Riesz components of the monogenic wavelet, 
\begin{eqnarray}
w_{\psi}^{(1)}(\bm{\xi};g)
 \label{firstreltheta}
&=&\ \  \cos(\theta)w_{\psi}^{(1)}(\bm{\xi}_0;g)+
\sin(\theta)w_{\psi}^{(2)}(\bm{\xi}_0;g), \nonumber \\
w_{\psi}^{(2)}(\bm{\xi};g)
&=&-\sin(\theta)w_{\psi}^{(1)}(\bm{\xi}_0;g)+
\cos(\theta)w_{\psi}^{(2)}(\bm{\xi}_0;g).
\end{eqnarray} 
 Using the above expressions (\ref{reltheta}), and (\ref{firstreltheta}), trivially follows that 
\begin{eqnarray}
 \left|w_{\psi}^{(+)}(\bm{\xi};g)\right|^2 &=&
w_{\psi}^{(r)2}(\bm{\xi};g)+w_{\psi}^{(1)2}(\bm{\xi};g)+w_{\psi}^{(2)2}(\bm{\xi};g)
= \left|w_{\psi}^{(+)}(\bm{\xi}_0;g)\right|^2.
\label{rotinvarianceproof}
\end{eqnarray}
Hence, if an isotropic real mother wavelet is used, the magnitude of the monogenic
wavelet coefficients is invariant to to the choice of orientation.

\subsection{Directional Selectivity \label{Dirsel}}
For simplicity we model the signal is directional in the Fourier domain,
i.e. via equation (\ref{dirsignaldef}). Note that the wavelet transform corresponds
to patterns at period $f.$ It then transpires that:
\begin{eqnarray*}
w_{\psi}^{(1)}\left(\bm{\xi};g \right)
&=&\int_{-\infty}^{\infty} \int_{-\infty}^{\infty}
W_{\psi}^{(1)}(\bm{\zeta};g) e^{2\bm{j}\pi\f_{b}^T\be}\;d^2\f_{b}\\
&\overset{(1)}{=}&\int_{0}^{\infty} \int_{0}^{2\pi} (-\bm{j})
f_b\cos(\phi_b-\theta)
\frac{\widetilde{G}(f_b)}{2}
\left[\delta(\phi_b-\phi_0-\pi)+\delta(\phi_b-\phi_0)\right]a\Psi^{(e)}
(af_b)\\
\nonumber
&&e^{2\bm{j}\pi \f_b^T\be}\;df_b\;d\phi\\
&=&\cos(\theta-\phi_0) \int_0^{\infty} f_b\widetilde{G}(f_b)
a\Psi^{(e)}
(af)\sin\left(2\pi f_b b\cos(\chi_b-\phi_0)\right)\;df_b,
\end{eqnarray*}
In a similar fashion,
\begin{eqnarray}
w_{\psi}^{(2)}\left(\bm{\xi};g \right)
&=&\sin(\theta-\phi_0) \int_0^{\infty} f_b\widetilde{G}(f_b)
a\Psi^{(e)}
(af_b)\sin\left(2\pi f_b b\cos(\chi_b-\phi_0)\right)\;df_b.
\label{sin1}
\end{eqnarray}
%
For (1) to hold, $\psi^{(r)}(\x)\equiv \psi^{(e)}(\x),$ {\em i.e.} is isotropic.
For a directional signal, when $w_{\psi}^{(1)}\left(\bm{\xi};g \right) \neq
0,$
\begin{equation}
\label{detlocalrot}
\theta-\phi_0=\tan^{-1}\left(\frac{w_{\psi}^{(2)}\left(\bm{\xi};g \right)}{
w_{\psi}^{(1)}\left(\bm{\xi};g \right)} \right).
\end{equation}
Thus the local directionality of a directional signal, can at any rotational
angle of the wavelet transform, be determined from the isotropic monogenic
wavelet transform.
Furthermore it is clear that
\begin{equation}
\label{detlocalrot2}
w_{\psi}^{(2)}\left(\bm{\xi}^{\dagger};g \right) =0,
\end{equation}
and hence by equation (\ref{RFTCWT}), the stated result follows.
%

\subsection{CWT of Phase-Shifted Signal \label{CWTINV}}
For brevity we denote the wavelet transform of the phase-shifted signals
by
$
w_{\psi,\theta_s}^{(\pm)}(\bm{\xi};g)=
w_{\psi}^{(\pm)}(\bm{\xi};\Lambda_{\theta_{s}} g).
$
We note that
\begin{eqnarray*}
\Lambda_{\theta_{s}} g(\x)&=&\left|g^{(+)}(\x)\right| \ \cos{( 2\pi\phi(\x)
- \theta_{s} )}\\
&=&\frac{1}{2}
\left[e^{  -\bm{e}_{\widetilde{\nu}_{-\theta,\theta}}(\x) \theta_{s} }\  \widetilde{g}^{(+)}_{-\theta,\theta}(\x)+
e^{  \bm{e}_{\widetilde{\nu}_{-\theta,\theta}}(\x) \theta_{s} }\  \widetilde{g}^{(-)}_{-\theta,\theta}(\x)
\right].
\end{eqnarray*}
Hence with the additional assumption that $\bm{e}_{\widetilde{\nu}_{-\theta,\theta}}(\x)=
\bm{e}_{\widetilde{\nu}_{-\theta,\theta}}(\be)$
is constant over the spatial width the wavelet transform is averaging over,
i.e. $\left|\x-\be\right|<a x_{i}$ we retrieve that 
\begin{eqnarray}
\nonumber
w_{\psi,\theta_s}^{(+)}(\bm{\xi};g)&=&
\frac{1}{2}\left[w_{\psi}^{(+)}(\bm{\xi};e^{-\bm{e}_{\widetilde{\nu}_{-\theta,\theta}}
(\x) \theta_{s} }\  \widetilde{g}^{(+)}_{-\theta,\theta}(\x))+
w_{\psi}^{(+)}(\bm{\xi};e^{\bm{e}_{\widetilde{\nu}_{-\theta,\theta}}(\x) 
\theta_{s} }\  \widetilde{g}^{(-)}_{-\theta,\theta}(\x))
\right]\\
\nonumber
&=&
\frac{1}{2}\left[w_{\psi}^{(+)}(\bm{\xi};e^{-\bm{e}_{\widetilde{\nu}_{-\theta,\theta}}(\be) \theta_{s} }\  \widetilde{g}^{(+)}_{-\theta,\theta}(\x))+w_{\psi}^{(+)}(\bm{\xi};
e^{\bm{e}_{\widetilde{\nu}_{-\theta,\theta}}(\be) \theta_{s} }\  \widetilde{g}^{(-)}_{-\theta,\theta}(\x))
\right]\\
\nonumber
&=&\frac{1}{2}\left[
e^{ -\bm{e}_{\widetilde{\nu}_{-\theta,\theta}}(\be)\theta_{s}}w_{\psi}^{(+)}(\bm{\xi};\  \widetilde{g}^{(+)}_{-\theta,\theta}(\x))+
e^{\bm{e}_{\widetilde{\nu}_{-\theta,\theta}}
(\be)\theta_{s}}w_{\psi}^{(+)}(\bm{\xi};  \widetilde{g}^{(-)}_{-\theta,\theta}(\x))
\right]\\
&=&\frac{1}{2}
e^{ -\bm{e}_{\widetilde{\nu}_{-\theta,\theta}}(\be)\theta_{s}}w_{\psi}^{(+)}(\bm{\xi};\
\widetilde{g}^{(+)}_{-\theta,\theta}(\x))=
e^{ -\bm{e}_{\widetilde{\nu}_{-\theta,\theta}}(\be)\theta_{s}}w_{\psi}^{(+)}(\bm{\xi};g(\x)),
\end{eqnarray}
where we used equation (\ref{FTCWTmnkl3p-}).
Also from equation (\ref{FTCWTmnkl3p}) we may note
\begin{eqnarray}
w_{\psi}^{(+)}(\bm{\xi};g(\x))&=&
\frac{1}{2}
w_{\psi}^{(+)}(\bm{\xi};\  \widetilde{g}^{(+)}_{\theta}(\x)),
\end{eqnarray}
and thus the result follows.

\end{document}